\documentclass[a4paper,10.5pt]{article}

\usepackage{stmaryrd}
\usepackage{bm}
\usepackage{cite}
\usepackage{color}
\usepackage{dsfont}
\usepackage{nicefrac}
\usepackage{multirow}
\usepackage{graphics}
\usepackage{graphicx}
\usepackage{epsfig}
\usepackage{subfigure}
\usepackage{pgf}
\usepackage{amsmath,amssymb,amsfonts}
\usepackage{wasysym}
\usepackage{yhmath}
\usepackage{booktabs}

\begin{document}

\title{Fiber bundle topology optimization for surface flows}

\author{Yongbo Deng$^1$\footnote{dengyb@ciomp.ac.cn (Y. Deng)}, Weihong Zhang$^2$\footnote{zhangwh@nwpu.edu.cn (W. Zhang)}, Jihong Zhu$^2$, \\
Yingjie Xu$^2$, Zhenyu Liu$^3$, Jan G. Korvink$^4$ \\
1 State Key Laboratory of Applied Optics, \\
Changchun Institute of Optics, Fine Mechanics and Physics, \\
Chinese Academy of Sciences, 130033 Changchun, China; \\
2 State IJR Center of Aerospace Design and Additive Manufacturing, \\
Northwestern Polytechnical University, 710072 Xi'an, China; \\
3 Changchun Institute of Optics, Fine Mechanics and Physics, \\
Chinese Academy of Sciences, 130033 Changchun, China; \\
4 Institute of Microstructure Technology (IMT), \\
Karlsruhe Institute of Technology (KIT), \\
Hermann-von-Helmholtzplatz 1, \\
76344 Eggenstein-Leopoldshafen, Germany.}

\maketitle

\abstract{This paper presents a topology optimization approach for the surface flows on variable design domains. Via this approach, the matching between the pattern of a surface flow and the 2-manifold used to define the pattern can be optimized, where the 2-manifold is implicitly defined on another fixed 2-manifold named as the base manifold. The fiber bundle topology optimization approach is developed based on the description of the topological structure of the surface flow by using the differential geometry concept of the fiber bundle. The material distribution method is used to achieve the evolution of the pattern of the surface flow. The evolution of the implicit 2-manifold is realized via a homeomorphous map. The design variable of the pattern of the surface flow and that of the implicit 2-manifold are regularized by two sequentially implemented surface-PDE filters. The two surface-PDE filters are coupled, because they are defined on the implicit 2-manifold and base manifold, respectively. The surface Navier-Stokes equations, defined on the implicit 2-manifold, are used to describe the surface flow. The fiber bundle topology optimization problem is analyzed using the continuous adjoint method implemented on the first-order Sobolev space. Several numerical examples have been provided to demonstrate this approach, where the combination of the viscous dissipation and pressure drop is used as the design objective.

}

\section{Introduction} \label{sec:IntroductionManifold}

Surface flows can greatly decrease the computational cost in the numerical design of the related fluidic structures.
The fluid flows in the channels attached on the walls of equipments can be described as surface flows on the curved surfaces corresponding to the outer shapes of fluidic structures. The streamsurfaces corresponding to the outer shapes of fluidic structures with complete-slip boundaries can be described as surface flows separated from the bulk flows, where the complete-slip boundaries can be approximated and achieved by chemically coating or physically structuring solid surfaces to derive the extreme hydrophobicity \cite{KwonLangmuir2009}, using the optimal control method to manipulate the boundary velocity of flows \cite{SritharanSIAM1998}, and producing vapor layers between the solid and liquid phases based on the Leidenfrost phenomenon \cite{ThimblebyPhysicsEducation1989}, etc.

\begin{figure}[!htbp]
\centering
\includegraphics[width=1\textwidth]{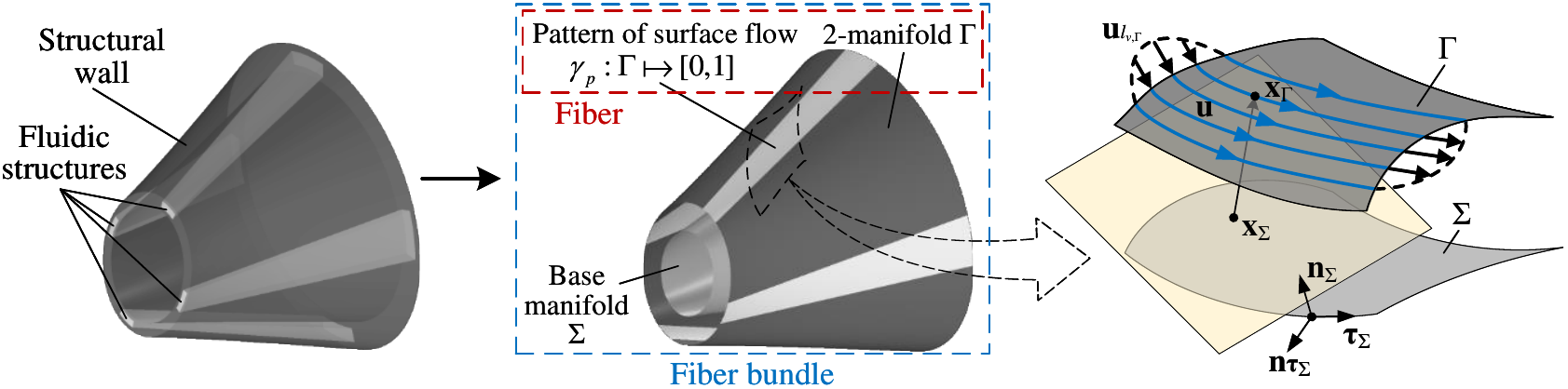}
\caption{Sketch for the fiber bundle of a surface flow, where $\Sigma$ is the base manifold, $\Gamma$ is the implicit 2-manifold used to define the pattern of the surface flow, $\gamma_p: \Gamma \mapsto \left[ 0,1 \right]$ is the pattern of the surface flow, $\mathbf{u}$ is the fluid velocity of the surface flow, $\mathbf{u}_{l_{v,\Gamma}}$ is the known fluid velocity at the boundary of $\Gamma$, $\mathbf{n}_\Gamma$ is the unitary normal vector of $\Gamma$, $\boldsymbol\tau_\Gamma$ is the unitary tangential vector at $\partial\Gamma$, $\mathbf{n}_{\boldsymbol\tau_\Gamma} = \mathbf{n}_\Gamma\times\boldsymbol\tau_\Gamma $ is the outward unitary normal at $\partial\Gamma$, $\mathbf{x}_\Sigma$ denotes a point on $\Sigma$, and $\mathbf{x}_\Gamma$ denotes a point on $\Gamma$. This paper focuses on the laminar surface flows with low and moderate Reynolds numbers to demonstrate the fiber bundle topology optimization approach, although the sketched surface flows can be turbulent with high Reynolds number.}\label{fig:ApplicationSorting}
\end{figure}

The topological structure of a surface flow can be described as the fiber bundle demonstrated in Fig. \ref{fig:ApplicationSorting}. Fiber bundle is a concept of differential geometry \cite{ChernDifferentialGeometry1999}. It is composed of the base manifold and the fiber defined on it, where the manifold represents the topological space locally homeomorphous to an Euclidean space. For the surface flow, the flow pattern together with its definition domain corresponds to the fiber of the fiber bundle. If there exists a 2-manifold homeomorphous to the fiber, it can be set as the base manifold of the fiber bundle. In computation, the existence of the base manifold can be ensured by presetting a fixed geometrical surface as the base manifold, then the fiber can be found on the preset base manifold. That means the definition domain of the pattern is an implicit 2-manifold defined on the preset base manifold, where the implicit 2-manifold is the fluid/solid interface corresponding to the outer-shape surface. The reason for this paper to use the concept of fiber bundle is to describe the topological structure of the surface flow as an ensemble instead of three separated components. Therefore, the task for the fiber bundle topology optimization of the surface flow is to find the optimized matching between the pattern and the implicit 2-manifold defined on the preset base manifold.

Topology optimization is currently regarded to be one of the most robust methodology for the determination of material distribution in structures that meet given structural performance criteria \cite{ChengOlhoffIJSS1981,Bendsoe1988,BendsoeAndSigmund2003}. With regard to flow problems, topology optimization has been implemented for Stokes flows \cite{StevenLiXie2000,Borrvall2003}, creeping fluid flows \cite{GuestProevost2006}, steady Navier-Stokes flows \cite{Gersborg-Hansen2006}, unsteady Navier-Stokes flows \cite{Kreissl2011,Deng2011}, flows with body forces \cite{DengCMAME2013,DengSMO2013}, turbulent flows \cite{DilgenCMAME2018,YoonCMAME2016}, two-phase flows of immiscible fluids \cite{DengLiuWuCiCP2017}, electroosmotic flows \cite{GregersenOkkels2009,DengKorvinkIJHMT2018} and flows of non-Newtonian fluids \cite{PingenCMA2010,AlonsoSMO2020}, etc; topology optimization for flow problems have been reviewed in \cite{AlexandersenFluids2020}. With regard to interfacial patterns, researches have been implemented for stiffness and multi-material structures \cite{VermaakSMO2014,SigmundJMPS1997,GaoZhangIJNME2011,
LuoAiaaJ2012,WangCMAME2004,ZhouWangSMO2007,PanagiotisVogiatzis2018}, layouts of shell structures \cite{KrogOlhoffComputersStructures1996,AnsolaComputersStructures2002,HassaniSMO2013, LochnerAldingerSchumacher2014,ClausenActaMechanicaSinica2017,DienemannStructMultidiscOptim2017, YanChengWangJSV2018}, electrode patterns of electroosmosis \cite{DengKorvinkIJHMT2018}, fluid-structure and fluid-particle interaction \cite{YoonFSIIJNME2010,LundgaardFSISMO2018,AndreasenFSISMO2019}, energy absorption \cite{AuligUlm2012}, cohesion \cite{MauteIJNME2017}, actuation \cite{MauteSMO2005} and wettability control \cite{DengMicrotextureCMAME2018,DengMicrotextureAMM2019,DengMicrotextureSMO2020}, etc.; topology optimization approaches implemented on 2-manifolds have also been developed with applications in elasticity, wettability control, heat transfer and electromagnetics \cite{HuoJAM2022,ZhangFengSMO2022,DengCMAME2020}; and the fiber bundle topology optimization approach has been developed for wettability control at fluid/solid interfaces \cite{DengMicrotextureSMO2020}; recently, topology optimization of surface flows has extended the design space of fluidic structures onto the 2-manifolds \cite{DengJCPTOOPSurfaceFlow2022}.

It is natural for one to ask if it is possible to implement the topology optimization to match the pattern of a surface flow and the implicit 2-manifold on which the pattern is defined. If such topology optimization can be achieved, the design space and design freedom will be further extended for flow problems by including the design domain for the pattern of the surface flow into the design space, where the design domain is the implicit 2-manifold. Therefore, this paper presents the fiber bundle topology optimization approach for the surface flow.

For the fiber bundle topology optimization approach, the material distribution method pioneered by \cite{Bendsoe1988} is used to determine the pattern of the surface flow. The implicit 2-manifold used to define the surface flow is described on the base manifold. Then, two sets of design variables are required for the pattern of the surface flow and the implicit 2-manifold, respectively. For the material distribution method, a porous medium model has been developed for Stokes flows \cite{Borrvall2003}. This model was then extended to implement topology optimization for steady and unsteady Navier-Stokes flows \cite{Gersborg-Hansen2006,Kreissl2011,Deng2011}. In this model, the porous medium was filled in the two/three-dimensional design domains. Correspondingly, an artificial Darcy friction was introduced into the force terms of the Stokes equations and Navier-Stokes equations. The impermeability of the porous medium was evolved in the topology optimization procedure to derive the fluidic structures. Inspired by the porous medium model, topology optimization for surface flows has been implemented by filling the porous medium onto fixed 2-manifolds, where an artificial Darcy friction is added to the surface Navier-Stokes equations \cite{DengJCPTOOPSurfaceFlow2022}. This paper inherits this model with the porous medium to implement the fiber bundle topology optimization for surface flows.

The remained sections of this paper are organized as follows.
In Section \ref{sec:MethodologyManifold}, a monolithic description of the fiber bundle topology optimization problem for a surface flow is presented.
In Section \ref{sec:NumericalImplementationSurfaceNSEqus}, numerical implementation for the iterative solution of the fiber bundle topology optimization problem is introduced.
In Section \ref{sec:NumericalExamplesMatchinOptSurfaceFlow}, numerical tests are provided to demonstrate the developed fiber bundle topology optimization approach.
In Sections \ref{sec:Conclusions} and \ref{sec:Acknowledgements}, the conclusion and acknowledgment of this paper are provided. In Section \ref{sec:Appendix}, details are provided for the adjoint analysis of the fiber bundle topology optimization problem.
All the mathematical descriptions are implemented in a Cartesian system.

\section{Methodology}\label{sec:MethodologyManifold}

In this section, the fiber bundle topology optimization problem is described to match the pattern of the surface flow and the implicit 2-manifold on which the surface flow is defined. The implicit 2-manifold is defined on the base manifold. The incompressible surface fluid is considered.

\subsection{Physical model and material interpolation} \label{sec:PorousMediumModelSurfaceNSEqus}

In the fiber bundle topology optimization for the surface flow, the porous medium model is utilized. In this model, the porous medium is filled onto the implicit 2-manifold. Correspondingly, the artificial Darcy friction is added to the surface Navier-Stokes equations. The artificial Darcy friction is derived based on the constitutive law of the porous medium. It is assumed to be proportional to the fluid velocity \cite{Borrvall2003,Gersborg-Hansen2006}:
\begin{equation}\label{equ:ArtificialDarcyFriction}
  \mathbf{b}_\Gamma = - \alpha \mathbf{u},~\forall \mathbf{x}_\Gamma \in \Gamma
\end{equation}
where $\alpha$ is the impermeability; $\Gamma$ is the implicit 2-manifold; and $\mathbf{x}_\Gamma$ denotes a point on $\Gamma$. When the porosity of the porous medium is zero, it corresponds to a solid material with infinite impermeability and zero fluid velocity caused by the infinite friction force. When the porosity is infinite, it corresponds to the structural void for the transport of the fluid with zero impermeability. Therefore, the impermeability can be described as
\begin{equation}\label{equ:ImpermeabilityInterp}
\left\{
\begin{split}
\alpha|_{\mathbf{x}_{\Gamma} \in \Gamma_D} & =
\left\{
\begin{split}
 & +\infty,~\gamma_p = 0 \\
 & 0      ,~\gamma_p = 1
\end{split}
\right. \\
\alpha|_{\mathbf{x}_{\Gamma} \in \Gamma_F} & = 0, ~ \gamma_p = 1 \\
\end{split}
\right.
\end{equation}
where $\gamma_p\in\left\{0,1\right\}$ is a binary distribution defined on $\Gamma$, with $0$ and $1$ representing the solid and fluid phases, respectively; $\Gamma_D$ is the design domain for the pattern of the surface flow, $\Gamma_F$ is the fluid domain with the material density enforced to be $\gamma_p = 1$, respectively, where $\Gamma_D$ and $\Gamma_F$ satisfy $\Gamma_D \cup \Gamma_F = \Gamma$ and $\Gamma_D \cap \Gamma_F = \emptyset$. Especially, $\Gamma$ is the design domain, when there is no enforced fluid domain, i.e., $\Gamma_F = \emptyset$ and $\Gamma = \Gamma_D$.

To avoid the numerical difficulty on solving a binary optimization problem, the binary variable $\gamma_p$ in the design domain is relaxed to vary continuously in $\left[0,1\right]$. The relaxed binary variable is referred to as the material density of the impermeability. Based on the description of the impermeability in Eq. \ref{equ:ImpermeabilityInterp}, the material interpolation of the impermeability can be implemented by using the convex and $q$-parameterized scheme \cite{Borrvall2003}:
\begin{equation}\label{equ:InterpolationForImpermeability}
\begin{split}
\alpha \left( \gamma_p \right) = \alpha_f + \left( \alpha_s - \alpha_f \right) q { 1 - \gamma_p \over q + \gamma_p }
\end{split}
\end{equation}
where $\alpha_s$ and $\alpha_f$ are the impermeability of the solid and fluid phases, respectively; $q$ is the parameter used to tune the convexity of this interpolation. For the fluid phase, the impermeability is zero, i.e., $\alpha_f = 0$. For the solid phase, $\alpha_s$ should be infinite theoretically; numerically, a finite value much larger than the fluid density $\rho$ is chosen for $\alpha_s$, to ensure the stability of the numerical implementation and approximate the solid phase with enough accuracy. Based on numerical tests, $q$ is valued as $1$ and $\alpha_s$ is chosen as $10^4\rho$ to satisfy $\alpha_s \gg \rho$ in this paper.

\subsection{Design variables}\label{sec:DesignVariablePattern}

In the fiber bundle topology optimization for the surface flow, two sets of design variables are required to be defined for the implicit 2-manifold and the pattern of the surface flow, respectively.

\subsubsection{Design variable for implicit 2-manifold}\label{subsec:DesignVariableImplicitManifold}

To describe the implicit 2-manifold, the design variable that takes continuous values in $\left[0,1\right]$ is defined on the base manifold. This design variable is used to describe the distribution of the normal displacement of the implicit 2-manifold relative to the base manifold. Equivalently, the pattern of the surface flow is defined on a variable design domain. Then, the result of the fiber bundle topology optimization can be regarded to be a two-order hierarchical structure composed of the base and secondary structures corresponding to the implicit 2-manifold and the pattern of the surface flow, respectively.

\begin{figure}[!htbp]
  \centering
  \includegraphics[width=0.8\textwidth]{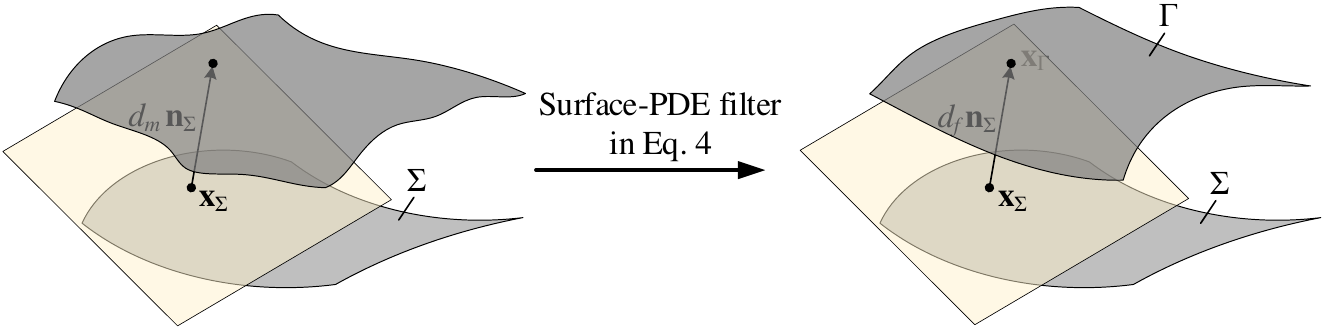}
  \caption{Sketch for the surface-PDE filter for the design variable of the implicit 2-manifold $\Gamma$ defined on the base manifold $\Sigma$.}\label{fig:SketchPDEFilterOptimalMatchingManifold}
\end{figure}

To control the smoothness of the implicit 2-manifold and ensure the well-poseness of the solution, a surface-PDE filter sketched in Fig. \ref{fig:SketchPDEFilterOptimalMatchingManifold} is imposed on the design variable of the implicit 2-manifold \cite{DengCMAME2020}:
      \begin{equation}\label{equ:PDEFilterzmBaseStructure1}
      \begin{split}
      & \left\{
        \begin{split}
          & - \mathrm{div}_\Sigma \left( r_m^2 \nabla_\Sigma d_f \right) + d_f = A_d \left( d_m - {1\over2} \right), ~ \forall \mathbf{x}_\Sigma \in \Sigma \\
          & \mathbf{n}_{\boldsymbol\tau_\Sigma} \cdot \nabla_\Sigma d_f = 0, ~ \forall \mathbf{x}_\Sigma \in \partial \Sigma \\
        \end{split}\right. \\
      \end{split}
      \end{equation}
where $d_m$ is the design variable for the implicit 2-manifold; $d_f$ is the filtered design variable; $r_m$ is the filter radius, and it is constant; $\Sigma$ is the base manifold used to define the implicit 2-manifold; $\mathbf{x}_\Sigma$ denotes a point on $\Sigma$; $\nabla_\Sigma$ and $\mathrm{div}_\Sigma$ are the tangential gradient operator and tangential divergence operator defined on $\Sigma$, respectively; $\mathbf{n}_{\boldsymbol\tau_\Sigma} = \mathbf{n}_\Sigma\times\boldsymbol\tau_\Sigma$ is the outward unitary conormal vector normal to $\partial\Sigma$ and tangent to $\Sigma$ at $\partial\Sigma$, with $\mathbf{n}_\Sigma$ and $\boldsymbol\tau_\Sigma$ representing the unitary normal vector on $\Sigma$ and the unitary tangential vector at $\partial\Sigma$, respectively; $A_d$ is a parameter used to specify the amplitude of the normal displacement of the implicit 2-manifold relative to the base manifold, and it is nonnegative ($A_d \geq 0$). Because $d_m$ is valued in $\left[0,1\right]$, $d_f$ is valued in $\left[ -A_d/2, A_d/2 \right]$.

After the filter operation, the implicit 2-manifold can be described by the filtered design variable as
\begin{equation}\label{equ:NormalDisplacementDistribution}
  \Gamma = \left\{ \mathbf{x}_\Gamma : \mathbf{x}_\Gamma = d_f \mathbf{n}_\Sigma + \mathbf{x}_\Sigma,~\forall \mathbf{x}_\Sigma \in \Sigma \right\}
\end{equation}
where $\Gamma$ is the implicit 2-manifold; $\mathbf{x}_\Gamma$ denotes a point on $\Gamma$. From Eq. \ref{equ:NormalDisplacementDistribution}, a differential homeomorphism can be determined corresponding to the bijection $d_f: \Sigma \mapsto \Gamma$ with $\mathbf{x}_\Gamma = d_f \mathbf{n}_\Sigma + \mathbf{x}_\Sigma $ for $\forall \mathbf{x}_\Sigma \in \Sigma$. Therefore, $\mathcal{H}\left(\Gamma\right)$ is homeomorphous to $\mathcal{H}\left(\Sigma\right)$. The Jacobian matrix of the homeomorphism (Eq. \ref{equ:NormalDisplacementDistribution}) for the implicit 2-manifold in the curvilinear coordinate system of the base manifold can be transformed into the following formulation:
\begin{equation}\label{equ:TransformedJacobian}
  {\partial \mathbf{x}_\Gamma \over \partial \mathbf{x}_\Sigma} = \nabla_\Sigma d_f \mathbf{n}_\Sigma^\mathrm{T} + d_f \nabla_\Sigma \mathbf{n}_\Sigma + \mathbf{I}, ~ \forall \mathbf{x}_\Sigma \in \Sigma
\end{equation}
with $\left| {\partial \mathbf{x}_\Gamma \over \partial \mathbf{x}_\Sigma} \right|$ representing its determinant.

The variational formulation of the surface-PDE filter in Eq. \ref{equ:PDEFilterzmBaseStructure1} is considered in the first order Sobolev space defined on $\Sigma$. It can be derived based on the Galerkin method as
\begin{equation}\label{equ:VariationalFormulationPDEFilterBaseManifold}
\begin{split}
     & \mathrm{Find}~d_f \in\mathcal{H}\left(\Sigma\right)~\mathrm{for}~d_m \in \mathcal{L}^2\left(\Sigma\right),~\mathrm{such~that} \\
     & \int_\Sigma r_m^2 \nabla_\Sigma d_f \cdot \nabla_\Sigma \tilde{d}_f + d_f\tilde{d}_f - A_d \left( d_m - {1\over2} \right) \tilde{d}_f \,\mathrm{d}\Sigma = 0,~\forall \tilde{d}_f \in \mathcal{H}\left(\Sigma\right)
\end{split}
\end{equation}
where $\tilde{d}_f$ is the test function of $d_f$; $\mathcal{H}\left(\Sigma\right)$ represents the first order Sobolev space defined on $\Sigma$; $\mathcal{L}^2\left(\Sigma\right)$ represents the second order Lebesque space defined on $\Sigma$.

\subsubsection{Design variable for pattern of surface flow}\label{subsec:DesignVariablePattern}

\begin{figure}[!htbp]
  \centering
  \includegraphics[width=1\textwidth]{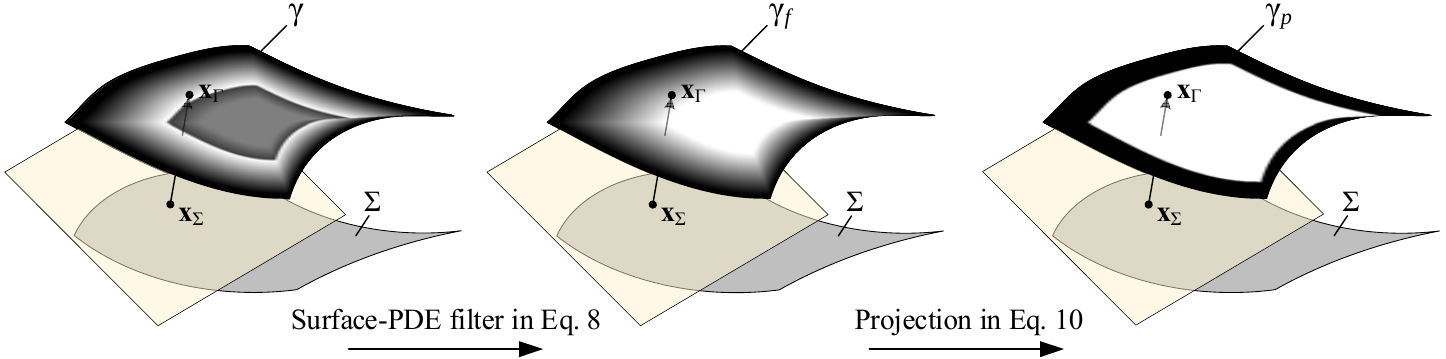}
  \caption{Sketch for the surface-PDE filter and projection operation for the design variable of the pattern of the surface flow.}\label{fig:SketchPDEFilterOptimalMatchingPattern}
\end{figure}

The pattern of the surface flow is represented by the material density defined on the implicit 2-manifold. The material density in Eqs. \ref{equ:ImpermeabilityInterp} and \ref{equ:InterpolationForImpermeability} is obtained by sequentially implementing the surface-PDE filter and the threshold projection on the design variable for the material density, as sketched in Fig. \ref{fig:SketchPDEFilterOptimalMatchingPattern}. This design variable is also valued continuously in $\left[0,1\right]$. Here, the threshold projection is used to remove the gray regions and control the minimum length scale in the derived pattern.

The surface-PDE filter for the design variable of the pattern is implemented by solving the following surface-PDE \cite{DengCMAME2020}:
      \begin{equation}\label{equ:PDEFilterGammaFilber}
        \left\{\begin{split}
        - \mathrm{div}_\Gamma \left( r_f^2 \nabla_\Gamma \gamma_f \right) + \gamma_f & = \gamma, ~\forall \mathbf{x}_\Gamma \in \Gamma \\
        \mathbf{n}_{\boldsymbol\tau_\Gamma} \cdot \nabla_\Gamma \gamma_f & = 0, ~\forall \mathbf{x}_\Gamma \in \partial\Gamma \\
        \end{split}\right.
      \end{equation}
where $\gamma$ is the design variable; $\gamma_f$ is the filtered design variable; $r_f$ is the filter radius, and it is constant; $\nabla_\Gamma$ and $\mathrm{div}_\Gamma$ are the tangential gradient operator and tangential divergence operator defined on the implicit 2-manifold $\Gamma$, respectively; $\mathbf{n}_{\boldsymbol\tau_\Gamma} = \mathbf{n}_\Gamma\times\boldsymbol\tau_\Gamma$ is the outward unitary conormal vector normal to $\partial\Gamma$ and tangent to $\Gamma$ at $\partial\Gamma$, with $\mathbf{n}_\Gamma$ and $\boldsymbol\tau_\Gamma$ representing the unitary normal vector on $\Gamma$ and the unitary tangential vector at $\partial\Gamma$, respectively. The threshold projection of the filtered design variable is implemented as \cite{WangStructMultidiscipOptim2011,GuestIntJNumerMethodsEng2004}
      \begin{equation}\label{equ:ProjectionGammaFilber}
        \gamma_p = { \tanh\left(\beta \xi\right) + \tanh\left(\beta \left(\gamma_f-\xi\right)\right) \over \tanh\left(\beta \xi\right) + \tanh\left(\beta \left(1-\xi\right)\right)}
      \end{equation}
where $\beta$ and $\xi$ are the parameters for the threshold projection, with values chosen based on numerical experiments \cite{GuestIntJNumerMethodsEng2004}.

The variational formulation of the surface-PDE filter is considered in the first order Sobolev space defined on $\Gamma$. It can be derived based on the Galerkin method as
\begin{equation}\label{equ:VariationalFormulationPDEFilter}
\begin{split}
     & \mathrm{Find}~\gamma_f \in\mathcal{H}\left(\Gamma\right)~\mathrm{for}~\gamma \in \mathcal{L}^2\left(\Gamma\right),~\mathrm{such~that} \\
     & \int_\Gamma r_f^2 \nabla_\Gamma \gamma_f \cdot \nabla_\Gamma \tilde{\gamma}_f + \gamma_f\tilde{\gamma}_f - \gamma \tilde{\gamma}_f \,\mathrm{d}\Gamma = 0,~\forall \tilde{\gamma}_f \in \mathcal{H}\left(\Gamma\right)
\end{split}
\end{equation}
where $\tilde{\gamma}_f$ is the test function of $\gamma_f$; $\mathcal{H}\left(\Gamma\right)$ represents the first order Sobolev space defined on $\Gamma$; $\mathcal{L}^2\left(\Gamma\right)$ represents the second order Lebesque space defined on $\Gamma$.

\subsubsection{Coupling of design variables}\label{subsec:CouplingDesignVariables}

The design variable introduced in Section \ref{subsec:DesignVariablePattern} for the pattern of the surface flow is defined on the implicit 2-manifold introduced in Section \ref{subsec:DesignVariableImplicitManifold}. Their coupling relation can be derived by transforming the tangential gradient operator $\nabla_\Gamma$, the tangential divergence operator $\mathrm{div}_\Gamma$ and the unitary normal $\mathbf{n}_\Gamma$ into the forms defined on the base manifold $\Sigma$.

The transformation of the tangential gradient operator $\nabla_\Gamma$ is implemented based on the following relation:
\begin{equation}\label{equ:TangentialOperatorRelation}
  \nabla_\Gamma = \nabla_\Sigma - \left(\mathbf{n}_\Gamma \cdot \nabla_\Sigma\right)\mathbf{n}_\Gamma = \mathbf{P}_\Gamma \nabla_\Sigma
\end{equation}
where $\mathbf{P}_\Gamma$ is the normal projector on the tangential space of $\Gamma$. The unitary normal vector $\mathbf{n}_\Gamma$ is transformed as
\begin{equation}\label{equ:UnitaryNormalGamma}
  \mathbf{n}_\Gamma^{\left( d_f \right)} = { \mathbf{n}_\Sigma - \nabla_\Sigma d_f \over \left\| \mathbf{n}_\Sigma - \nabla_\Sigma d_f \right\|_2}
\end{equation}
where $\left\| \cdot \right\|_2$ is the 2-norm of a vector. In Eq. \ref{equ:UnitaryNormalGamma}, the transformed unitary normal vector is distinguished from the original form by using the filtered design variable $d_f$ as the superscript, and this identification method is used in the following for the other transformed operators and variables. The normal projector $\mathbf{P}_\Gamma$ is sequentially transformed as
\begin{equation}\label{equ:NormalProjector}
  \mathbf{P}_\Gamma^{\left( d_f \right)} = \mathbf{I} - \mathbf{n}_\Gamma\mathbf{n}_\Gamma^\mathrm{T} = \mathbf{I} - { \left( \mathbf{n}_\Sigma - \nabla_\Sigma d_f \right) \left( \mathbf{n}_\Sigma - \nabla_\Sigma d_f \right)^\mathrm{T} \over \left( \mathbf{n}_\Sigma - \nabla_\Sigma d_f \right)^2}
\end{equation}
where $\mathbf{I}$ is the two-dimensional unitary tensor; and the superscript $\mathrm{T}$ represents the transposition operation of a vector or tensor.
The tangential gradient operator $\nabla_\Gamma$ can then be transformed as
\begin{equation}\label{equ:TransformedTangentialOperator}
  \nabla_\Gamma^{\left( d_f \right)} g = \nabla_\Sigma g - { \left( \mathbf{n}_\Sigma - \nabla_\Sigma d_f \right) \cdot \nabla_\Sigma g \over \left( \mathbf{n}_\Sigma - \nabla_\Sigma d_f \right)^2} \left( \mathbf{n}_\Sigma - \nabla_\Sigma d_f \right),~\forall g \in \mathcal{H}\left( \Sigma \right).
\end{equation}
Based on the transformed tangential gradient operator, the tangential divergence operator $\mathrm{div}_\Gamma$ can be transformed as
\begin{equation}\label{equ:TransformedDivergenceOperator}
\begin{split}
  \mathrm{div}_\Gamma^{\left( d_f \right)} \mathbf{g}
  = \mathrm{tr}\left( \nabla_\Gamma^{\left( d_f \right)} \mathbf{g} \right)
  = \mathrm{tr} \left( \nabla_\Sigma \mathbf{g} - { \left( \mathbf{n}_\Sigma - \nabla_\Sigma d_f \right) \cdot \nabla_\Sigma \mathbf{g} \over \left( \mathbf{n}_\Sigma - \nabla_\Sigma d_f \right)^2} \left( \mathbf{n}_\Sigma - \nabla_\Sigma d_f \right) \right), ~ \forall \mathbf{g} \in \left( \mathcal{H}\left( \Sigma \right) \right)^3
\end{split}
\end{equation}
where $\mathrm{tr}$ is the trace operator used to extract the trace of a tensor.

Because the tangential gradient operator $\nabla_\Gamma$ depends on $d_f$, its first-order variational to $d_f$ can be derived as
\begin{equation}\label{equ:FirstOrderVarisForTangentialOperators}
\begin{split}
& \mathrm{For} ~ \forall g \in \mathcal{H}\left( \Sigma \right), \\
& \nabla_\Gamma^{\left(d_f, \tilde{d}_f\right)} g = { \nabla_\Sigma \tilde{d}_f \cdot \nabla_\Sigma g \over \left( \mathbf{n}_\Sigma - \nabla_\Sigma d_f \right)^2} \left( \mathbf{n}_\Sigma - \nabla_\Sigma d_f \right) + { \left( \mathbf{n}_\Sigma - \nabla_\Sigma d_f \right) \cdot \nabla_\Sigma g \over \left( \mathbf{n}_\Sigma - \nabla_\Sigma d_f \right)^2} \nabla_\Sigma \tilde{d}_f \\
  & ~~~~~~ - 2 { \left[ \left( \mathbf{n}_\Sigma - \nabla_\Sigma d_f \right) \cdot \nabla_\Sigma g \right] \left[ \left( \mathbf{n}_\Sigma - \nabla_\Sigma d_f \right) \cdot \nabla_\Sigma \tilde{d}_f \right] \over \left( \mathbf{n}_\Sigma - \nabla_\Sigma d_f \right)^4} \left( \mathbf{n}_\Sigma - \nabla_\Sigma d_f \right), ~ \forall \tilde{d}_f \in \mathcal{H}\left( \Sigma \right).
\end{split}
\end{equation}
Similarly, the first-order variational of $\mathrm{div}_\Gamma$ to $d_f$ can be derived as
\begin{equation}\label{equ:FirstOrderVarisForDivergenceOperators}
\begin{split}
 & \mathrm{For} ~ \forall \mathbf{g} \in \left(\mathcal{H} \left( \Sigma \right)\right)^3, \\
  & \mathrm{div}_\Gamma^{\left(d_f,\tilde{d}_f\right)} \mathbf{g} = ~ \mathrm{tr} \Bigg( { \nabla_\Sigma \tilde{d}_f \cdot \nabla_\Sigma \mathbf{g} \over \left( \mathbf{n}_\Sigma - \nabla_\Sigma d_f \right)^2} \left( \mathbf{n}_\Sigma - \nabla_\Sigma d_f \right) + { \left( \mathbf{n}_\Sigma - \nabla_\Sigma d_f \right) \cdot \nabla_\Sigma \mathbf{g} \over \left( \mathbf{n}_\Sigma - \nabla_\Sigma d_f \right)^2} \nabla_\Sigma \tilde{d}_f \\
  & ~~~~~~ - 2 { \left[ \left( \mathbf{n}_\Sigma - \nabla_\Sigma d_f \right) \cdot \nabla_\Sigma \mathbf{g} \right] \left[ \left( \mathbf{n}_\Sigma - \nabla_\Sigma d_f \right) \cdot \nabla_\Sigma \tilde{d}_f \right] \over \left( \mathbf{n}_\Sigma - \nabla_\Sigma d_f \right)^4} \left( \mathbf{n}_\Sigma - \nabla_\Sigma d_f \right) \Bigg), ~ \forall \tilde{d}_f \in \mathcal{H}\left( \Sigma \right).
\end{split}
\end{equation}

Because $d_f$ is differential homeomorphism, it can induce the Riemannian metric. Then, the differential on the base manifold and implicit 2-manifold satisfies
\begin{equation}\label{equ:DiffRiemannian}
\left\{\begin{split}
  & \mathrm{d}\Gamma = \left| {\partial \mathbf{x}_\Gamma \over \partial \mathbf{x}_\Sigma} \right| \left\| {\partial \mathbf{x}_\Gamma \over \partial \mathbf{x}_\Sigma} \mathbf{n}_\Gamma^{\left( d_f \right)} \right\|_2^{-1} \mathrm{d}\Sigma \\
  & \mathrm{d}l_{\partial\Gamma} = \left\| \boldsymbol{\tau}_\Gamma \right\|_2 \left\| \left( \partial \mathbf{x}_\Gamma \over \partial \mathbf{x}_\Sigma \right)^{-1} \boldsymbol{\tau}_\Gamma \right\|_2^{-1} \mathrm{d}l_{\partial\Sigma} \\
\end{split}\right.,
\end{equation}
where $\mathrm{d}l_{\partial\Gamma}$ and $\mathrm{d}l_{\partial\Sigma}$ are the differential of the boundary curves of $\Gamma$ and $\Sigma$, respectively.

Based on the transformed tangential gradient operator in Eq. \ref{equ:TransformedTangentialOperator} and the homeomorphism between $\mathcal{H}\left(\Gamma\right)$ and $\mathcal{H}\left(\Sigma\right)$ described in Eq. \ref{equ:NormalDisplacementDistribution}, the coupling relation between the two sets of design variables can be derived by instituting Eq. \ref{equ:TransformedTangentialOperator} into Eq. \ref{equ:VariationalFormulationPDEFilter}:
\begin{equation}\label{equ:CoupledVariationalPDEFilter}
\begin{split}
     & \mathrm{Find}~\gamma_f \in\mathcal{H}\left(\Sigma\right)~\mathrm{for}~\gamma \in \mathcal{L}^2\left(\Sigma\right),~\mathrm{such~that} \\
     & \int_\Sigma \left( r_f^2 \nabla_\Gamma^{\left( d_f \right)} \gamma_f \cdot \nabla_\Gamma^{\left( d_f \right)} \tilde{\gamma}_f + \gamma_f\tilde{\gamma}_f - \gamma \tilde{\gamma}_f \right) \left| {\partial \mathbf{x}_\Gamma \over \partial \mathbf{x}_\Sigma} \right| \left\| {\partial \mathbf{x}_\Gamma \over \partial \mathbf{x}_\Sigma} \mathbf{n}_\Gamma^{\left( d_f \right)} \right\|_2^{-1} \,\mathrm{d}\Sigma = 0, ~ \forall \tilde{\gamma}_f \in \mathcal{H}\left(\Sigma\right)
\end{split}
\end{equation}
where the tangential gradient operator $\nabla_\Gamma$ on $\Gamma$ is replaced to be its transformed form $\nabla_\Gamma^{\left( d_f \right)}$ in Eq. \ref{equ:TransformedTangentialOperator}.

\subsection{Surface Navier-Stokes equations defined on implicit 2-manifold} \label{sec:SurfaceNSEqus}

The governing equations for the motion of a Newtonian surface fluid can be formulated intrinsically on a 2-manifold of codimension one in an Euclidian space.
Based on the conservation laws of momentum and mass, the surface Navier-Stokes equations can be derived to describe the incompressible surface flows \cite{ArroyoPRE2009,BrennerElsevier2013,RahimiSoftMatter2013}:
\begin{equation}\label{equ:UnsteadyNSequOnManifolds}
\begin{split}
\left.\begin{split}
\rho \left( \mathbf{u} \cdot \nabla_\Gamma \right) \mathbf{u} - \mathbf{P}_\Gamma \,\mathrm{div}_\Gamma \left[ \eta \left( \nabla_\Gamma \mathbf{u} + \nabla_\Gamma \mathbf{u}^\mathrm{T} \right) \right] + \nabla_\Gamma p & = - \alpha \mathbf{u} \\
- \mathrm{div}_\Gamma \mathbf{u} & = 0 \\
\mathbf{u} \cdot \mathbf{n}_\Gamma & = 0 \\
\end{split}\right\}~\forall \mathbf{x}_\Gamma \in \Gamma
\end{split}
\end{equation}
where $\mathbf{u}$ is the fluid velocity; $p$ is the fluid pressure; $\rho$ is the fluid density; $\eta$ is the dynamic viscosity; $\mathbf{u} \cdot \mathbf{n}_\Gamma = 0$ is the tangential constraint of the fluid velocity. The tangential constraint is imposed, because the fluid spatially flows on the 2-manifold $\Gamma$ and the fluid velocity is a vector in the tangential space of $\Gamma$.

To solve the surface Navier-Stokes equations, the fluid velocity and pressure are required to be specified at some boundaries, interfaces or points of the 2-manifold $\Gamma$:
\begin{equation}\label{equ:BNDConditionSurfaceNSEqu}
\left\{
\begin{split}
& \mathbf{u} = \mathbf{u}_{l_{v,\Gamma}}, ~ \forall \mathbf{x}_\Gamma \in l_{v,\Gamma} ~~~ \left( \mathrm{Inlet~or~interfacial~boundary~condition} \right) \\
& \mathbf{T}_\Gamma \cdot \mathbf{n}_{\boldsymbol\tau_\Gamma} = \mathbf{0}, ~\forall \mathbf{x}_\Gamma \in l_{s,\Gamma} ~~~ \left( \mathrm{Open~boundary~condition} \right) \\
& p = p_{0,\Gamma}, ~ \forall \mathbf{x}_\Gamma \in \mathcal{P}_\Gamma ~~~ \left( \mathrm{Point~condition} \right)
\end{split}
\right.
\end{equation}
where $\mathbf{u}_{l_{v,\Gamma}}$ is a known distribution of the fluid velocity depending on the specified fluid velocity $\mathbf{u}_{l_{v,\Sigma}}$ at $l_{v,\Sigma}$ representing a boundary or interface curve of $\Sigma$; $l_{v,\Gamma}$ satisfies $l_{v,\Gamma}\subset\partial\Gamma$ when $l_{v,\Gamma}$ is a boundary curve of $\Gamma$, and it satisfies $l_{v,\Gamma}\subset\Gamma$ when $l_{v,\Gamma}$ is an interface curve of $\Gamma$; $l_{s,\Gamma}$ is the boundary curve with open boundary condition, and it satisfies $l_{s,\Gamma}\subset\partial\Gamma$; $p_{0,\Gamma}$ is the known fluid pressure depending on the specified fluid pressure $p_{0,\Sigma}$ at $\mathcal{P}_\Sigma$ representing a finite point set on $\Sigma$; $\mathcal{P}_\Gamma$ is a finite point set on $\Gamma$. In Eq. \ref{equ:BNDConditionSurfaceNSEqu}, when the known fluid velocity is $\mathbf{0}$, the inlet or interfacial boundary condition degenerates into the no-slip boundary condition:
\begin{equation}\label{equ:NoSlipFluidVelocity}
  \mathbf{u} = \mathbf{0}, ~ \forall \mathbf{x}_\Gamma \in l_{v0,\Gamma}
\end{equation}
where $\mathbf{u}_{l_{v,\Gamma}}$ is equal to $\mathbf{0}$ on $l_{v0,\Gamma}\subset l_{v,\Gamma}$, and $l_{v0,\Gamma}$ is the no-slip part of the boundary curve.

The variational formulation of the surface Navier-Stokes equations is considered in the functional spaces without containing the tangential constraint of the fluid velocity. The tangential constraint of the fluid velocity is imposed by using the Lagrangian multiplier \cite{FriesIJNMF2018,ReutherPOF2018}. Based on the Galerkin method, the variational formulation of the surface Navier-Stokes equations can be derived as
\begin{equation}\label{equ:VariationalFormulationSurfaceNavierStokesEqus}
\begin{split}
  & \mathrm{Find}~\left\{\begin{split}
    & \mathbf{u}\in\left(\mathcal{H}\left(\Gamma\right)\right)^3~\mathrm{with} ~ \mathbf{u} = \mathbf{u}_{l_{v,\Gamma}}, ~ \forall \mathbf{x}_\Gamma \in l_{v,\Gamma} \\
  & p\in \mathcal{H}\left(\Gamma\right)~\mathrm{with}~p = p_{0,\Gamma}, ~{ \forall \mathbf{x}_\Gamma \in \mathcal{P}_\Gamma} \\
  & \lambda\in\mathcal{L}^2\left(\Gamma\right)~\mathrm{with}~ \lambda=0, ~ \forall \mathbf{x}_\Gamma \in l_{v,\Gamma}\\
  \end{split}\right., ~ \mathrm{such~that} \\
  &\int_\Gamma \rho \left( \mathbf{u} \cdot \nabla_\Gamma \right) \mathbf{u} \cdot \tilde{\mathbf{u}} + {\eta\over2} \left( \nabla_\Gamma \mathbf{u} + \nabla_\Gamma \mathbf{u}^\mathrm{T} \right) : \left( \nabla_\Gamma \tilde{\mathbf{u}} + \nabla_\Gamma \tilde{\mathbf{u}}^\mathrm{T} \right) - p\,\mathrm{div}_\Gamma \tilde{\mathbf{u}} - \tilde{p} \,\mathrm{div}_\Gamma \mathbf{u} \\
  & + \alpha \mathbf{u} \cdot \tilde{\mathbf{u}} + \lambda \tilde{\mathbf{u}} \cdot \mathbf{n}_\Gamma + \tilde{\lambda} \mathbf{u} \cdot \mathbf{n}_\Gamma \,\mathrm{d}\Gamma = 0, \\
  & \forall \tilde{\mathbf{u}} \in\left(\mathcal{H}\left(\Gamma\right)\right)^3,~\forall \tilde{p} \in \mathcal{H}\left(\Gamma\right), ~ \forall \tilde{\lambda} \in \mathcal{L}^2\left(\Gamma\right)
\end{split}
\end{equation}
where $\lambda$ is the Lagrange multiplier used to impose the tangential constraint of the fluid velocity; $\tilde{\mathbf{u}}$, $\tilde{p}$ and $\tilde{\lambda}$ are the test functions of $\mathbf{u}$, $p$ and $\lambda$, respectively. The Lagrangian multiplier in Eq. \ref{equ:VariationalFormulationSurfaceNavierStokesEqus} is used to impose the tangential constraint of the fluid velocity and acts as a distributed force in the normal direction of $\Gamma$. Such distributed force cancels out the centrifugal, Coriolis and Euler forces exerted on the fluid particles in the normal direction of $\Gamma$ to satisfy the tangential constraint.

Because $l_{v,\Gamma}$ is homeomorphous to $l_{v,\Sigma}$, $\mathcal{P}_\Gamma$ is homeomorphous to $\mathcal{P}_\Sigma$ and $l_{v,\Sigma}$ and $\mathcal{P}_\Sigma$ are fixed, $\mathbf{u}_{l_{v,\Gamma}}$ and $p_{0,\Gamma}$ are homeomorphous to $\mathbf{u}_{l_{v,\Sigma}}$ and $p_{0,\Sigma}$, respectively.
Then, based on the coupling relations in Section \ref{subsec:CouplingDesignVariables}, the variational formulation in Eq. \ref{equ:VariationalFormulationSurfaceNavierStokesEqus} can be transformed into the form defined on the base manifold $\Sigma$:
\begin{equation}\label{equ:TransformedVariationalFormulationSurfaceNSEqus}
\begin{split}
  & \mathrm{Find}~\left\{\begin{split}
    & \mathbf{u}\in\left(\mathcal{H}\left(\Sigma\right)\right)^3~\mathrm{with}~ \mathbf{u} = \mathbf{u}_{l_{v,\Sigma}}, ~ \forall \mathbf{x}_\Sigma \in l_{v,\Sigma} \\
  & p\in \mathcal{H}\left(\Sigma\right)~\mathrm{with}~p = p_{0,\Sigma}, ~{ \forall \mathbf{x}_\Sigma \in \mathcal{P}_\Sigma} \\
  & \lambda\in\mathcal{L}^2\left(\Sigma\right)~\mathrm{with}~ \lambda=0, ~ \forall \mathbf{x}_\Sigma \in l_{v,\Sigma}\\
  \end{split}\right.,~ \mathrm{such~that}\\
  & \int_\Sigma \Bigg[ \rho \left( \mathbf{u} \cdot \nabla_\Gamma^{\left(d_f\right)} \right) \mathbf{u} \cdot \tilde{\mathbf{u}} + {\eta\over2} \left( \nabla_\Gamma^{\left(d_f\right)} \mathbf{u} + \nabla_\Gamma^{\left(d_f\right)} \mathbf{u}^\mathrm{T} \right) : \left( \nabla_\Gamma^{\left(d_f\right)} \tilde{\mathbf{u}} + \nabla_\Gamma^{\left(d_f\right)} \tilde{\mathbf{u}}^\mathrm{T} \right) \\
  & - p \, \mathrm{div}_\Gamma^{\left( d_f \right)} \tilde{\mathbf{u}} - \tilde{p} \, \mathrm{div}_\Gamma^{\left( d_f \right)} \mathbf{u} + \alpha \mathbf{u} \cdot \tilde{\mathbf{u}} + \lambda \tilde{\mathbf{u}} \cdot { \mathbf{n}_\Sigma - \nabla_\Sigma d_f \over \left\| \mathbf{n}_\Sigma - \nabla_\Sigma d_f \right\|_2} \\
  & + \tilde{\lambda} \mathbf{u} \cdot { \mathbf{n}_\Sigma - \nabla_\Sigma d_f \over \left\| \mathbf{n}_\Sigma - \nabla_\Sigma d_f \right\|_2} \Bigg] \left| {\partial \mathbf{x}_\Gamma \over \partial \mathbf{x}_\Sigma} \right| \left\| {\partial \mathbf{x}_\Gamma \over \partial \mathbf{x}_\Sigma} \mathbf{n}_\Gamma^{\left( d_f \right)} \right\|_2^{-1} \,\mathrm{d}\Sigma = 0, \\
  & \forall \tilde{\mathbf{u}} \in\left(\mathcal{H}\left(\Sigma\right)\right)^3,~\forall \tilde{p} \in \mathcal{H}\left(\Sigma\right), ~ \forall \tilde{\lambda} \in \mathcal{L}^2\left(\Sigma\right).
\end{split}
\end{equation}

\subsection{Design objective in general form} \label{sec:GeneralDesignObjective}

\begin{figure}[!htbp]
  \centering
  \includegraphics[width=0.88\textwidth]{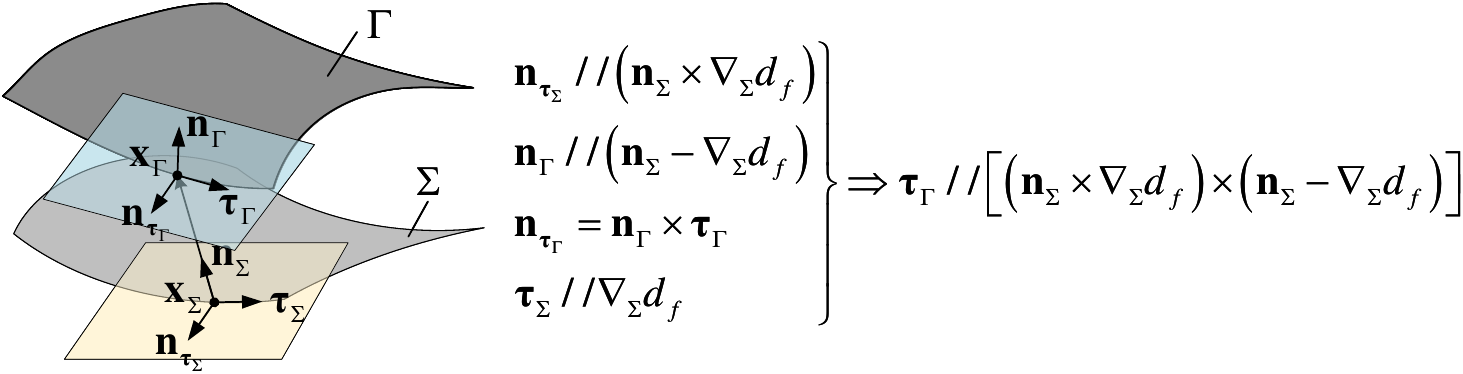}
  \caption{Sketch for relation among the unitary tangential vector $\boldsymbol\tau_\Gamma$ at $\partial\Gamma$, the unitary normal vector $\mathbf{n}_\Sigma$ on $\Sigma$ and the tangential gradient $\nabla_\Sigma d_f$.}\label{fig:DemonForTangentialVectorGamma}
\end{figure}

The design objective of the fiber bundle topology optimization problem for the surface flow is considered in the following general form:
\begin{equation}\label{equ:GeneralDesignObjective}
  J = \int_\Gamma A \left( \mathbf{u}, \nabla_\Gamma \mathbf{u}, p; \gamma_p \right) \,\mathrm{d}\Gamma + \int_{\partial\Gamma} B\left(\mathbf{u},p \right)\,\mathrm{d}l_{\partial\Gamma}
\end{equation}
where $A$ and $B$ are the integrands of the design objective. Based on the coupling relations in Section \ref{subsec:CouplingDesignVariables}, the design objective in Eq. \ref{equ:GeneralDesignObjective} can be transformed into the form defined on the base manifold $\Sigma$:
\begin{equation}\label{equ:TransformedGeneralDesignObjective}
\begin{split}
  J = & \int_\Sigma A \left( \mathbf{u}, \nabla_\Gamma^{\left( d_f \right)} \mathbf{u}, p; \gamma_p \right) \left| {\partial \mathbf{x}_\Gamma \over \partial \mathbf{x}_\Sigma} \right| \left\| {\partial \mathbf{x}_\Gamma \over \partial \mathbf{x}_\Sigma} \mathbf{n}_\Gamma^{\left( d_f \right)} \right\|_2^{-1} \,\mathrm{d}\Sigma \\
  & + \int_{\partial\Sigma} B\left(\mathbf{u},p \right) \left\| \boldsymbol{\tau}_\Gamma \right\|_2 \left\| \left( \partial \mathbf{x}_\Gamma \over \partial \mathbf{x}_\Sigma \right)^{-1} \boldsymbol{\tau}_\Gamma \right\|_2^{-1} \,\mathrm{d}l_{\partial\Sigma}.
\end{split}
\end{equation}
In Eq. \ref{equ:TransformedGeneralDesignObjective}, the unitary tangential vector $\boldsymbol\tau_\Gamma$ at $\partial\Gamma$ satisfies the relation sketched in Fig. \ref{fig:DemonForTangentialVectorGamma}:
\begin{equation}\label{equ:TauGamma}
  \boldsymbol\tau_\Gamma // \left[ \left( \mathbf{n}_\Sigma \times \nabla_\Sigma d_f \right) \times \left( \mathbf{n}_\Sigma - \nabla_\Sigma d_f\right) \right],
\end{equation}
Eq. \ref{equ:TransformedGeneralDesignObjective} can then be transformed into
\begin{equation}\label{equ:FurtherTransformedGeneralDesignObjective}
\begin{split}
  J = & \int_\Sigma A \left( \mathbf{u}, \nabla_\Gamma^{\left( d_f \right)} \mathbf{u}, p; \gamma_p \right) \left| {\partial \mathbf{x}_\Gamma \over \partial \mathbf{x}_\Sigma} \right| \left\| {\partial \mathbf{x}_\Gamma \over \partial \mathbf{x}_\Sigma} \mathbf{n}_\Gamma^{\left( d_f \right)} \right\|_2^{-1} \,\mathrm{d}\Sigma \\
  & + \int_{\partial\Sigma} B\left( \mathbf{u},p \right) \left\| \left( \mathbf{n}_\Sigma \times \nabla_\Sigma d_f \right) \times \left( \mathbf{n}_\Sigma - \nabla_\Sigma d_f\right) \right\|_2 \\
  & \left\| \left( \partial \mathbf{x}_\Gamma \over \partial \mathbf{x}_\Sigma \right)^{-1} \left[ \left( \mathbf{n}_\Sigma \times \nabla_\Sigma d_f \right) \times \left( \mathbf{n}_\Sigma - \nabla_\Sigma d_f\right)\right] \right\|_2^{-1} \,\mathrm{d}l_{\partial\Sigma}.
\end{split}
\end{equation}
Based on the transformed design objective, the adjoint analysis of the fiber bundle topology optimization problem can then be implemented on the functional spaces defined on the base manifold.

\subsection{Fiber bundle topology optimization problem} \label{sec:MatchingOptimizationSurfaceNSEqus}

\begin{figure}[!htbp]
  \centering
  \includegraphics[width=1\textwidth]{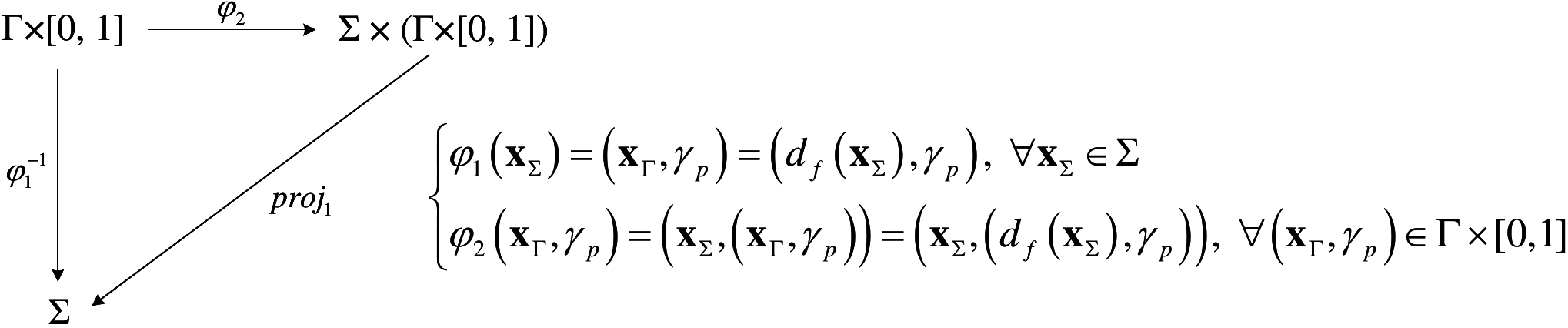}
  \caption{Diagram for the fiber bundle composed of the base manifold, the implicit 2-manifold and the pattern of the surface flow.}\label{fig:DiagramFiberBundle1}
\end{figure}

The fiber bundle of the surface flow is composed of the base manifold together with the implicit 2-manifold and the pattern, where $\Sigma$ is the base manifold and $\Gamma \times \left[0,1\right]$ is the fiber, respectively. It can be expressed as $\left(\Sigma \times \left(\Gamma \times \left[0,1\right]\right), \Sigma, proj_1, \Gamma \times \left[0,1\right] \right)$ with the diagram shown in Fig. \ref{fig:DiagramFiberBundle1}, where $proj_1$ is the natural projection $proj_1: \Sigma \times \left(\Gamma \times \left[0,1\right]\right) \mapsto \Sigma$ satisfying $ proj_1\left(\mathbf{x}_\Sigma,\left(\mathbf{x}_\Gamma, \gamma_p \right)\right) = proj_1\left(\mathbf{x}_\Sigma, \left( d_f\left(\mathbf{x}_\Sigma \right), \gamma_p \right)\right) = \mathbf{x}_\Sigma$ for $\forall \mathbf{x}_\Sigma\in\Sigma$, $\varphi_1$ is the homeomorphous map $\varphi_1 : \Sigma \mapsto \Gamma \times \left[0,1\right]$ satisfying $\varphi_1 \left( \mathbf{x}_\Sigma \right) = \left( \mathbf{x}_\Gamma, \gamma_p \right) = \left( d_f \left( \mathbf{x}_\Sigma \right) , \gamma_p \right)$ for $\forall \mathbf{x}_\Sigma\in\Sigma$, and $\varphi_2$ is the homeomorphous map $\varphi_2 : \Gamma \times \left[0,1\right] \mapsto \Sigma \times \left(\Gamma \times \left[0,1\right]\right)$ satisfying $\varphi_2 \left( \mathbf{x}_\Gamma, \gamma_p \right) = \left( \mathbf{x}_\Sigma, \left( \mathbf{x}_\Gamma, \gamma_p \right) \right) = \left( \mathbf{x}_\Sigma, \left( d_f \left(\mathbf{x}_\Sigma\right), \gamma_p \right) \right)$ for $\forall \left( \mathbf{x}_\Gamma, \gamma_p \right)\in\Gamma\times\left[0,1\right]$.

Based on the above introduction, the fiber bundle topology optimization problem of the surface flow can be constructed as
\begin{equation}\label{equ:VarProToopSurfaceFlows}
\begin{split}
  & \mathrm{Find}~ \left\{\begin{split}
  & \gamma: \Gamma \mapsto \left[0,1\right] \\
  & d_m: \Sigma \mapsto \left[0,1\right]\end{split}\right.~ \mathrm{for} ~
  \left(\Sigma \times \left(\Gamma \times \left[0,1\right]\right), \Sigma, proj_1, \Gamma \times \left[0,1\right] \right) \\
  & \mathrm{to} ~ \mathrm{minimize}~{J \over J_0}~ \mathrm{with} ~ J = \int_\Gamma A \left( \mathbf{u}, \nabla_\Gamma \mathbf{u}, p; \gamma_p \right) \,\mathrm{d}\Gamma + \int_{\partial\Gamma} B\left(\mathbf{u},p \right)\,\mathrm{d}l_{\partial\Gamma} \\
  & \mathrm{constrained} ~ \mathrm{by} \\
  & \left\{\begin{split}
  & \begin{split}
  & \left\{\begin{split}
    & \left.\begin{split}
       & \rho \left( \mathbf{u} \cdot \nabla_\Gamma \right) \mathbf{u} - \mathbf{P}_\Gamma \,\mathrm{div}_\Gamma \left[ \eta \left( \nabla_\Gamma \mathbf{u} + \nabla_\Gamma \mathbf{u}^\mathrm{T} \right) \right] + \nabla_\Gamma p = - \alpha \mathbf{u} \\
       & - \mathrm{div}_\Gamma \mathbf{u} = 0 \\
       & \mathbf{u} \cdot \mathbf{n}_\Gamma = 0 \\
    \end{split}\right\} ~ \forall \mathbf{x}_\Gamma \in \Gamma \\
    \end{split}\right.
    \end{split} \\
  & \left\{\begin{split}
        & - \mathrm{div}_\Gamma \left( r_f^2 \nabla_\Gamma \gamma_f \right) + \gamma_f = \gamma,~\forall \mathbf{x}_\Gamma \in \Gamma \\
        & \mathbf{n}_{\boldsymbol\tau_\Gamma} \cdot \nabla_\Gamma \gamma_f = 0,~\forall \mathbf{x}_\Gamma \in \partial\Gamma \\
    \end{split}\right. \\
  & \left\{
        \begin{split}
          & - \mathrm{div}_\Sigma \left( r_m^2 \nabla_\Sigma d_f \right) + d_f = A_d \left( d_m - {1\over2} \right), ~ \forall \mathbf{x}_\Sigma \in \Sigma \\
          & \mathbf{n}_{\boldsymbol\tau_\Sigma} \cdot \nabla_\Sigma d_f = 0, ~ \forall \mathbf{x}_\Sigma \in \partial \Sigma \\
        \end{split}\right. \\
  & \Gamma = \left\{ \mathbf{x}_\Gamma : \mathbf{x}_\Gamma = d_f \mathbf{n}_\Sigma + \mathbf{x}_\Sigma,~\forall \mathbf{x}_\Sigma \in \Sigma \right\} \\
  & \gamma_p = { \tanh\left(\beta \xi\right) + \tanh\left(\beta \left(\gamma_f-\xi\right)\right) \over \tanh\left(\beta \xi\right) + \tanh\left(\beta \left(1-\xi\right)\right)} \\
  & s \leq s_0, ~\mathrm{with} ~ s = {1 \over \left|\Gamma\right|} \int_\Gamma \gamma_p \,\mathrm{d}\Gamma, ~ \left|\Gamma\right| = \int_\Gamma 1 \,\mathrm{d}\Gamma ~~~ (\mathrm{Area} ~\mathrm{constraint}) \\
  & \left|v - v_0\right| \leq 10^{-3}, ~ \mathrm{with} ~ v = {1 \over \left|\Sigma\right|} \int_\Sigma d_f \,\mathrm{d}\Sigma, ~ \left|\Sigma\right| = \int_\Sigma 1 \,\mathrm{d}\Sigma ~~~ (\mathrm{Volume} ~\mathrm{constraint}) \\
\end{split}\right. \\
\end{split}
\end{equation}
where $J_0$ is the value of the design objective corresponding to the initial distribution of the design variables;
to regularize this optimization problem, area and volume constraints are imposed on the pattern of the surface flow and implicit 2-manifold, respectively; $s$ is the area fraction of the pattern of the surface flow; $v$ is volume fraction of spacial domain enclosed by the implicit 2-manifold and the base manifold; $s_0 \in \left(0,1\right)$ and $v_0 \in \left[0,1\right)$ are the specified area and volume fractions, respectively.

The coupling relations among the variables, functions, tangential divergence operator and tangential gradient operator in the fiber bundle topology optimization problem are illustrated by the following arrow chart:
\[\begin{array}{cccccccc}
 \textcolor{blue}{d_m} & \xrightarrow{\mathrm{Eq.~}\ref{equ:PDEFilterzmBaseStructure1}} & d_f & \xrightarrow{\mathrm{Eq.~}\ref{equ:VarProToopSurfaceFlows}} & \textcolor[rgb]{0.50,0.00,0.00}{v} & \\
 & & \bigg\downarrow\vcenter{\rlap{\small{Eq.~\ref{equ:TransformedTangentialOperator}}}} & & \\
 & & \left\{ \mathrm{div}_\Gamma, \nabla_\Gamma, \mathbf{n}_\Gamma \right\} & \xrightarrow{\mathrm{Eq.~}\ref{equ:UnsteadyNSequOnManifolds}} & \left\{ \mathbf{u},~p,~\lambda \right\} & \xrightarrow{\mathrm{Eq.~}\ref{equ:VarProToopSurfaceFlows}} & \textcolor[rgb]{0.50,0.00,0.00}{J} \\
 & & \bigg\downarrow\vcenter{\rlap{\small{Eq.~\ref{equ:PDEFilterGammaFilber}}}} & &  \bigg\uparrow\vcenter{\rlap{\small{Eq.~\ref{equ:UnsteadyNSequOnManifolds}}}}   \\
 \textcolor{blue}{\gamma} & \xrightarrow{\mathrm{Eq.~}\ref{equ:PDEFilterGammaFilber}} & \gamma_f & \xrightarrow{\mathrm{Eq.~}\ref{equ:ProjectionGammaFilber}} & \gamma_p & \xrightarrow{\mathrm{Eq.~}\ref{equ:VarProToopSurfaceFlows}} & \textcolor[rgb]{0.50,0.00,0.00}{s} \\
\end{array}\]\\
where the design variables $d_m$ and $\gamma$, marked in blue, are the inputs; the design objective $J$, the area fraction $s$ and the volume fraction $v$, marked in red, are the outputs.

\subsection{Adjoint analysis} \label{sec:AdjointAnalysisSurfaceNSEqus}

The fiber bundle topology optimization problem in Eq. \ref{equ:VarProToopSurfaceFlows} can be solved by using a gradient information-based iterative procedure, where the adjoint sensitivities are used to determine the relevant gradient information. The adjoint analysis is implemented for the design objective and the area and volume constraints to derive the adjoint sensitivities. The details for the adjoint analysis have been provided in the appendix in Section \ref{sec:Appendix}.

Based on the continuous adjoint analysis method \cite{HinzeSpringer2009}, the adjoint sensitivity of the design objective $J$ is derived as
\begin{equation}\label{equ:AdjSensitivityGaDm}
\begin{split}
\delta J = \int_\Sigma - \gamma_{fa} \tilde{\gamma} \left| {\partial \mathbf{x}_\Gamma \over \partial \mathbf{x}_\Sigma} \right| \left\| {\partial \mathbf{x}_\Gamma \over \partial \mathbf{x}_\Sigma} \mathbf{n}_\Gamma^{\left( d_f \right)} \right\|_2^{-1} - A_d d_{fa} \tilde{d}_m \,\mathrm{d}\Sigma,~ \forall \left( \tilde{\gamma}, \tilde{d}_m \right) \in \left(\mathcal{L}^2\left(\Sigma\right)\right)^2
\end{split}
\end{equation}
where $\gamma_{fa}$ and $d_{fa}$ are the adjoint variables of the filtered design variables $\gamma_f$ and $d_f$, respectively; $\delta$ is the first-order variational operator. The adjoint variables can be derived from the adjoint equations in the variational formulations. The variational formulation for the adjoint equations of the surface Naiver-Stokes equations is derived as
  \begin{equation}\label{equ:AdjSurfaceNavierStokesEqusJObjective}
  \begin{split}
  & \mathrm{Find} ~ \left\{\begin{split}
  & \mathbf{u}_a \in\left(\mathcal{H}\left(\Sigma\right)\right)^3~\mathrm{with}~ \mathbf{u}_a = \mathbf{0}, ~ {\forall \mathbf{x} \in l_{v,\Sigma} } \\
  & p_a \in \mathcal{H}\left(\Sigma\right)~\mathrm{with}~ p_a=0,~ \forall \mathbf{x} \in \mathcal{P}_\Sigma \\
  & \lambda_a \in \mathcal{L}^2\left(\Sigma\right)~\mathrm{with}~ \lambda_a = 0,~ \forall \mathbf{x} \in l_{v,\Sigma} \\
  \end{split}\right.,~\mathrm{such~that} \\
  & \int_\Sigma \bigg[ {\partial A \over \partial \mathbf{u}} \cdot \tilde{\mathbf{u}}_a + {\partial A \over \partial \nabla_\Gamma^{\left( d_f \right)} \mathbf{u} } : \nabla_\Gamma^{\left( d_f \right)} \tilde{\mathbf{u}}_a + {\partial A \over \partial p } \tilde{p}_a + \rho \left( \tilde{\mathbf{u}}_a \cdot \nabla_\Gamma^{\left( d_f \right)} \right) \mathbf{u} \cdot \mathbf{u}_a + \rho \left( \mathbf{u} \cdot \nabla_\Gamma^{\left( d_f \right)} \right) \tilde{\mathbf{u}}_a \cdot \mathbf{u}_a \\
  & + {\eta\over2} \left( \nabla_\Gamma^{\left( d_f \right)} \tilde{\mathbf{u}}_a + \nabla_\Gamma^{\left( d_f \right)} \tilde{\mathbf{u}}_a^\mathrm{T} \right) : \left( \nabla_\Gamma^{\left( d_f \right)} \mathbf{u}_a + \nabla_\Gamma^{\left( d_f \right)} \mathbf{u}_a^\mathrm{T} \right) - \tilde{p}_a \, \mathrm{div}_\Gamma^{\left( d_f \right)} \mathbf{u}_a - p_a \, \mathrm{div}_\Gamma^{\left( d_f \right)} \tilde{\mathbf{u}}_a \\
  & + \alpha \tilde{\mathbf{u}}_a \cdot \mathbf{u}_a + \left( \tilde{\lambda}_a \mathbf{u}_a + \lambda_a \tilde{\mathbf{u}}_a \right) \cdot {\nabla_\Sigma d_f + \mathbf{n}_\Sigma \over \left\| \nabla_\Sigma d_f + \mathbf{n}_\Sigma \right\|_2 } \bigg] \left| {\partial \mathbf{x}_\Gamma \over \partial \mathbf{x}_\Sigma} \right| \left\| {\partial \mathbf{x}_\Gamma \over \partial \mathbf{x}_\Sigma} \mathbf{n}_\Gamma^{\left( d_f \right)} \right\|_2^{-1} \,\mathrm{d}\Sigma \\
  & + \int_{\partial\Sigma} \left( {\partial B \over \partial \mathbf{u}} \cdot \tilde{\mathbf{u}}_a + {\partial B \over \partial p} \tilde{p}_a \right) \left\| \left( \mathbf{n}_\Sigma \times \nabla_\Sigma d_f \right) \times \left( \mathbf{n}_\Sigma - \nabla_\Sigma d_f\right) \right\|_2 \\
  & \left\| \left( \partial \mathbf{x}_\Gamma \over \partial \mathbf{x}_\Sigma \right)^{-1} \left[ \left( \mathbf{n}_\Sigma \times \nabla_\Sigma d_f \right) \times \left( \mathbf{n}_\Sigma - \nabla_\Sigma d_f\right)\right] \right\|_2^{-1} \,\mathrm{d}l_{\partial\Sigma} = 0, \\
  & \forall \tilde{\mathbf{u}}_a \in\left(\mathcal{H}\left(\Sigma\right)\right)^3, ~ \tilde{p}_a \in \mathcal{H}\left(\Sigma\right), ~ \tilde{\lambda}_a \in \mathcal{L}^2\left(\Sigma\right)
  \end{split}
  \end{equation}
where $\mathbf{u}_a$, $p_a$ and $\lambda_a$ are the adjoint variables of $\mathbf{u}$, $p$ and $\lambda$, respectively; $\tilde{\mathbf{u}}_a$, $\tilde{p}_a$ and $\tilde{\lambda}_a$ are the test functions of $\mathbf{u}_a$, $p_a$ and $\lambda_a$, respectively.
The variational formulations for the adjoint equations of the surface-PDE filters for $\gamma$ and $d_m$ are derived as
\begin{equation}\label{equ:AdjPDEFilterJObjectiveGa}
\begin{split}
  & \mathrm{Find}~\gamma_{fa}\in\mathcal{H}\left(\Sigma\right),~\mathrm{such~that} \\
  & \int_\Sigma \left( {\partial A \over \partial \gamma_p } {\partial \gamma_p \over \partial \gamma_f } \tilde{\gamma}_{fa} + {\partial \alpha \over \partial \gamma_p} {\partial \gamma_p \over \partial \gamma_f} \mathbf{u} \cdot \mathbf{u}_a \tilde{\gamma}_{fa} + r_f^2 \nabla_\Gamma^{\left( d_f \right)} \gamma_{fa} \cdot \nabla_\Gamma^{\left( d_f \right)} \tilde{\gamma}_{fa} + \gamma_{fa} \tilde{\gamma}_{fa} \right) \\
  & \left| {\partial \mathbf{x}_\Gamma \over \partial \mathbf{x}_\Sigma} \right| \left\| {\partial \mathbf{x}_\Gamma \over \partial \mathbf{x}_\Sigma} \mathbf{n}_\Gamma^{\left( d_f \right)} \right\|_2^{-1} \,\mathrm{d}\Sigma = 0, ~ \forall \tilde{\gamma}_{fa} \in \mathcal{H}\left(\Sigma\right)
\end{split}
\end{equation}
and
\begin{equation}\label{equ:AdjPDEFilterJObjectiveDm}
\begin{split}
  & \mathrm{Find}~d_{fa}\in\mathcal{H}\left(\Sigma\right),~\mathrm{such~that} \\
  & \int_\Sigma \Bigg[ {\partial A \over \partial \nabla_\Gamma^{\left( d_f \right)} \mathbf{u} } : \nabla_\Gamma^{\left( d_f, \tilde{d}_{fa} \right)} \mathbf{u} + \rho \left( \mathbf{u} \cdot \nabla_\Gamma^{\left( d_f, \tilde{d}_{fa} \right)} \right) \mathbf{u} \cdot \mathbf{u}_a + {\eta\over 2} \left( \nabla_\Gamma^{\left( d_f, \tilde{d}_{fa} \right)} \mathbf{u} + \nabla_\Gamma^{\left( d_f, \tilde{d}_{fa} \right)} \mathbf{u}^\mathrm{T} \right) \\
  & : \left( \nabla_\Gamma^{\left( d_f \right)} \mathbf{u}_a + \nabla_\Gamma^{\left( d_f \right)} \mathbf{u}_a^\mathrm{T} \right) + {\eta\over2} \left( \nabla_\Gamma^{\left( d_f \right)} \mathbf{u} + \nabla_\Gamma^{\left( d_f \right)} \mathbf{u}^\mathrm{T} \right) : \left( \nabla_\Gamma^{\left( d_f, \tilde{d}_{fa} \right)} \mathbf{u}_a + \nabla_\Gamma^{\left( d_f, \tilde{d}_{fa} \right)} \mathbf{u}_a^\mathrm{T} \right) \\
  & - p \, \mathrm{div}_\Gamma^{\left( d_f, \tilde{d}_{fa} \right)} \mathbf{u}_a - p_a \, \mathrm{div}_\Gamma^{\left( d_f, \tilde{d}_{fa} \right)} \mathbf{u} + \left( \lambda \mathbf{u}_a + \lambda_a \mathbf{u} \right) \\
  & \cdot \Bigg({\nabla_\Sigma \tilde{d}_{fa} \over \left\| \nabla_\Sigma d_f + \mathbf{n}_\Sigma \right\|_2 } - {\nabla_\Sigma d_f + \mathbf{n}_\Sigma \over \left( \nabla_\Sigma d_f + \mathbf{n}_\Sigma \right)^2 } { \left( \nabla_\Sigma d_f + \mathbf{n}_\Sigma \right) \cdot \nabla_\Sigma \tilde{d}_{fa} \over \left\| \nabla_\Sigma d_f + \mathbf{n}_\Sigma \right\|_2 } \Bigg) \\
  & + r_f^2 \bigg( \nabla_\Gamma^{\left( d_f, \tilde{d}_{fa} \right)} \gamma_f \cdot \nabla_\Gamma^{\left( d_f \right)} \gamma_{fa} + \nabla_\Gamma^{\left( d_f \right)} \gamma_f \cdot \nabla_\Gamma^{\left( d_f, \tilde{d}_{fa} \right)} \gamma_{fa} \bigg) \Bigg] \left| {\partial \mathbf{x}_\Gamma \over \partial \mathbf{x}_\Sigma} \right| \left\| {\partial \mathbf{x}_\Gamma \over \partial \mathbf{x}_\Sigma} \mathbf{n}_\Gamma^{\left( d_f \right)} \right\|_2^{-1} \\
  & + \Bigg[ A + \rho \left( \mathbf{u} \cdot \nabla_\Gamma^{\left( d_f \right)} \right) \mathbf{u} \cdot \mathbf{u}_a + {\eta\over2} \Big( \nabla_\Gamma^{\left( d_f \right)} \mathbf{u} + \nabla_\Gamma^{\left( d_f \right)} \mathbf{u}^\mathrm{T} \Big) : \left( \nabla_\Gamma^{\left( d_f \right)} \mathbf{u}_a + \nabla_\Gamma^{\left( d_f \right)} \mathbf{u}_a^\mathrm{T} \right) \\
  & - p \, \mathrm{div}_\Gamma^{\left( d_f \right)} \mathbf{u}_a - p_a \, \mathrm{div}_\Gamma^{\left( d_f \right)} \mathbf{u} + \alpha \mathbf{u} \cdot \mathbf{u}_a + \left( \lambda \mathbf{u}_a + \lambda_a \mathbf{u} \right) \cdot {\nabla_\Sigma d_f + \mathbf{n}_\Sigma \over \left\| \nabla_\Sigma d_f + \mathbf{n}_\Sigma \right\|_2} \\
  & + \left( r_f^2 \nabla_\Gamma^{\left( d_f \right)} \gamma_f \cdot \nabla_\Gamma^{\left( d_f \right)} \gamma_{fa} + \gamma_f \gamma_{fa} - \gamma \gamma_{fa} \right) \Bigg] \Bigg( {\partial \left| {\partial \mathbf{x}_\Gamma \over \partial \mathbf{x}_\Sigma} \right| \left\| {\partial \mathbf{x}_\Gamma \over \partial \mathbf{x}_\Sigma} \mathbf{n}_\Gamma^{\left( d_f \right)} \right\|_2^{-1} \over \partial d_f} \tilde{d}_{fa} \\
  & + {\partial \left| {\partial \mathbf{x}_\Gamma \over \partial \mathbf{x}_\Sigma} \right| \left\| {\partial \mathbf{x}_\Gamma \over \partial \mathbf{x}_\Sigma} \mathbf{n}_\Gamma^{\left( d_f \right)} \right\|_2^{-1} \over \partial \nabla_\Sigma d_f} \cdot \nabla_\Sigma \tilde{d}_{fa} \Bigg) + r_m^2 \nabla_\Sigma d_{fa} \cdot \nabla_\Sigma \tilde{d}_{fa} + d_{fa} \tilde{d}_{fa} \,\mathrm{d}\Sigma \\
  & + \int_{\partial\Sigma} B { \partial \left\| \left( \mathbf{n}_\Sigma \times \nabla_\Sigma d_f \right) \times \left( \mathbf{n}_\Sigma - \nabla_\Sigma d_f\right) \right\|_2 \left\| \left( \partial \mathbf{x}_\Gamma \over \partial \mathbf{x}_\Sigma \right)^{-1} \left[ \left( \mathbf{n}_\Sigma \times \nabla_\Sigma d_f \right) \times \left( \mathbf{n}_\Sigma - \nabla_\Sigma d_f\right)\right] \right\|_2^{-1} \over \partial d_f } \tilde{d}_{fa} \\
  & + B { \partial \left\| \left( \mathbf{n}_\Sigma \times \nabla_\Sigma d_f \right) \times \left( \mathbf{n}_\Sigma - \nabla_\Sigma d_f\right) \right\|_2 \left\| \left( \partial \mathbf{x}_\Gamma \over \partial \mathbf{x}_\Sigma \right)^{-1} \left[ \left( \mathbf{n}_\Sigma \times \nabla_\Sigma d_f \right) \times \left( \mathbf{n}_\Sigma - \nabla_\Sigma d_f\right)\right] \right\|_2^{-1} \over \partial \nabla_\Sigma d_f } \\
  & \cdot \nabla_\Sigma \tilde{d}_{fa} \,\mathrm{d}l_{\partial\Sigma} = 0, ~ \forall \tilde{d}_{fa} \in \mathcal{H}\left(\Sigma\right)
\end{split}
\end{equation}
where $\tilde{\gamma}_{fa}$ and $\tilde{d}_{fa}$ are the test functions of $\gamma_{fa}$ and $d_{fa}$, respectively.

For the area constraint, the adjoint sensitivity of the area $s$ is derived as
\begin{equation}\label{equ:AdjSensitivityGaDmAreaConstr}
\begin{split}
\delta s = \delta { s \left| \Gamma \right| \over \left| \Gamma \right| } = { 1 \over \left| \Gamma \right| } \delta \left( s \left| \Gamma \right| \right) - { s \over \left| \Gamma \right| } \delta \left| \Gamma \right|.
\end{split}
\end{equation}
In Eq. \ref{equ:AdjSensitivityGaDmAreaConstr}, the adjoint sensitivity $\delta \left( s \left| \Gamma \right| \right)$ can be derived based on the adjoint analysis of $ s \left| \Gamma \right| = \int_\Gamma \gamma_p \,\mathrm{d}\Gamma$:
\begin{equation}\label{equ:AdjSensAreaConstr1}
\begin{split}
  & \delta \left( s \left| \Gamma \right| \right) = \int_\Sigma - \gamma_{fa} \tilde{\gamma} \left| {\partial \mathbf{x}_\Gamma \over \partial \mathbf{x}_\Sigma} \right| \left\| {\partial \mathbf{x}_\Gamma \over \partial \mathbf{x}_\Sigma} \mathbf{n}_\Gamma^{\left( d_f \right)} \right\|_2^{-1} - A_d d_{fa} \tilde{d}_m \,\mathrm{d}\Sigma, ~ \forall \left( \tilde{\gamma}, \tilde{d}_m \right) \in \left(\mathcal{L}^2\left(\Sigma\right)\right)^2.
\end{split}
\end{equation}
In Eq. \ref{equ:AdjSensAreaConstr1}, the adjoint variables $\gamma_{fa}$ and $d_{fa}$ are derived by solving the variational formulations for the adjoint equations of the surface-PDE filters for $\gamma$ and $d_m$, respectively:
\begin{equation}\label{equ:AdjEquAreaConstr1Ga}
\begin{split}
  & \mathrm{Find} ~ \gamma_{fa} \in \mathcal{H}\left( \Sigma \right),~ \mathrm{such~that} \\
  & \int_\Sigma \left( {\partial \gamma_p \over \partial \gamma_f} \tilde{\gamma}_{fa} + r_f^2 \nabla_\Gamma^{\left( d_f \right)} \gamma_{fa} \cdot \nabla_\Gamma^{\left( d_f \right)} \tilde{\gamma}_{fa} + \tilde{\gamma}_{fa} \gamma_{fa} \right) \left| {\partial \mathbf{x}_\Gamma \over \partial \mathbf{x}_\Sigma} \right| \left\| {\partial \mathbf{x}_\Gamma \over \partial \mathbf{x}_\Sigma} \mathbf{n}_\Gamma^{\left( d_f \right)} \right\|_2^{-1} \,\mathrm{d}\Sigma = 0, ~ \forall \tilde{\gamma}_{fa} \in \mathcal{H}\left(\Sigma\right)
\end{split}
\end{equation}
and
\begin{equation}\label{equ:AdjEquAreaConstr1Dm}
\begin{split}
  & \mathrm{Find} ~ d_{fa} \in \mathcal{H}\left( \Sigma \right),~ \mathrm{such~that} \\
  & \int_\Sigma r_f^2 \left( \nabla_\Gamma^{\left( d_f, \tilde{d}_{fa} \right)} \gamma_f \cdot \nabla_\Gamma^{\left( d_f \right)} \gamma_{fa} + \nabla_\Gamma^{\left( d_f \right)} \gamma_f \cdot \nabla_\Gamma^{\left( d_f, \tilde{d}_{fa} \right)} \gamma_{fa} \right) \left| {\partial \mathbf{x}_\Gamma \over \partial \mathbf{x}_\Sigma} \right| \left\| {\partial \mathbf{x}_\Gamma \over \partial \mathbf{x}_\Sigma} \mathbf{n}_\Gamma^{\left( d_f \right)} \right\|_2^{-1} \\
  & + \left( \gamma_p + r_f^2 \nabla_\Gamma^{\left( d_f \right)} \gamma_f \cdot \nabla_\Gamma^{\left( d_f \right)} \gamma_{fa} + \gamma_f \gamma_{fa} - \gamma \gamma_{fa} \right) \Bigg( { \partial \left| {\partial \mathbf{x}_\Gamma \over \partial \mathbf{x}_\Sigma} \right| \left\| {\partial \mathbf{x}_\Gamma \over \partial \mathbf{x}_\Sigma} \mathbf{n}_\Gamma^{\left( d_f \right)} \right\|_2^{-1} \over \partial d_f} \tilde{d}_{fa} \\
  & + { \partial \left| {\partial \mathbf{x}_\Gamma \over \partial \mathbf{x}_\Sigma} \right| \left\| {\partial \mathbf{x}_\Gamma \over \partial \mathbf{x}_\Sigma} \mathbf{n}_\Gamma^{\left( d_f \right)} \right\|_2^{-1} \over \partial \nabla_\Sigma d_f} \cdot \nabla_\Sigma \tilde{d}_{fa} \Bigg) + r_m^2 \nabla_\Sigma d_{fa} \cdot \nabla_\Sigma \tilde{d}_{fa} + d_{fa} \tilde{d}_{fa} \,\mathrm{d}\Sigma = 0, ~ \forall \tilde{d}_{fa} \in \mathcal{H}\left(\Sigma\right).
\end{split}
\end{equation}
The adjoint sensitivity $\delta \left| \Gamma \right|$ in Eq. \ref{equ:AdjSensitivityGaDmAreaConstr} can be derived based on the adjoint analysis of $ \left| \Gamma \right| = \int_\Gamma 1 \,\mathrm{d}\Gamma$:
\begin{equation}\label{equ:AdjSensAreaConstr11}
\begin{split}
  \delta \left| \Gamma \right| = & \int_\Sigma - \gamma_{fa} \tilde{\gamma} \left| {\partial \mathbf{x}_\Gamma \over \partial \mathbf{x}_\Sigma} \right| \left\| {\partial \mathbf{x}_\Gamma \over \partial \mathbf{x}_\Sigma} \mathbf{n}_\Gamma^{\left( d_f \right)} \right\|_2^{-1} - A_d d_{fa} \tilde{d}_m \,\mathrm{d}\Sigma, ~ \forall \left( \tilde{\gamma}, \tilde{d}_m \right) \in \left(\mathcal{L}^2\left(\Sigma\right)\right)^2.
\end{split}
\end{equation}
In Eq. \ref{equ:AdjSensAreaConstr11}, the adjoint variables $\gamma_{fa}$ and $d_{fa}$ are derived by solving the variational formulations for the adjoint equations of the surface-PDE filters for $\gamma$ and $d_m$, respectively:
\begin{equation}\label{equ:AdjEquAreaConstr11Ga}
\begin{split}
  & \mathrm{Find} ~ \gamma_{fa} \in \mathcal{H}\left( \Sigma \right),~ \mathrm{such~that} \\
  & \int_\Sigma \left( r_f^2 \nabla_\Gamma^{\left(d_f\right)} \gamma_{fa} \cdot \nabla_\Gamma^{\left(d_f\right)} \tilde{\gamma}_{fa} + \gamma_{fa} \tilde{\gamma}_{fa} \right) \left| {\partial \mathbf{x}_\Gamma \over \partial \mathbf{x}_\Sigma} \right| \left\| {\partial \mathbf{x}_\Gamma \over \partial \mathbf{x}_\Sigma} \mathbf{n}_\Gamma^{\left( d_f \right)} \right\|_2^{-1} \,\mathrm{d}\Sigma = 0, ~ \forall \tilde{\gamma}_{fa} \in \mathcal{H}\left(\Sigma\right)
\end{split}
\end{equation}
and
\begin{equation}\label{equ:AdjEquAreaConstr11Dm}
\begin{split}
  & \mathrm{Find} ~ d_{fa} \in \mathcal{H}\left( \Sigma \right),~ \mathrm{such~that} \\
  & \int_\Sigma r_f^2 \left( \nabla_\Gamma^{\left(d_f, \tilde{d}_{fa} \right)} \gamma_f \cdot \nabla_\Gamma^{\left(d_f\right)} \gamma_{fa} + \nabla_\Gamma^{\left(d_f\right)} \gamma_f \cdot \nabla_\Gamma^{\left(d_f, \tilde{d}_{fa} \right)} \gamma_{fa} \right) \left| {\partial \mathbf{x}_\Gamma \over \partial \mathbf{x}_\Sigma} \right| \left\| {\partial \mathbf{x}_\Gamma \over \partial \mathbf{x}_\Sigma} \mathbf{n}_\Gamma^{\left( d_f \right)} \right\|_2^{-1} \\
  & + \left( 1 + r_f^2 \nabla_\Gamma^{\left(d_f\right)} \gamma_f \cdot \nabla_\Gamma^{\left(d_f\right)} \gamma_{fa} + \gamma_f \gamma_{fa} - \gamma \gamma_{fa} \right) \Bigg( { \partial \left| {\partial \mathbf{x}_\Gamma \over \partial \mathbf{x}_\Sigma} \right| \left\| {\partial \mathbf{x}_\Gamma \over \partial \mathbf{x}_\Sigma} \mathbf{n}_\Gamma^{\left( d_f \right)} \right\|_2^{-1} \over \partial d_f} \tilde{d}_{fa} \\
  & + { \partial \left| {\partial \mathbf{x}_\Gamma \over \partial \mathbf{x}_\Sigma} \right| \left\| {\partial \mathbf{x}_\Gamma \over \partial \mathbf{x}_\Sigma} \mathbf{n}_\Gamma^{\left( d_f \right)} \right\|_2^{-1} \over \partial \nabla_\Sigma d_f} \cdot \nabla_\Sigma \tilde{d}_{fa} \Bigg) + r_m^2 \nabla_\Sigma d_{fa} \cdot \nabla_\Sigma \tilde{d}_{fa} + d_{fa} \tilde{d}_{fa} \,\mathrm{d}\Sigma = 0, ~ \forall \tilde{d}_{fa} \in \mathcal{H}\left(\Sigma\right).
\end{split}
\end{equation}

For the volume constraint, the adjoint sensitivity of the volume $v$ is derived as
\begin{equation}\label{equ:AdjSensitivityGaDmVolConstr}
\begin{split}
\delta v = \int_\Sigma - A_d d_{fa} \tilde{d}_m \,\mathrm{d}\Sigma, ~ \forall \tilde{d}_m \in \mathcal{L}^2\left(\Sigma\right).
\end{split}
\end{equation}
In Eq. \ref{equ:AdjSensitivityGaDmVolConstr}, the adjoint variable $d_{fa}$ is derived by solving the variational formulation for the adjoint equation of the surface-PDE filter for $d_m$:
\begin{equation}\label{equ:AdjEquVolConstrDm}
\begin{split}
  & \mathrm{Find} ~ d_{fa} \in \mathcal{H}\left( \Sigma \right),~ \mathrm{such~that} \\
  & \int_\Sigma {1 \over \left| \Sigma \right|} \tilde{d}_{fa} + r_m^2 \nabla_\Sigma d_{fa} \cdot  \nabla_\Sigma \tilde{d}_{fa} + d_{fa} \tilde{d}_{fa} \,\mathrm{d}\Sigma = 0, ~ \forall \tilde{d}_{fa} \in \mathcal{H}\left(\Sigma\right).
\end{split}
\end{equation}

After the derivation of the adjoint sensitivities in Eqs. \ref{equ:AdjSensitivityGaDm}, \ref{equ:AdjSensitivityGaDmAreaConstr} and \ref{equ:AdjSensitivityGaDmVolConstr}, the design variables $\gamma$ and $d_m$ can be evolved iteratively to determine the fiber bundle of the surface flow.

\section{Numerical implementation} \label{sec:NumericalImplementationSurfaceNSEqus}

The fiber bundle topology optimization problem in Eq. \ref{equ:VarProToopSurfaceFlows} is solved by using an iterative procedure described as the pseudocode in Tab. \ref{tab:IterativeProcedureSurfaceFlow}, where a loop is included for the iterative solution. The surface finite element method is utilized to solve the variational formulations of the relevant PDEs and adjoint equations.
On the details for the surface finite element solution, one can refer to \cite{DziukActaNumerica2013}. Especially, when the surface finite element method is used to solve the surface flow problems on the implicit 2-manifold filled with the porous medium, the Lagrange multiplier method is used to enforce the tangential constraints of the fluid velocity \cite{FriesIJNMF2018,ReutherPOF2018}. To avoid the numerical singularity caused by the null value of the denominator, the 2-norm of a vector function is approximated in the numerical implementation as $\left\| \mathbf{f} \right\|_2 \rightarrow \left( \mathbf{f}^2 + \epsilon_0 \right)^{1/2} $, where $\mathbf{f}$ is the vector function and $\epsilon_0$ is the value of floating point precision.

\begin{figure}[!htbp]
  \centering
  \includegraphics[width=0.8\textwidth]{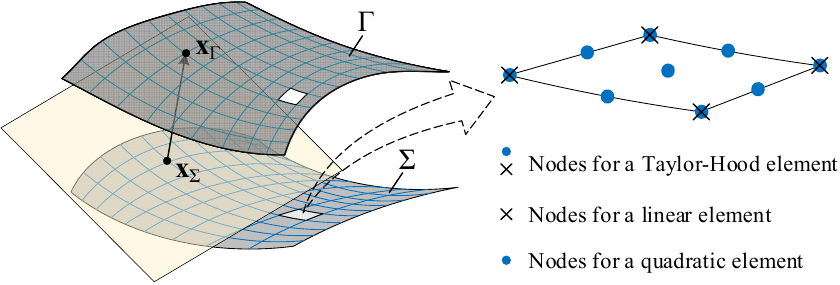}
  \caption{Sketch for the meshes of the Taylor-Hood, linear and quadratic elements of the quadrangular-element based discretization of the base manifold $\Sigma$ and the mapping meshes on the implicit 2-manifold $\Gamma$.}\label{fig:ElementNodes}
\end{figure}

To ensure the well-posedness of the numerical solution of the variational formulations of the surface Navier-Stokes equations and their adjoint equations (Eqs. \ref{equ:VariationalFormulationSurfaceNavierStokesEqus} and \ref{equ:AdjSurfaceNavierStokesEqusJObjective}), the Taylor-Hood elements satisfying the inf-sup condition are used \cite{ElmanFEMFlow2006}. The linear elements are used to interpolate the design variable of the pattern of the surface flow and solve the variational formulation of the surface-PDE filter for this design variable and the corresponding adjoint equation. The quadratic elements are used to interpolate the design variable for the implicit 2-manifold and solve the variational formulation of the surface-PDE filter for this design variable and the corresponding adjoint equation. The meshes of the Taylor-Hood, linear and quadratic elements of the quadrangular-element based discretization of the base manifold have been sketched in Fig. \ref{fig:ElementNodes}, including the mapping meshes on the implicit 2-manifold.

\begin{table}[!htbp]
\centering
\begin{tabular}{l}
  \hline
  \textbf{Algorithm}: iterative solution of Eq. \ref{equ:VarProToopSurfaceFlows} \\
  \hline
  Set $\mathbf{u}_{l_v}$ $p_0$, $\rho$, $\eta$, $A_d$, $v_0$, and $s_0$;\\
  Set $\left\{
  \begin{array}{l}
    \gamma \leftarrow s_0 \\
    d_m \leftarrow v_0+1/2
  \end{array}
  \right.$, $\left\{
  \begin{array}{l}
    r_f = 1/50 \\
    r_m = 2/25 \\
  \end{array}
  \right.$, $\left\{
  \begin{array}{l}
    n_{\max} \leftarrow 240 \\
    n_i \leftarrow 1
  \end{array}
  \right.$, \\
  ~~~~~~ $\left\{
  \begin{array}{l}
    \xi \leftarrow 0.5 \\
    \beta \leftarrow 1
  \end{array}
  \right.$, $\left\{
  \begin{array}{l}
    \alpha_{\min} \leftarrow 0 \\
    \alpha_{\max} \leftarrow 10^4 \rho \\
    q \leftarrow 1
  \end{array}
  \right.$; \\
  \hline
  \textbf{loop} \\
          \hspace{1em} Solve $d_f$ from Eq. \ref{equ:PDEFilterzmBaseStructure1} and compute $v$; \\
          \hspace{1em} Solve $\gamma_f$ from Eq. \ref{equ:VariationalFormulationPDEFilter}; \\
          \hspace{1em} Project $\gamma_f$ to derive $\gamma_p$ and compute $s$; \\
          \hspace{1em} Solve $\mathbf{u}$, $p$ and $\lambda$ from Eq. \ref{equ:VariationalFormulationSurfaceNavierStokesEqus}, and evaluate $J / J_0$; \\
          \hspace{1em} Solve $\mathbf{u}_a$, $p_a$ and $\lambda_a$ from Eq. \ref{equ:AdjSurfaceNavierStokesEqusJObjective}; \\
          \hspace{1em} Solve $\gamma_{fa}$ and $d_{fa}$ from Eqs. \ref{equ:AdjPDEFilterJObjectiveGa} and \ref{equ:AdjPDEFilterJObjectiveDm}; \\
          \hspace{1em} Evaluate $\delta J $ from Eq. \ref{equ:AdjSensitivityGaDm}; \\
          \hspace{1em} Solve $\gamma_{fa}$ and $d_{fa}$ from Eqs. \ref{equ:AdjEquAreaConstr1Ga} and \ref{equ:AdjEquAreaConstr1Dm}; \\
          \hspace{1em} Evaluate $\delta \left( s \left| \Gamma \right| \right)$ from Eq. \ref{equ:AdjSensAreaConstr1}; \\
          \hspace{1em} Solve $\gamma_{fa}$ and $d_{fa}$ from Eqs. \ref{equ:AdjEquAreaConstr11Ga} and \ref{equ:AdjEquAreaConstr11Dm}; \\
          \hspace{1em} Evaluate $\delta \left| \Gamma \right|$ from Eq. \ref{equ:AdjSensAreaConstr11}; \\
          \hspace{1em} Evaluate $\delta s$ in Eq. \ref{equ:AdjSensitivityGaDmAreaConstr} based on $\delta \left( s \left| \Gamma \right| \right)$ and $\delta \left| \Gamma \right|$; \\
          \hspace{1em} Solve $d_{fa}$ from Eq. \ref{equ:AdjEquVolConstrDm}; \\
          \hspace{1em} Evaluate $\delta v$ from Eq. \ref{equ:AdjSensitivityGaDmVolConstr}; \\
          \hspace{1em} Update $\gamma$ and $d_m$ based on $\delta J$, $\delta s$ and $\delta v$; \\
          \hspace{1em} \textbf{if} $\mod\left(n_i,30\right)==0$ \\
          \hspace{2em} $\beta \leftarrow 2\beta$; \\
          \hspace{1em} \textbf{end} \textbf{if} \\
          \hspace{1em} \textbf{if} $ \left( n_i==n_{\max} \right) $ or $\left\{
          \begin{array}{l}
            \beta == 2^7 \\
            {1\over5}\sum_{m=0}^4 \left| J_{n_i-m} - J_{n_i-\left(m+1\right)} \right|\Big/J_0 \leq 10^{-3} \\
            s \leq s_0 \\
            \left|v-v_0\right| \leq 10^{-3}
          \end{array}
          \right.$ \\
          \hspace{2em} break; \\
          \hspace{1em} \textbf{end} \textbf{if} \\
          \hspace{1em} $n_i \leftarrow n_i+1$ \\
  \textbf{end} \textbf{loop} \\
  \hline
\end{tabular}
\caption{Pseudocode used to solve the fiber bundle topology optimization problem for the surface flow. In the iterative solution loop, $n_i$ is the loop-index; $n_{\max}$ is the maximal value of $n_i$; $J_{n_i-m}$ and $J_{n_i-\left(m+1\right)}$ are the values of $J$ in the $\left(n_i-m\right)$-th and $\left(n_i-\left(m+1\right)\right)$-th iterations; and $\mod$ is the operator used to take the remainder. }\label{tab:IterativeProcedureSurfaceFlow}
\end{table}

In the iterative procedure, the projection parameter $\beta$ with the initial value of $1$ is doubled after every $30$ iterations; the loop is stopped when the maximal iteration number is reached, or if the averaged variation of the design objective in continuous 5 iterations and the residuals of the area and volume constraints are simultaneously satisfied. The design variable is updated by using the method of moving asymptotes \cite{SvanbergIntJNumerMethodsEng1987}.

\section{Results and discussion} \label{sec:NumericalExamplesMatchinOptSurfaceFlow}

In this section, the fiber bundle topology optimization is carried out for the surface flows defined on several different base manifolds, including the flat surfaces for the bending channel and the four-terminal device, the curved surfaces deformed from a square to a sphere and the ones deformed from a cylinder to a M\"{o}bius.

The design objective is set to be the combination of the power of the viscous dissipation and pressure drop between the inlet and outlet:
\begin{equation}\label{equ:DesignObjOptimalMatchSurfaceFlows}
\begin{split}
  J = & ~ \omega \int_\Gamma {\eta \over 2} \left( \nabla_\Gamma \mathbf{u} + \nabla_\Gamma \mathbf{u}^\mathrm{T} \right) : \left( \nabla_\Gamma \mathbf{u} + \nabla_\Gamma \mathbf{u}^\mathrm{T} \right) + \alpha \mathbf{u}^2 \, \mathrm{d}\Gamma \\
  & ~ + \left( 1- \omega \right) \left( \int_{l_{v,\Gamma} \setminus l_{v0,\Gamma}} p \, \mathrm{d} l_{\partial\Gamma} - \int_{l_{s,\Gamma}} p \, \mathrm{d} l_{\partial\Gamma} \right),
\end{split}
\end{equation}
where $\omega$ is the weight of the viscous dissipation and it is valued to be $9/10$ and the weight of the pressure drop is hence $1/10$. The density and dynamic viscosity of the fluid are assigned to be unitary. The surface flows are driven by the boundary velocity at the inlets, in the forms of parabolic distribution as the functions of arc length and with the magnitude set to be $U_0$ with $U_0 = \sup_{\forall\mathbf{x}\in l_{v,\Gamma}}{\left\|\mathbf{u}_{l_{v,\Gamma}}\right\|_2}$. The outlets are set to be open boundaries. The remained boundaries are in the type of no slip.

\subsection{Bending channel} \label{subsec:BendingSurfaceFlowsOptimalMatching}

\begin{figure}[!htbp]
  \centering
  \includegraphics[width=0.8\textwidth]{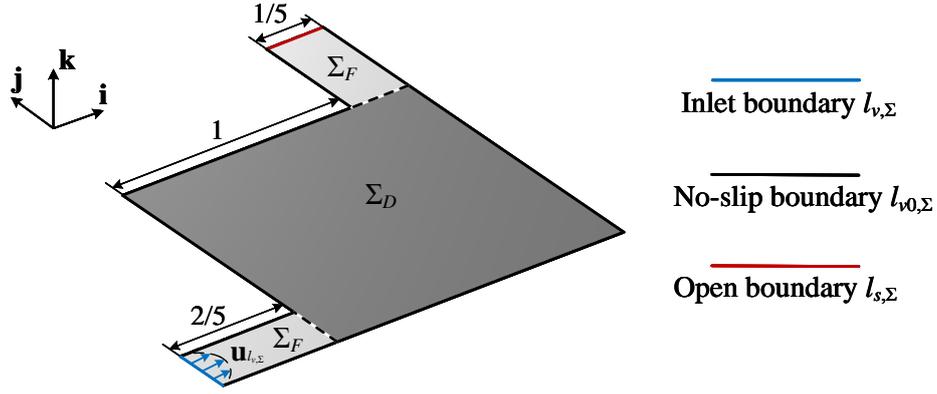}
  \caption{Sketch for the base manifold for the fiber bundle topology optimization of the bending channel, where $\Sigma_D$ is the design domain and $\Sigma_F$ is the channel domains.}\label{fig:BendingChannelDesignDom}
\end{figure}

\begin{figure}[!htbp]
  \centering
  \includegraphics[width=1\textwidth]{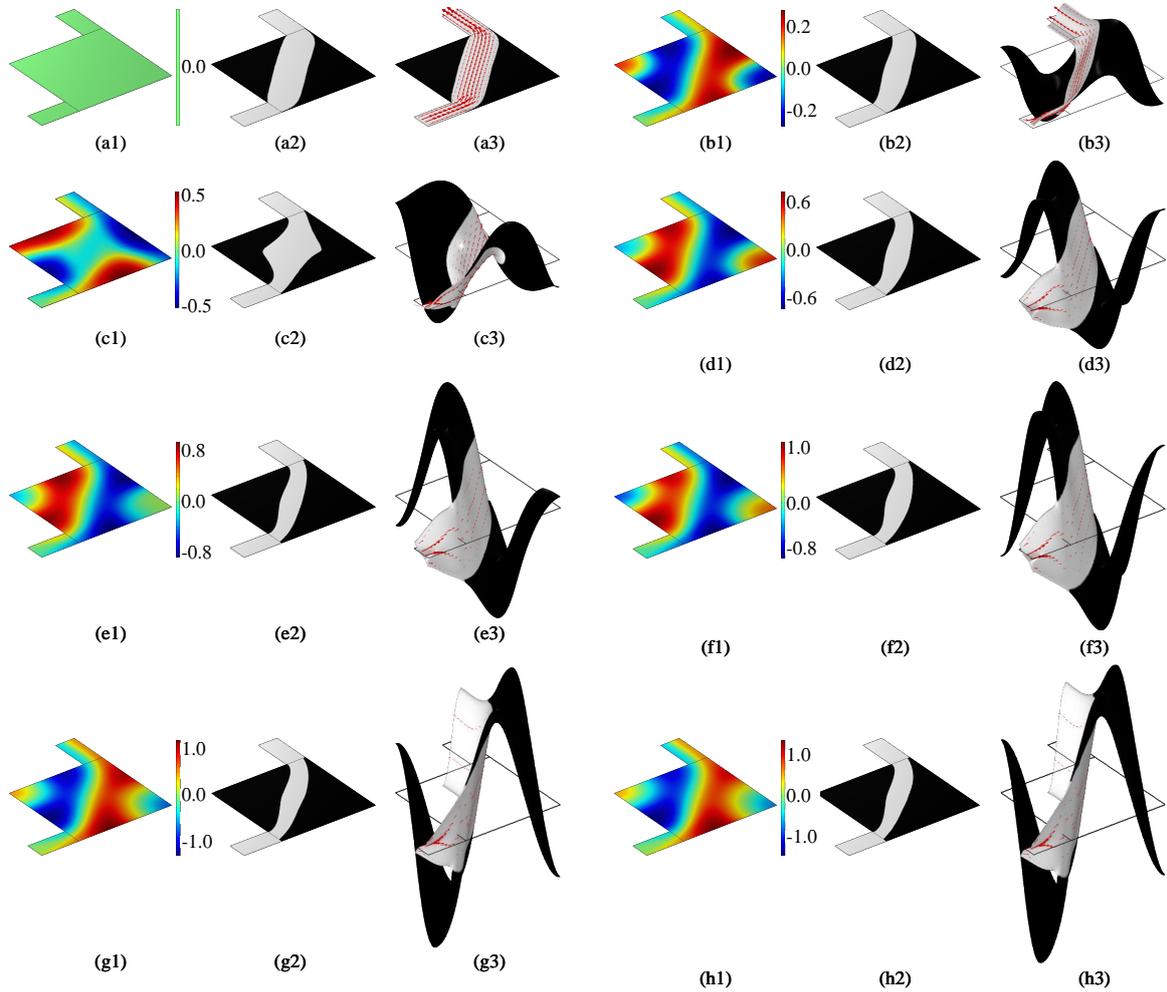}
  \caption{Fiber bundle topology optimization of the bending channel with different values of the magnitude parameter $A_d$, including the distribution of the filtered design variables, the pattern of the bending flow projected on the base manifold and the fiber bundle composed of the base manifold together with the implicit 2-manifold and pattern of the surface flow, where the red arrows represent the distribution of the fluid velocity. (a1)$\sim$(a3) are the results for $A_d = 0$; (b1)$\sim$(b3) are the results for $A_d = 1$; (c1)$\sim$(c3) are the results for $A_d = 2$; (d1)$\sim$(d3) are the results for $A_d = 3$; (e1)$\sim$(e3) are the results for $A_d = 4$; (f1)$\sim$(f3) are the results for $A_d = 5$; (g1)$\sim$(g3) are the results for $A_d = 6$; (h1)$\sim$(h3) are the results for $A_d = 7$.}\label{fig:BendingChannelAdTest}
\end{figure}

For the bending channel, the fiber bundle topology optimization is implemented on the flat surface $\Sigma$ composed of the design domain $\Sigma_D$ and the fluid domain $\Sigma_F$ as shown in Fig. \ref{fig:BendingChannelDesignDom}.
For different values of the magnitude parameter $A_d$, the optimized fiber bundles and their components are derived as shown in Fig. \ref{fig:BendingChannelAdTest}(a$\sim$h) including the distribution of the velocity vectors, where the area and volume fractions are set to be $s_0 = 0.3$ and $v_0 = 0$, respectively. Especially, the fiber bundle topology optimization problem degenerates into the topology optimization problem for the bending flow on the flat surface as shown in Fig. \ref{fig:BendingChannelAdTest}(a), when the magnitude parameter $A_d$ is set to be $0$. By setting the volume fraction to be $0$, the implicit 2-manifold $\Gamma$ is derived with the same absolute values of the positive part and negative part of the enclosed volumes at the two sides of $\Sigma$, respectively. This can be confirmed from the distribution of the filtered design variable $d_f$ shown in Fig. \ref{fig:BendingChannelAdTest}(a1, b1, c1, d1, e1, f1, g1 and h1).

The objective values for the optimized results derived in Fig. \ref{fig:BendingChannelAdTest} have been listed in Tab. \ref{tab:BendingChannelAdObjValues}. From Tab. \ref{tab:BendingChannelAdObjValues}, it can be concluded that higher value of $A_d$ is helpful to decrease the viscous dissipation and pressure drop of the bending flow, because the design space of the fiber bundle topology optimization problem in Eq. \ref{equ:VarProToopSurfaceFlows} can be enlarged by increasing the value of the magnitude parameter in Eq. \ref{equ:PDEFilterzmBaseStructure1}.

\begin{table}[!htbp]
\centering
\begin{tabular}{ccccc}
  \toprule
          Fig. \ref{fig:BendingChannelAdTest}a
        & Fig. \ref{fig:BendingChannelAdTest}b
        & Fig. \ref{fig:BendingChannelAdTest}c
        & Fig. \ref{fig:BendingChannelAdTest}d \\
          $ \left( A_d = 0 \right) $ & $ \left( A_d = 1 \right) $ & $ \left( A_d = 2 \right) $ & $ \left( A_d = 3 \right) $ \\
  \midrule
  $3.6022\times10^1$ & $1.6361\times10^1$ & $9.6811\times10^0$ & $5.7958\times10^0$ \\
  \bottomrule
\end{tabular}
\centering
\begin{tabular}{ccccc}
  \toprule
          Fig. \ref{fig:BendingChannelAdTest}e
        & Fig. \ref{fig:BendingChannelAdTest}f
        & Fig. \ref{fig:BendingChannelAdTest}g
        & Fig. \ref{fig:BendingChannelAdTest}h \\
          $ \left( A_d = 4 \right) $ & $ \left( A_d = 5 \right) $ & $ \left( A_d = 6 \right) $ & $ \left( A_d = 7 \right) $ \\
  \midrule
  $4.7057\times10^0$ & $4.0221\times10^0$ & $3.6217\times10^0$ & $3.3652\times10^0$ \\
  \bottomrule
\end{tabular} \\
\caption{Converged values of the design objective for the fiber bundles derived by sequentially setting the amplitude parameter $A_d$ to be the elements of $\left\{0, 1, 2, 3, 4, 5, 6, 7\right\}$.}\label{tab:BendingChannelAdObjValues}
\end{table}

\begin{figure}[!htbp]
  \centering
  \includegraphics[width=1\textwidth]{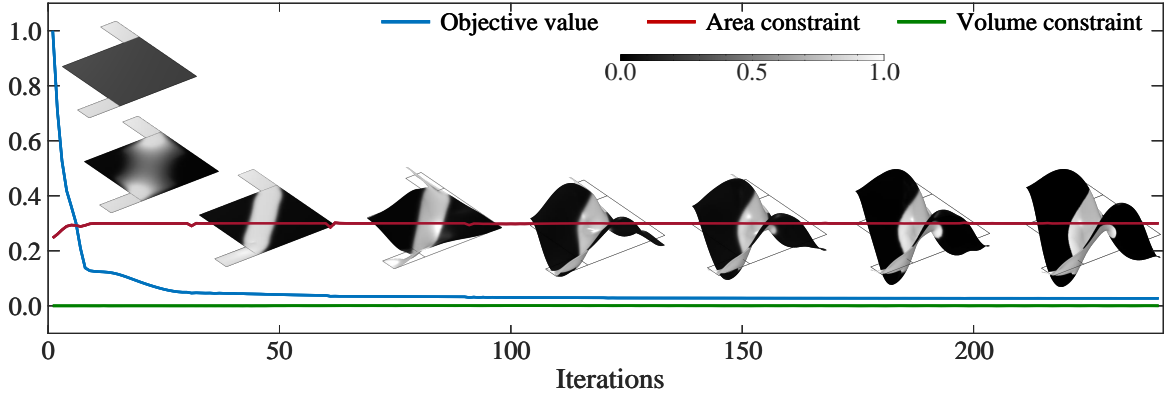}
  \caption{Convergent histories for the result of the bending channel in Fig. \ref{fig:BendingChannelAdTest}c including the snapshots for the evolution of the fiber bundles.}\label{fig:ConvergHistoryBendingAd}
\end{figure}

The convergent histories of the objective values and area and volume constraints have been plotted in Fig. \ref{fig:ConvergHistoryBendingAd} for the results in Fig. \ref{fig:BendingChannelAdTest}c with the magnitude parameter $A_d=2$, including the snapshots for the evolution of the fiber bundles. From the convergent histories, the robust convergence of the numerical solution of the fiber bundle topology optimization problem can be confirmed for the bending flow. In the convergent histories, there are jumps of the objective values and area constraints, and those phenomena are caused by updating the projection parameter $\beta$ in Eq. \ref{equ:ProjectionGammaFilber}. Meanwhile, the convergent histories of the volume constraints are smooth, without jumping phenomenon. This is because that no projection operation is imposed to regularize the design variable of the implicit 2-manifolds.

\begin{figure}[!htbp]
  \centering
  \includegraphics[width=1\textwidth]{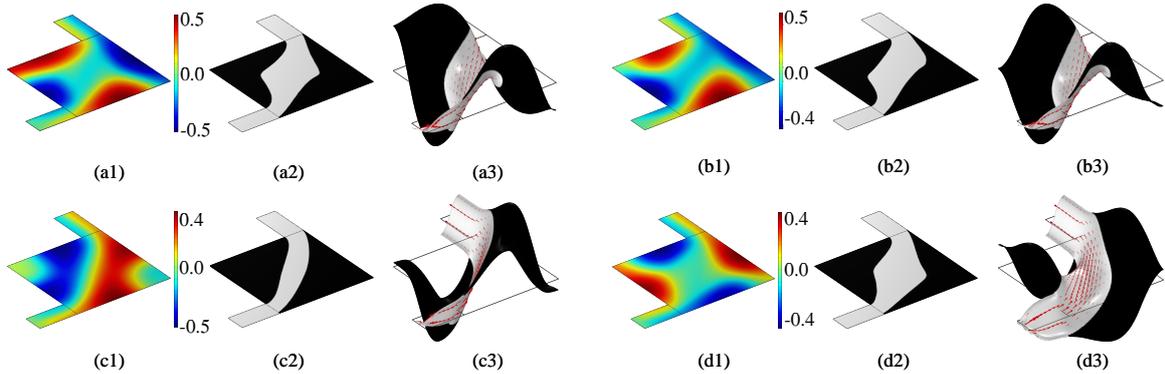}
  \caption{Fiber bundle topology optimization of the bending channel with different values of the velocity magnitude $U_0$ at the inlet, including the distribution of the filtered design variable $d_f$, the pattern of the bending channel $\gamma_p$ projected on $\Sigma$ and the fiber bundle composed of the pattern and 2-manifolds, where the red arrows represent the distribution of the fluid velocity. (a1)$\sim$(a3) are the results for $U_0 = 1\times10^0$; (b1)$\sim$(b3) are the results for $U_0 = 2\times10^2$; (c1)$\sim$(c3) are the results for $U_0 = 5\times10^2$; (d1)$\sim$(d3) are the results for $U_0 = 8\times10^2$.}\label{fig:BendingChannelRe}
\end{figure}

By choosing the magnitude parameter to be $A_d = 2$, the fiber bundle topology optimization problem is further investigated for different values of the velocity magnitude at the inlet of the bending flow. The optimized results are derived as shown in Fig. \ref{fig:BendingChannelRe} by setting $U_0$ to be the elements of $\left\{1\times10^0, 2\times10^2, 5\times10^2, 8\times10^2\right\}$, sequentially. Because larger value of $U_0$ corresponds to stronger Reynolds effect of the surface flow, increasing the velocity magnitude at the inlet can strengthen the convection of the surface flow. Therefore, different fiber bundles are derived as shown in Fig. \ref{fig:BendingChannelRe}.

To confirm the optimality of the derived fiber bundles of the bending flow, the results in Fig. \ref{fig:BendingChannelRe} are cross compared by computing the objective values as listed in Tab. \ref{tab:OptimalityConfirmationBendingChannel}. From the lowest value (marked in bold) of the design objective in every row of Tab. \ref{tab:OptimalityConfirmationBendingChannel}, the optimized performance of the derived fiber bundles can be confirmed.

\begin{table}[!htbp]
\centering
\begin{tabular}{l|cccc}
  \toprule
        & Fig. \ref{fig:BendingChannelRe}a
        & Fig. \ref{fig:BendingChannelRe}b
        & Fig. \ref{fig:BendingChannelRe}c
        & Fig. \ref{fig:BendingChannelRe}d \\
  \midrule
  $ U_0 = 1\times10^0 $ & $\mathbf{9.6811\times10^0}$ & $9.8902\times10^0$ & $9.7068\times10^0$ & $1.2096\times10^1$ \\
  \midrule
  $ U_0 = 2\times10^2 $ & $3.6185\times10^5$ & $\mathbf{3.4643\times10^5}$ & $3.6381\times 10^5$ & $4.2392\times10^5$ \\
  \midrule
  $ U_0 = 5\times10^2 $ & $2.9563\times10^6$ & $2.8859\times10^6$ & $\mathbf{2.4548\times 10^6}$ & $2.7867\times10^6$ \\
  \midrule
  $ U_0 = 8\times10^2 $ & $9.5576\times10^6$ & $9.3152\times10^6$ & $7.6127\times 10^6$ & $\mathbf{7.2915\times10^6}$ \\
  \bottomrule
\end{tabular} \\
\caption{Values of the design objective in Eq. \ref{equ:DesignObjOptimalMatchSurfaceFlows} for the fiber bundles in Fig. \ref{fig:BendingChannelRe}. The optimized entries have been noted in bold.}\label{tab:OptimalityConfirmationBendingChannel}
\end{table}

\subsection{Four-terminal device} \label{subsec:FourTerminalDevicesOptimalMatching}

\begin{figure}[!htbp]
  \centering
  \includegraphics[width=1\textwidth]{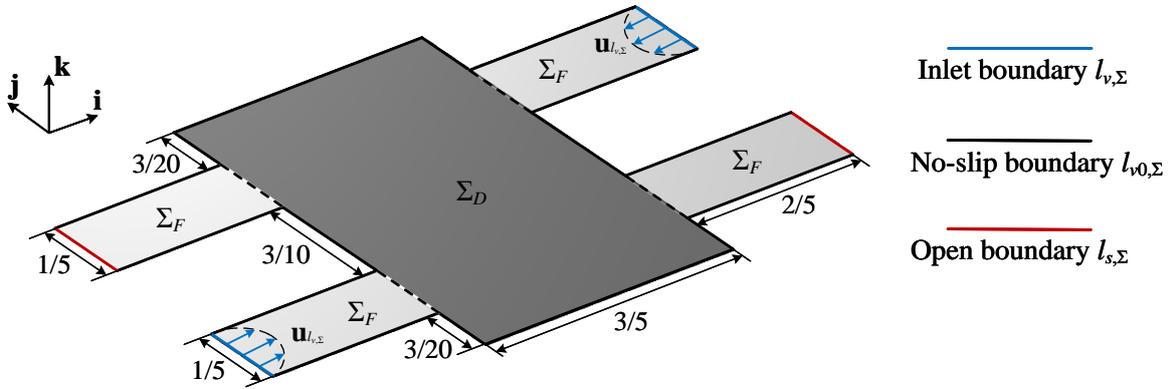}
  \caption{Sketch for the base manifold for the fiber bundle topology optimization of the four-terminal device, where $\Sigma_D$ is the design domain and $\Sigma_F$ is the channel domains.}\label{fig:FourTerminalDeviceDesignDom}
\end{figure}

For the fiber bundle topology optimization problem with the flat surface as its base manifold, the four-terminal device is further investigated by setting the magnitude parameter $A_d$ to be $0$ and $2$, respectively. The computational domain is set as the flat surface shown in Fig. \ref{fig:FourTerminalDeviceDesignDom} composed of the design domain $\Sigma_D$ and the channel domains $\Sigma_F$. By setting the area and volume fractions to be $s_0 = 0.4$ and $v_0 = 0$, the optimized results are derived as listed in Tab. \ref{tab:FourTerminalDevicesResults} for different velocity magnitude at the inlets.

\begin{table}[!htbp]
\centering
\begin{tabular}{l|cccc}
  \toprule
        & $U_0 = 5.50\times10^2$
        & $U_0 = 5.75\times10^2$
        & $U_0 = 6.25\times10^2$ \\
  \midrule
  $ A_d = 0 $ & \includegraphics[width=0.28\textwidth]{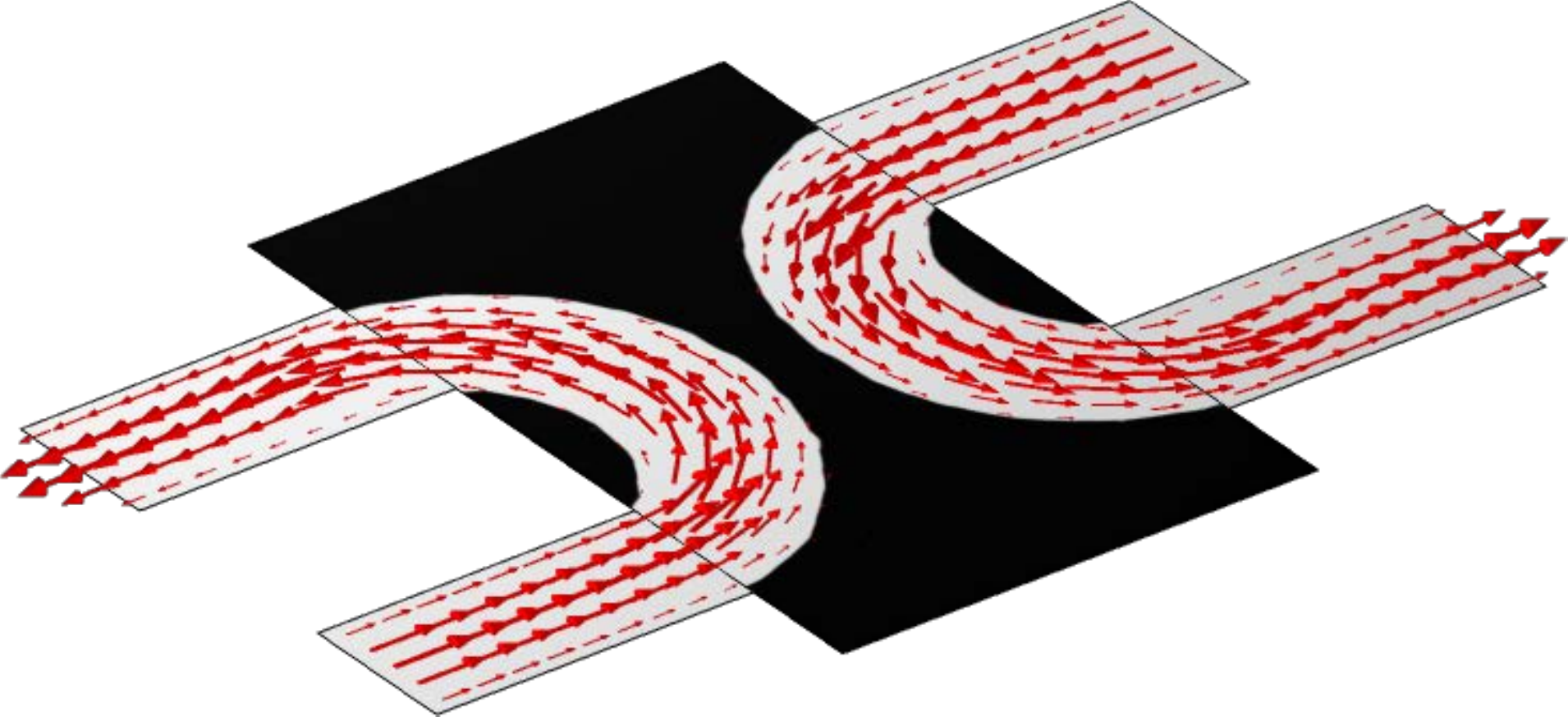} & \includegraphics[width=0.28\textwidth]{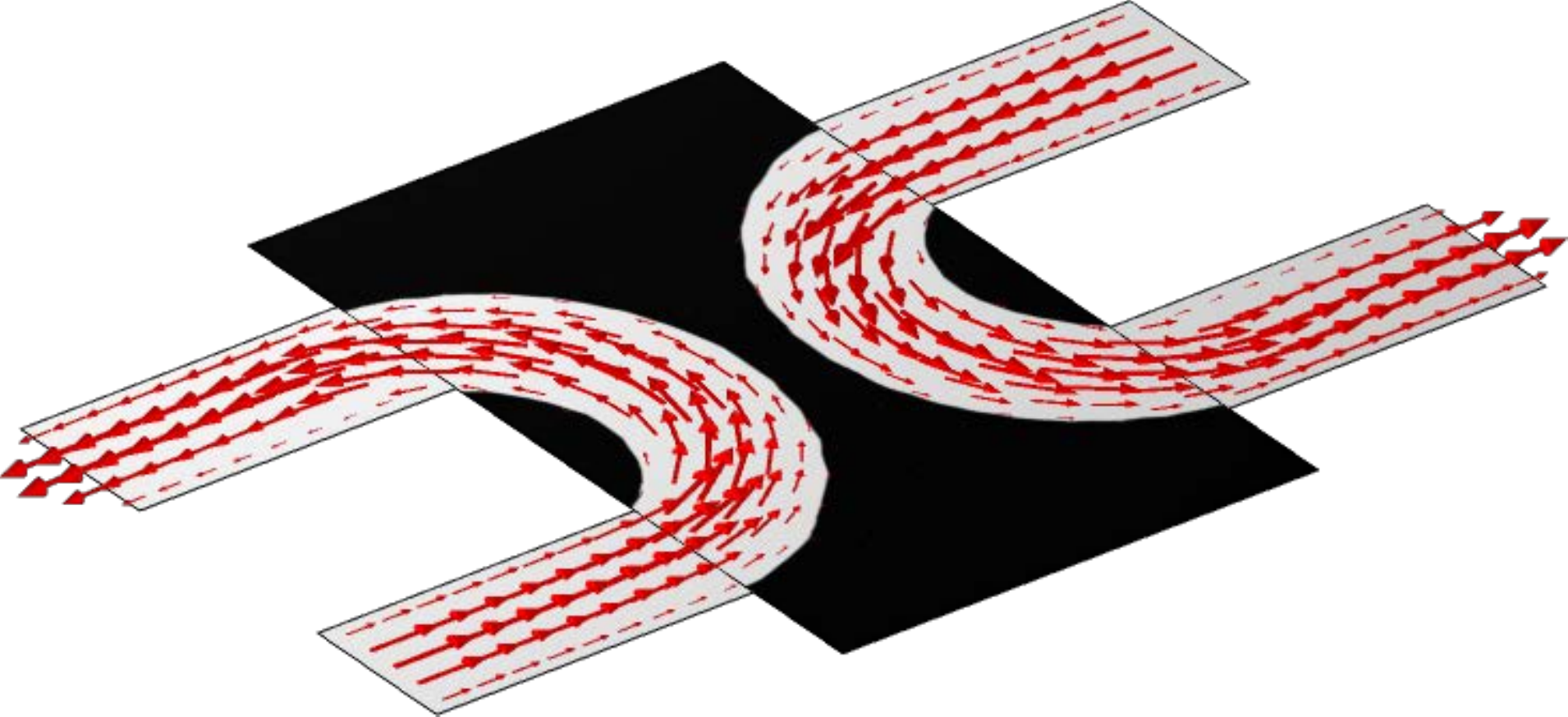} & \includegraphics[width=0.28\textwidth]{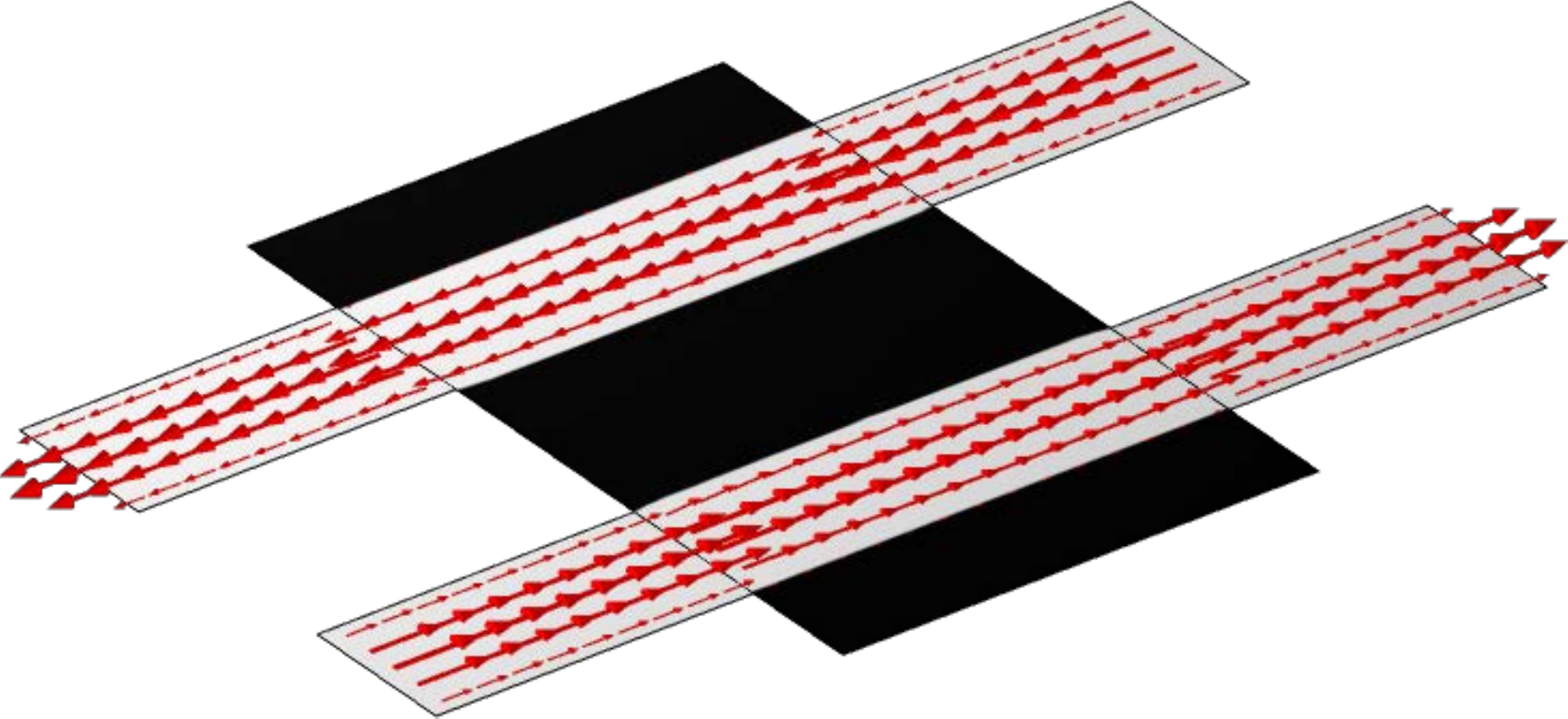} \\
  \midrule
  $ A_d = 2 $ & \includegraphics[width=0.28\textwidth]{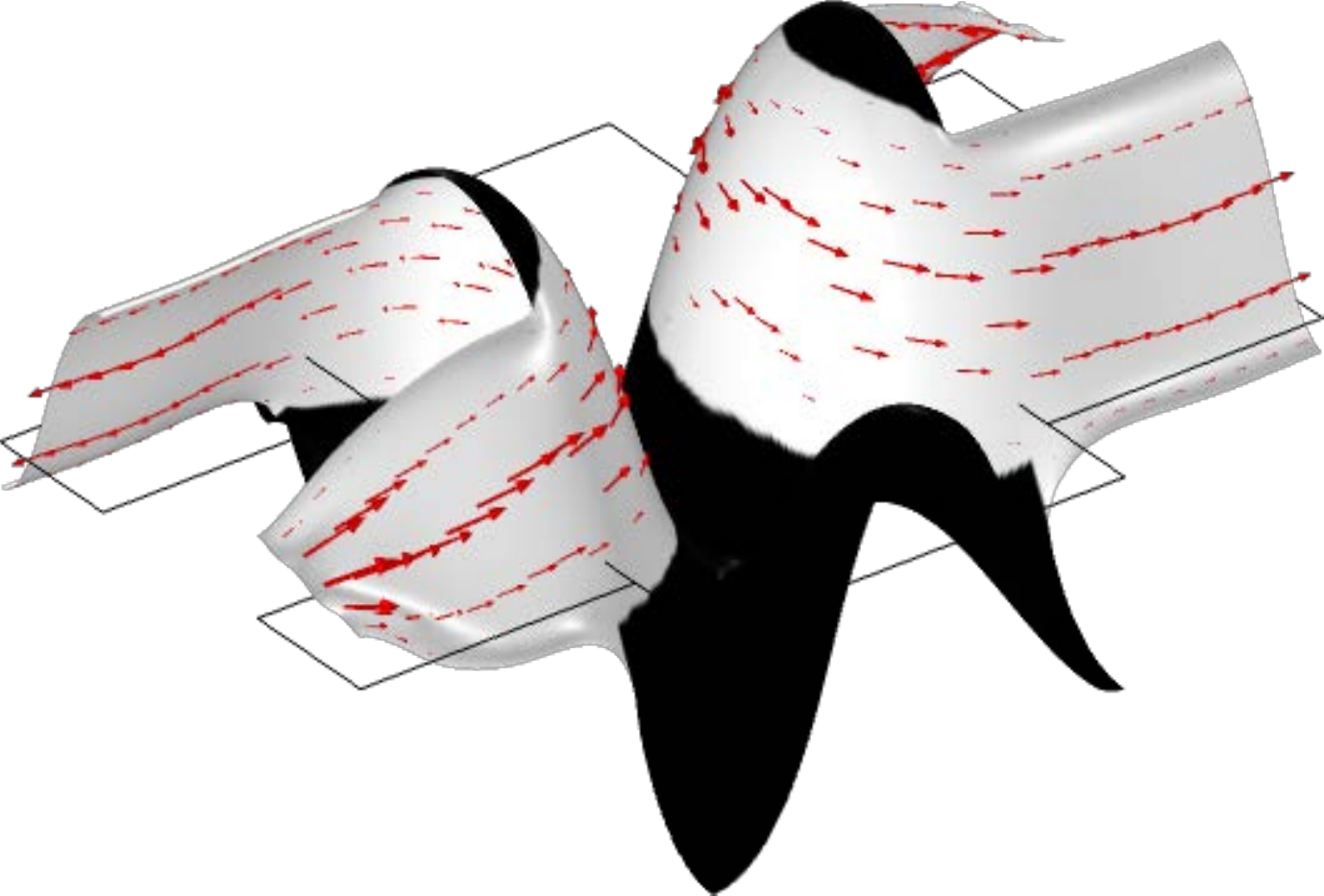} & \includegraphics[width=0.28\textwidth]{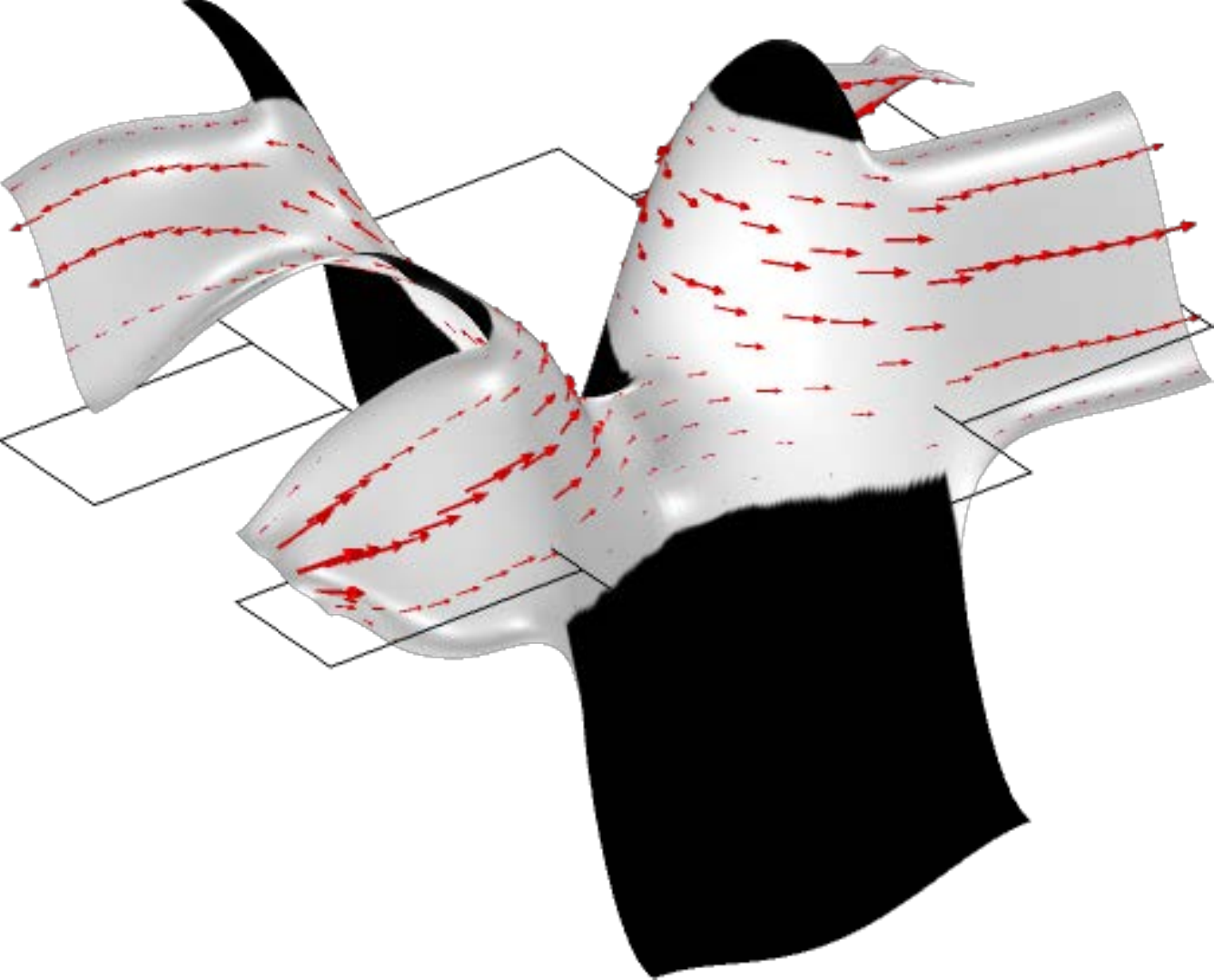} & \includegraphics[width=0.28\textwidth]{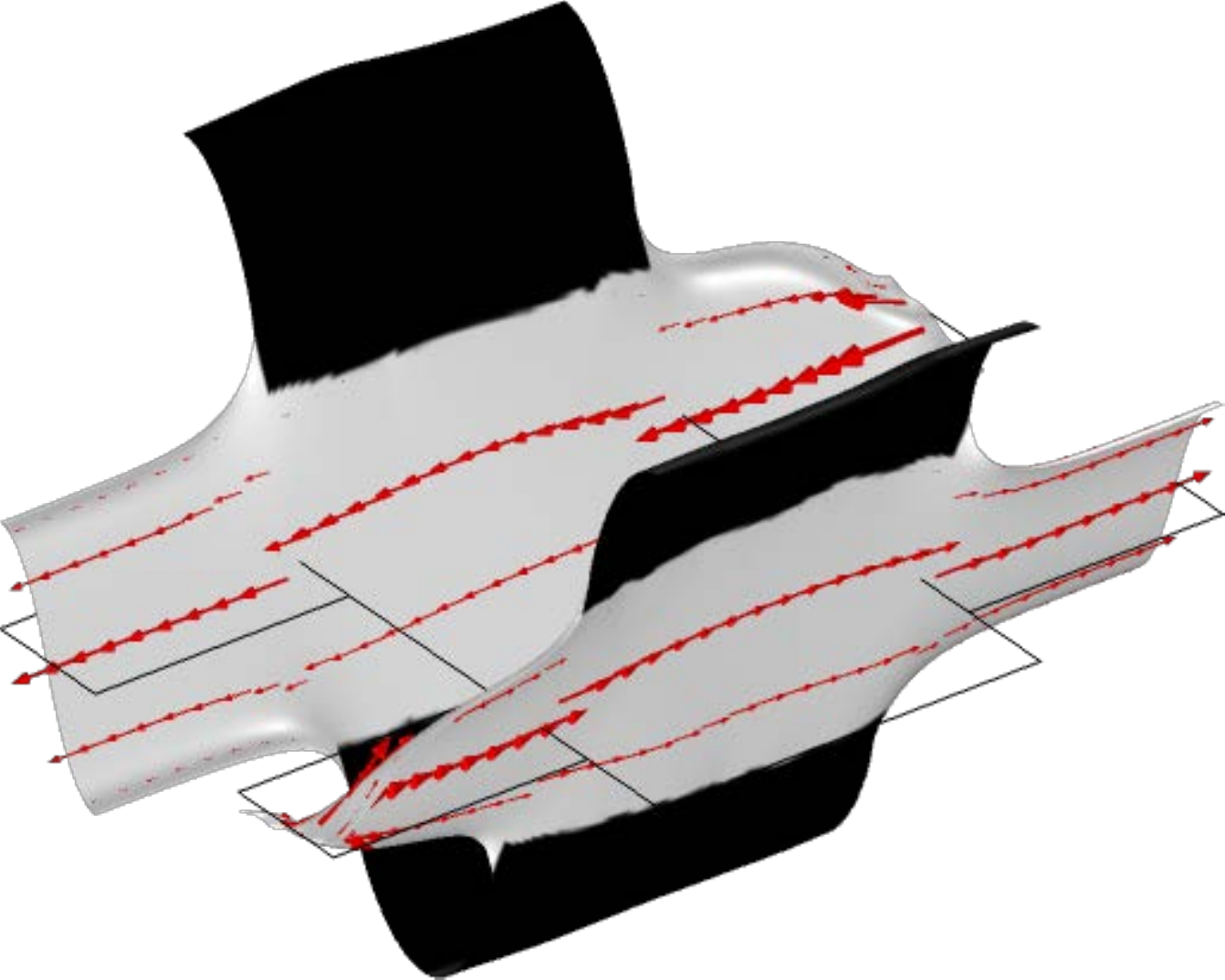} \\
              & \includegraphics[width=0.28\textwidth]{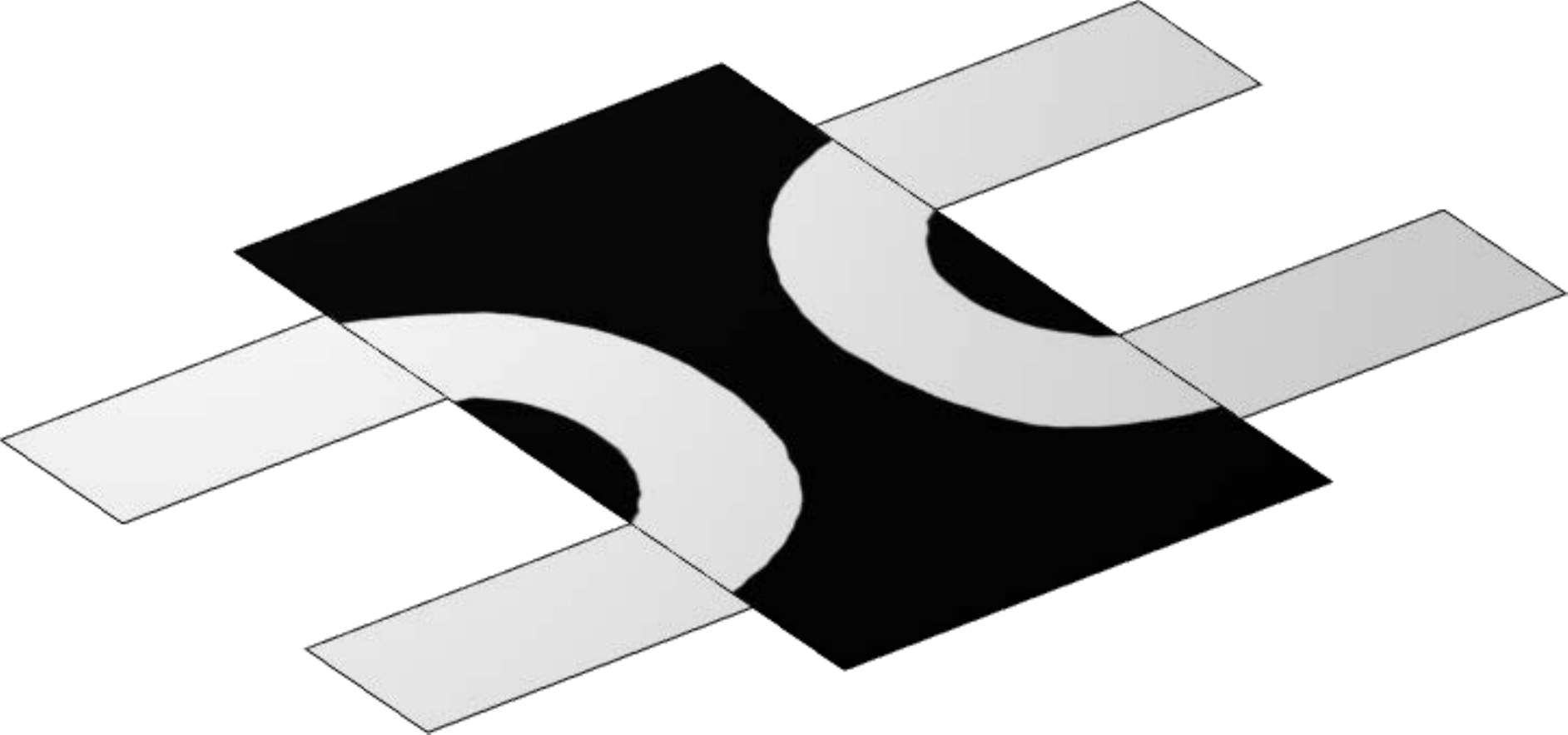} & \includegraphics[width=0.28\textwidth]{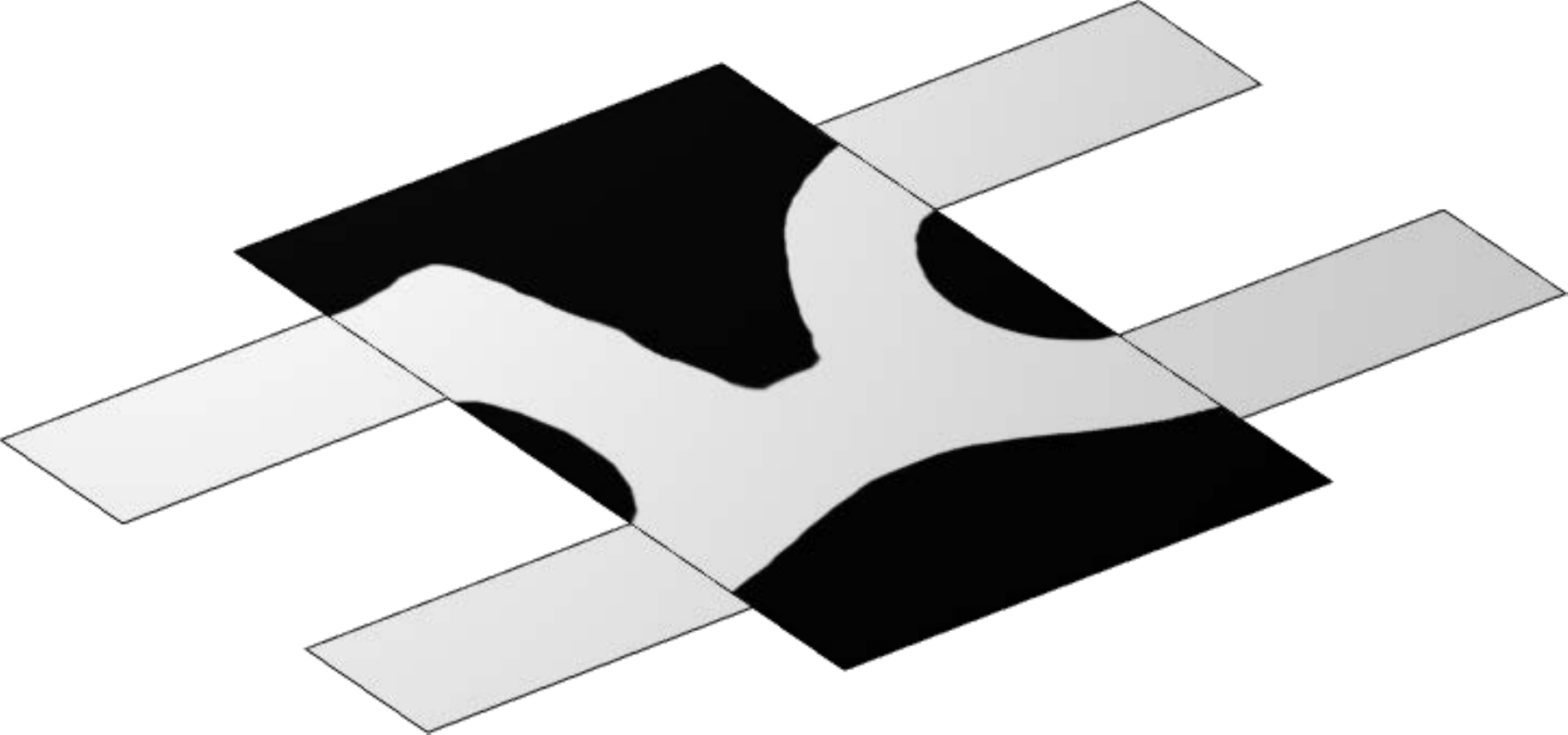} & \includegraphics[width=0.28\textwidth]{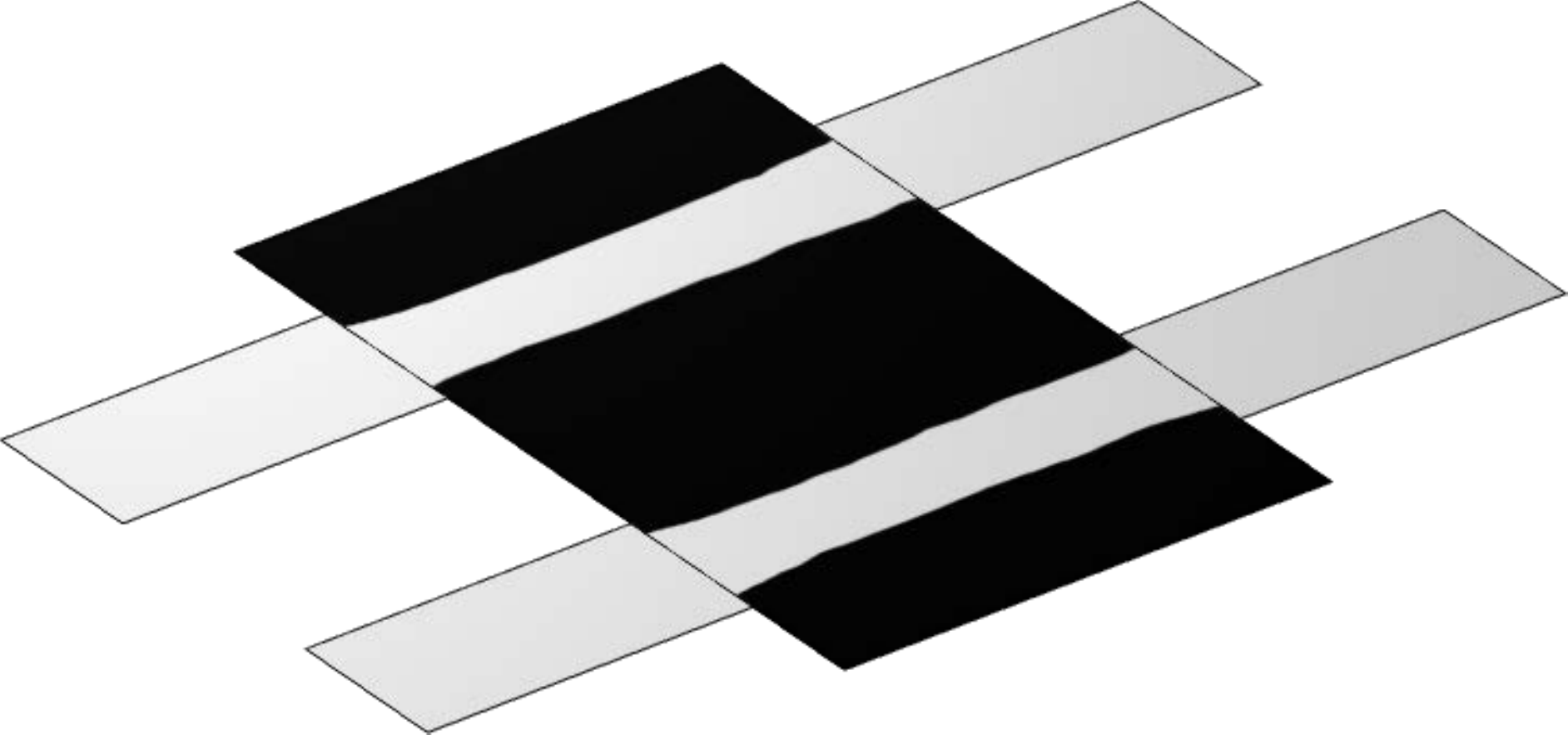} \\
  \bottomrule
\end{tabular} \\
\caption{Fiber bundle topology optimization of the four-terminal device for different values of the velocity magnitude $U_0$ at the inlets, where the red arrows represent the distribution of the fluid velocity.}\label{tab:FourTerminalDevicesResults}
\end{table}

In Tab. \ref{tab:FourTerminalDevicesResults}, the four-terminal device has the topology of double bending channels, when the surface flow has relatively weak Reynolds effect; the topology changes into double straight channels, as the Reynolds effect is strengthened. From the comparison of the results corresponding to $A_d=0$ and $A_d=2$, it can be concluded that the increase of the magnitude parameter can speed up the change of the optimized topology of the four-terminal device, when the Reynolds effect is strengthened. This is because that the increase of the magnitude parameter enlarges the design space of the four-terminal device, i.e., the characteristic size and area of the channels are increased and the averaged velocity is decreased at the inlets, then the gradient of the velocity decreases, and hence the viscous dissipation and pressure drop decrease.

\subsection{Flows on deformed surfaces} \label{subsec:SurfaceDeformedSurfacesOptimalMatching}

\begin{figure}[!htbp]
  \centering
  \includegraphics[width=0.85\textwidth]{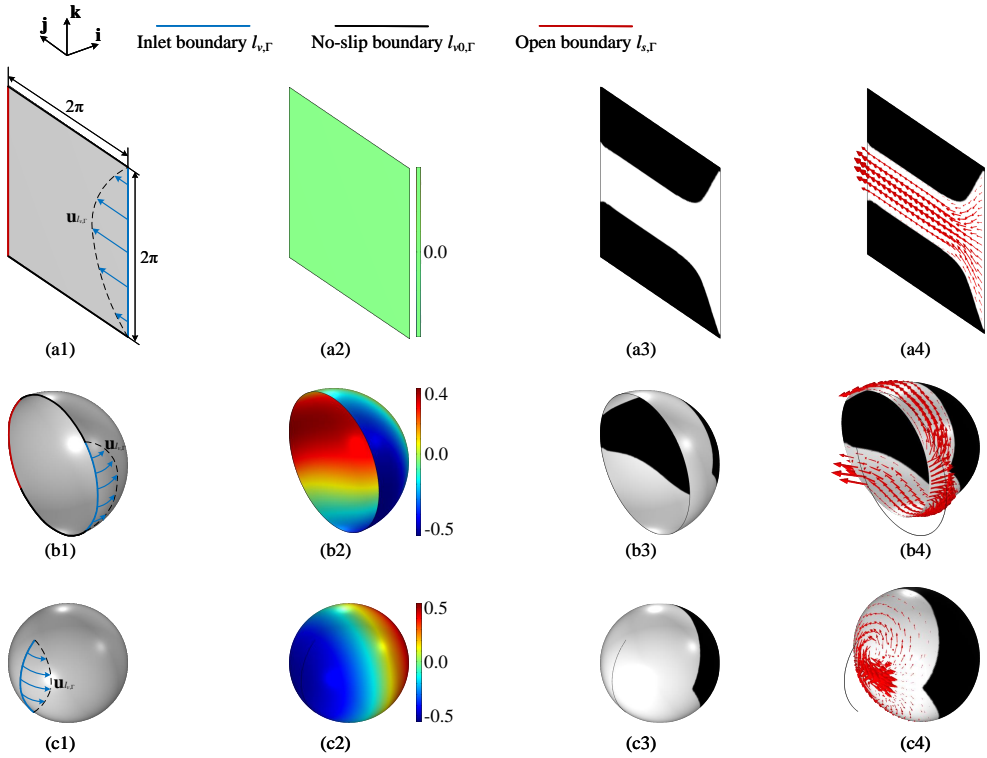}
  \caption{Fiber bundle topology optimization for the surface flows on the base manifolds deformed from a square to a sphere with keeping the conservation of area, where the base manifolds are sketched in Figs. \ref{fig:SquareToSphere}(a1), \ref{fig:SquareToSphere}(b1) and \ref{fig:SquareToSphere}(c1), the distribution of the filtered design variables for the implicit 2-manifolds are shown in Figs. \ref{fig:SquareToSphere}(a2), \ref{fig:SquareToSphere}(b2) and \ref{fig:SquareToSphere}(c2), the projected patterns on the base manifolds are shown in Figs. \ref{fig:SquareToSphere}(a3), \ref{fig:SquareToSphere}(b3) and \ref{fig:SquareToSphere}(c3), and the fiber bundles are derived as shown in Figs. \ref{fig:SquareToSphere}(a4), \ref{fig:SquareToSphere}(b4) and \ref{fig:SquareToSphere}(c4) with the red arrows representing the distribution of the fluid velocity.}\label{fig:SquareToSphere}
\end{figure}

For continuously deformed base manifolds with keeping the conservation of area, the fiber bundle topology optimization is implemented for the surface flows. By deforming a square to a sphere as shown in Fig. \ref{fig:SquareToSphere}(a1$\sim$c1), the optimized fiber bundles are derived as shown in Fig. \ref{fig:SquareToSphere}(a2$\sim$a4), \ref{fig:SquareToSphere}(b2$\sim$b4) and \ref{fig:SquareToSphere}(c2$\sim$c4) including the distribution of the velocity vectors, where the area and volume fractions and the magnitude parameter are set to be $s_0 = 0.4$, $v_0 = 0$ and $A_d = 2$, respectively.

\begin{figure}[!htbp]
  \centering
  \includegraphics[width=0.65\textwidth]{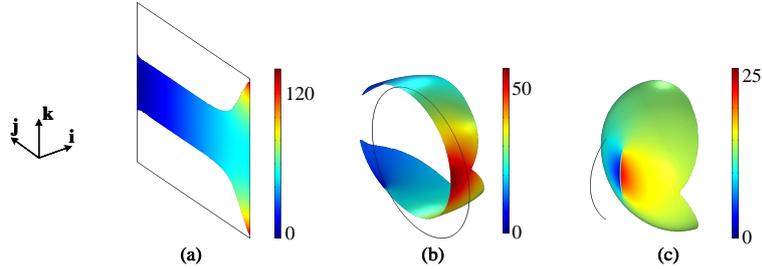}
  \caption{Distribution of the fluid pressure for the surface flows on the fiber bundles derived as shown in Figs. \ref{fig:SquareToSphere}(a4), \ref{fig:SquareToSphere}(b4) and \ref{fig:SquareToSphere}(c4).}\label{fig:PressureSquareToSphere}
\end{figure}

In Fig. \ref{fig:SquareToSphere}, the fiber bundle for the surface flow on the square is composed of the base manifold together with the pattern of the flat diffuser and the implicit 2-manifold coinciding with the base manifold. The flat diffuser is consistent with the previously reported results derived by using topology optimization \cite{Borrvall2003}. When the square deforms into the shape of a semi-sphere, the pattern of the surface flow spits into two branches; and the implicit 2-manifold shrinks to straighten the channels corresponding to the pattern of the surface flow. When the semi-sphere further deforms into a sphere, the two branches merges to remove one part of the no-slip boundary, and the surface flow evolves into the enclosed mode with two vortexes. The underlying mechanism for the evolution of the fiber bundles along with the deformation of base manifolds is that the fluid is prone to moving in the short path and widening the channel, and detachment of the no-slip boundary can help to decrease the viscous dissipation and pressure drop. Additionally, the viscous dissipation and pressure drop of the surface flows with the patterns in the optimized fiber bundles decreases, along with the base manifold deforming from a square to a sphere. This can be conformed from the pressure distribution in Fig. \ref{fig:PressureSquareToSphere}.

\begin{figure}[!htbp]
  \centering
  \includegraphics[width=1\textwidth]{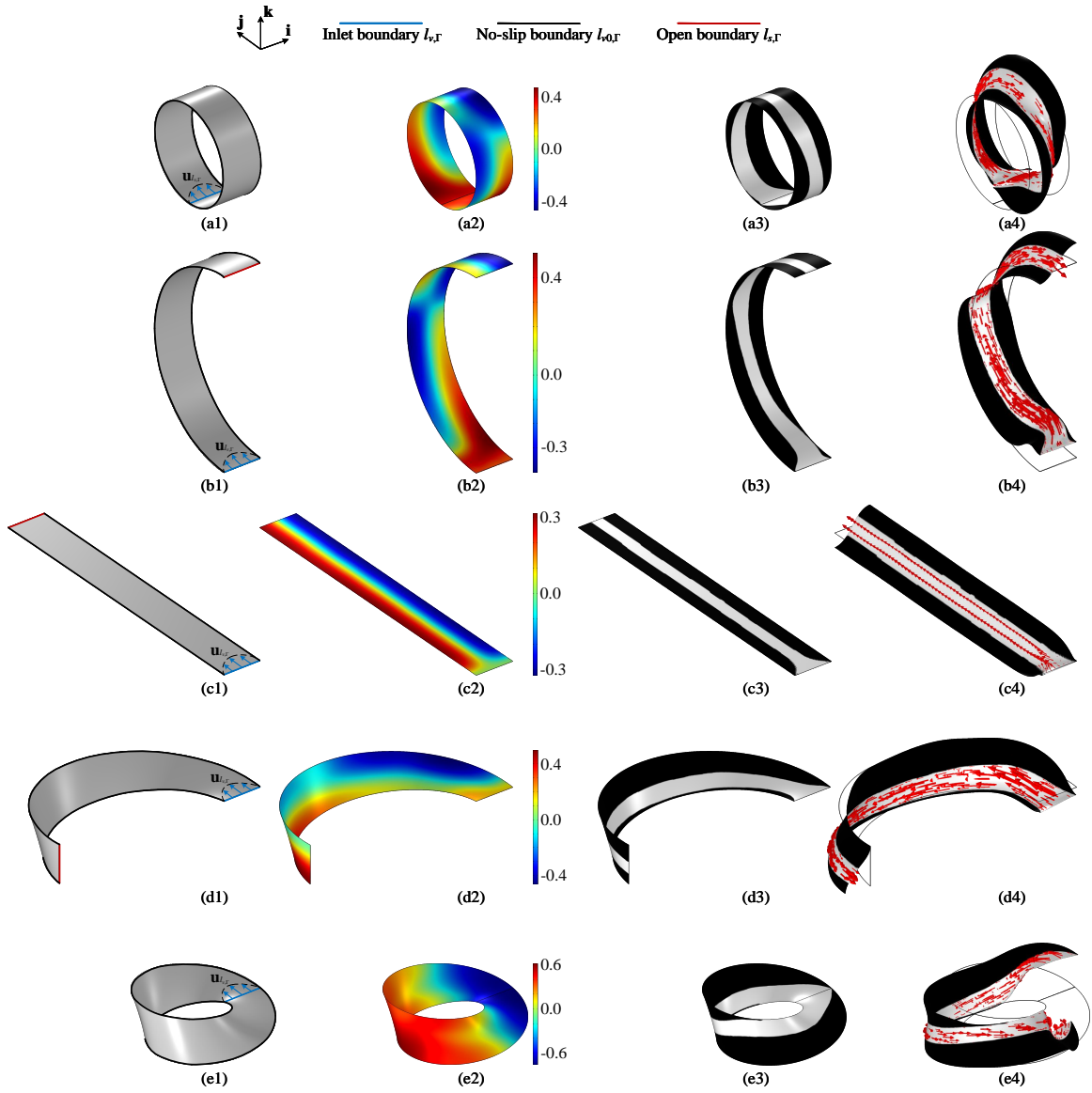}
  \caption{Fiber bundle topology optimization for the surface flows on the base manifolds deformed from a cylinder to a M\"{o}bius with keeping the conservation of area, where the base manifolds are sketched in Figs. \ref{fig:CylinderToMoebius}(a1), \ref{fig:CylinderToMoebius}(b1), \ref{fig:CylinderToMoebius}(c1), \ref{fig:CylinderToMoebius}(d1) and \ref{fig:CylinderToMoebius}(e1), the distribution of the filtered design variables for the implicit 2-manifolds are shown in Figs. \ref{fig:CylinderToMoebius}(a2), \ref{fig:CylinderToMoebius}(b2), \ref{fig:CylinderToMoebius}(c2), \ref{fig:CylinderToMoebius}(d2) and \ref{fig:CylinderToMoebius}(e2), the projected patterns on the base manifolds are shown in Figs. \ref{fig:CylinderToMoebius}(a3), \ref{fig:CylinderToMoebius}(b3), \ref{fig:CylinderToMoebius}(c3), \ref{fig:CylinderToMoebius}(d3) and \ref{fig:CylinderToMoebius}(e3), and the fiber bundles are derived as shown in Figs. \ref{fig:CylinderToMoebius}(a4), \ref{fig:CylinderToMoebius}(b4), \ref{fig:CylinderToMoebius}(c4), \ref{fig:CylinderToMoebius}(d4) and \ref{fig:CylinderToMoebius}(e4) with the red arrows representing the distribution of the fluid velocity.}\label{fig:CylinderToMoebius}
\end{figure}

Further, the fiber bundle topology optimization is implemented on the base manifolds derived by deforming a cylinder firstly to a strip and then to a M\"{o}bius as shown in Fig. \ref{fig:CylinderToMoebius}(a1$\sim$e1) with the sizes marked on the strip, to minimize the viscous dissipation and pressure drop for the surface flows. The area and volume fractions and magnitude parameter are set without change. The patterns of the surface flows and the implicit 2-manifolds of the optimized fiber bundles are derived as shown in Figs. \ref{fig:CylinderToMoebius}(a2$\sim$a4), \ref{fig:CylinderToMoebius}(b2$\sim$b4), \ref{fig:CylinderToMoebius}(c2$\sim$c4), \ref{fig:CylinderToMoebius}(d2$\sim$d4) and \ref{fig:CylinderToMoebius}(e2$\sim$e4), including the distribution of the velocity vectors represented by the red arrows. Especially, the derived implicit 2-manifold on the M\"{o}bius is broken by setting the inlet simultaneously to be the outlet with the same known velocity distribution. The destination of such setting is to make the derived implicit 2-manifold be orientable and remove the singularity of the normal direction on the non-orientable M\"{o}bius.

In Fig. \ref{fig:CylinderToMoebius}, the derived pattern and implicit 2-manifold defined on a cylinder is composed of a circle channel and a curved surface with asymmetry. Then, the cylinder is opened and sequentially evolved into the shapes of semi-cylinder, strip, semi-M\"{o}bius until being enclosed again into the shape of M\"{o}bius, where the conservation of area is kept. The derived patterns and implicit 2-manifolds defined on the base manifolds corresponding to the sequential evolution of the opened cylinder are all composed of channel-shaped patterns and curved surfaces with asymmetry. The asymmetry assists the derived implicit 2-manifolds to satisfy the volume constraint with the volume fraction of $0$, and it is caused by the asymmetrical property of the convection of the surface flows. The asymmetry is advantageous to shorten the path of the surface flows, then to decrease the viscous dissipation and pressure drop.

\begin{figure}[!htbp]
  \centering
  \includegraphics[width=1\textwidth]{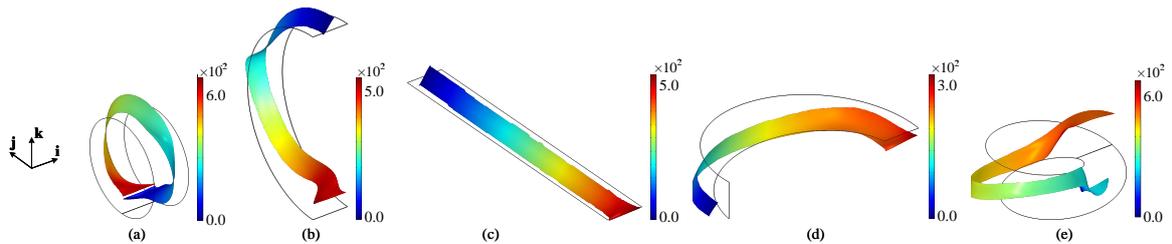}
  \caption{Distribution of the fluid pressure for the surface flows on the fiber bundles derived as shown in Figs. \ref{fig:CylinderToMoebius}(a4), \ref{fig:CylinderToMoebius}(b4), \ref{fig:CylinderToMoebius}(c4), \ref{fig:CylinderToMoebius}(d4) and \ref{fig:CylinderToMoebius}(e4).}\label{fig:PressureCylinderToMoebius}
\end{figure}

The pressure distribution in the derived fiber bundles has been provided in Fig. \ref{fig:PressureCylinderToMoebius}, which shows that the deformation of the base manifold from the cylinder to the M\"{o}bius is advantageous to decrease the viscous dissipation and pressure drop. This is because that the deformation can shorten the path of surface flows with the optimized patterns on the derived implicit 2-manifolds. Additionally, the case of M\"{o}bius has similar pressure drop to that of cylinder. Therefore, the orientability of 2-manifolds can provide similar effectivity on the design domain for fiber bundle topology optimization to minimize the viscous dissipation and pressure drop of a surface flow.

\section{Conclusions}\label{sec:Conclusions}

A fiber bundle topology optimization approach for the surface flow has been developed to match the implicit 2-manifold and the pattern defined on it, where the surface flow is described by the surface Navier-Stokes equations defined on the implicit 2-manifold. The material distribution method is used to implement the topology evolution of the pattern, where an artificial Darcy friction force of the porous model is added to the surface Navier-Stokes equations. The implicit 2-manifold is evolved based on the homeomorphous map between it and the base manifold. Continuous adjoint analysis method has been used to analyze the fiber bundle topology optimization problem.

Numerical tests have been presented to demonstrate this approach, including the fiber bundle topology optimization of bending channel, four-terminal device and fluid channels on continuously deformed base manifolds. The bending channel has been optimized to present the effect of the magnitude parameter used to determine the design space of the implicit 2-manifold, and the results show that the increase of the magnitude parameter can enlarge the design space of the fiber bundle for the surface flow. The Reynolds effect has been demonstrated by the fiber bundle topology optimization of the bending channel and four-terminal device, by setting different velocity magnitude at the inlets. In the results of the bending channel, the valley and slope-shaped implicit 2-manifolds are derived to shorten the fluid path. In the results of the four-terminal device, the magnitude parameter of the implicit 2-manifold can speed up the topology change from double bending channels to double straight channels as the Reynolds effect is strengthened, compared with the results of topology optimization for the surface flow on the flat surface corresponding to the degenerated case with the null value of the magnitude parameter. The fiber bundle topology optimization has also been implemented on the base manifolds derived by deforming a square into a sphere and deforming a cylinder into a M\"{o}bius, where the area conservation is kept during the deformation. The derived results show that the non-orientable base manifolds have similar performance to the orientable ones on shortening the fluid path and minimizing the viscous dissipation and pressure drop.

The presented fiber bundle topology optimization approach includes the design domain into the design space of fluidic structures. This approach achieves the topology optimization for fluid flows on the variable design domain. It provides a topology optimization method for the conformal design of fluidic channels, where the channel topology and the outer shapes of structural walls can be optimized simultaneously to achieve the matching optimization. Especially, the fiber bundle topology optimization problem will degenerate into the topology optimization problem for the fluid flow on a flat surface, if the null value is chosen for the magnitude parameter of the surface-PDE filter and the flat surface is set as the base manifold. Therefore, the presented fiber bundle topology optimization is the generalization of topology optimization for two-dimensional flow problems. This paper focuses on the laminar surface flows. In the future, it can be promoted for the turbulent surface flows.

\section{Acknowledgements}\label{sec:Acknowledgements}

The authors acknowledge the support of the National Natural Science Foundation of China (No. 51875545), the Innovation Grant of Changchun Institute of Optics, Fine Mechanics and Physics (CIOMP), the Youth Innovation Promotion Association of the Chinese Academy of Sciences (No. 2018253) and the Fund of State Key Laboratory of Applied Optics (SKLAO). They are also grateful to Prof. K. Svanberg of KTH for supplying the codes for the method of moving asymptotes.

\section{Appendix}\label{sec:Appendix}

This section provides the details for the adjoint analysis of the fiber bundle topology optimization problem in Eq. \ref{equ:VarProToopSurfaceFlows}.

\subsection{Adjoint analysis for design objective} \label{sec:AdjointAnalysisDesignObjective}

Based on the transformed design objective in Eq. \ref{equ:FurtherTransformedGeneralDesignObjective}, the variational formulations of the surface-PDE filters in Eqs. \ref{equ:VariationalFormulationPDEFilterBaseManifold} and \ref{equ:VariationalFormulationPDEFilter} and the surface Navier-Stokes equations in Eq. \ref{equ:TransformedVariationalFormulationSurfaceNSEqus},  the augmented Lagrangian of the design objective in Eq. \ref{equ:VarProToopSurfaceFlows} can be derived as
\begin{equation}\label{equ:AugmentedLagrangianMatchOptimization}
\begin{split}
  \hat{J} = & \int_\Sigma A \left| {\partial \mathbf{x}_\Gamma \over \partial \mathbf{x}_\Sigma} \right| \left\| {\partial \mathbf{x}_\Gamma \over \partial \mathbf{x}_\Sigma} \mathbf{n}_\Gamma^{\left( d_f \right)} \right\|_2^{-1} \,\mathrm{d}\Sigma + \int_{\partial\Sigma} B \left\| \left( \mathbf{n}_\Sigma \times \nabla_\Sigma d_f \right) \times \left( \mathbf{n}_\Sigma - \nabla_\Sigma d_f\right) \right\|_2 \\
  & \left\| \left( \partial \mathbf{x}_\Gamma \over \partial \mathbf{x}_\Sigma \right)^{-1} \left[ \left( \mathbf{n}_\Sigma \times \nabla_\Sigma d_f \right) \times \left( \mathbf{n}_\Sigma - \nabla_\Sigma d_f\right)\right] \right\|_2^{-1} \mathrm{d}l_{\partial\Sigma} + \int_\Sigma \Bigg[ \rho \left( \mathbf{u} \cdot \nabla_\Gamma^{\left( d_f \right)} \right) \mathbf{u} \cdot \mathbf{u}_a \\
  & + {\eta\over2} \left( \nabla_\Gamma^{\left( d_f \right)} \mathbf{u} + \nabla_\Gamma^{\left( d_f \right)} \mathbf{u}^\mathrm{T} \right) : \left( \nabla_\Gamma^{\left( d_f \right)} \mathbf{u}_a + \nabla_\Gamma^{\left( d_f \right)} \mathbf{u}_a^\mathrm{T} \right) - p \, \mathrm{div}_\Gamma^{\left( d_f \right)} \mathbf{u}_a  - p_a \mathrm{div}_\Gamma^{\left( d_f \right)} \mathbf{u} \\
  & + \alpha \mathbf{u} \cdot \mathbf{u}_a + \left( \lambda \mathbf{u}_a + \lambda_a \mathbf{u} \right) \cdot {\nabla_\Sigma d_f + \mathbf{n}_\Sigma \over \left\| \nabla_\Sigma d_f + \mathbf{n}_\Sigma \right\|_2} \Bigg] \left| {\partial \mathbf{x}_\Gamma \over \partial \mathbf{x}_\Sigma} \right| \left\| {\partial \mathbf{x}_\Gamma \over \partial \mathbf{x}_\Sigma} \mathbf{n}_\Gamma^{\left( d_f \right)} \right\|_2^{-1} \,\mathrm{d}\Sigma \\
  & + \int_\Sigma \Big( r_f^2 \nabla_\Gamma^{\left( d_f \right)} \gamma_f \cdot \nabla_\Gamma^{\left( d_f \right)} \gamma_{fa} + \gamma_f \gamma_{fa} - \gamma \gamma_{fa} \Big) \left| {\partial \mathbf{x}_\Gamma \over \partial \mathbf{x}_\Sigma} \right| \left\| {\partial \mathbf{x}_\Gamma \over \partial \mathbf{x}_\Sigma} \mathbf{n}_\Gamma^{\left( d_f \right)} \right\|_2^{-1} \,\mathrm{d}\Sigma \\
  & + \int_\Sigma r_m^2 \nabla_\Sigma d_f \cdot \nabla_\Sigma d_{fa} + d_f d_{fa} - A_d \left( d_m - {1\over2} \right) d_{fa} \,\mathrm{d}\Sigma
\end{split}
\end{equation}
with
\begin{equation}\label{equ:ConstraintForAugmentedLagrangian}
  \left\{\begin{split}
  & \mathbf{u}_a \in\left(\mathcal{H}\left(\Sigma\right)\right)^3~\mathrm{with}~ \mathbf{u}_a = \mathbf{0}, ~ {\forall \mathbf{x} \in l_{v,\Sigma} } \\
  & p_a \in \mathcal{H}\left(\Sigma\right)~\mathrm{with}~ p_a = 0,~ \forall \mathbf{x} \in \mathcal{P}_\Sigma \\
  & \lambda_a \in \mathcal{L}^2\left(\Sigma\right)~\mathrm{with}~ \lambda_a = 0,~ \forall \mathbf{x} \in l_{v,\Sigma}
  \end{split}\right..\\
\end{equation}

Based on the transformed operators in Eqs. \ref{equ:TransformedTangentialOperator} and \ref{equ:TransformedDivergenceOperator} and their first order variationals in Eqs. \ref{equ:FirstOrderVarisForTangentialOperators} and \ref{equ:FirstOrderVarisForDivergenceOperators}, together with the first order variational of the 2-norm of a vector function
\begin{equation}\label{equ:FirstOrderVariOfLengthnGamma}
\begin{split}
  \delta \left( \left\| \mathbf{f} \right\|_2 \right)^2
  = 2 \left\| \mathbf{f} \right\|_2 \delta \left\| \mathbf{f} \right\|_2 = \delta \mathbf{f}^2 = 2 \mathbf{f} \cdot \delta \mathbf{f} \Rightarrow \delta \left\| \mathbf{f} \right\|_2
  = { \mathbf{f} \over \left\| \mathbf{f} \right\|_2 } \cdot \delta  \mathbf{f}
\end{split}
\end{equation}
with $\mathbf{f}$ representing the vector function, the first order variational of the augmented Lagrangian in Eq. \ref{equ:AugmentedLagrangianMatchOptimization} can be derived as
\begin{equation}\label{equ:1stVarAugmentedLagrangianMatchOptimization}
\begin{split}
  \delta \hat{J} = & \int_\Sigma \left[ {\partial A \over \partial \mathbf{u}} \cdot \delta \mathbf{u} + {\partial A \over \partial \nabla_\Gamma^{\left( d_f \right)} \mathbf{u} } : \left( \nabla_\Gamma^{\left( d_f \right)} \delta \mathbf{u} + \nabla_\Gamma^{\left( d_f, \delta d_f \right)} \mathbf{u} \right) + {\partial A \over \partial p } \delta p + {\partial A \over \partial \gamma_p } {\partial \gamma_p \over \partial \gamma_f } \delta \gamma_f \right] \left| {\partial \mathbf{x}_\Gamma \over \partial \mathbf{x}_\Sigma} \right| \left\| {\partial \mathbf{x}_\Gamma \over \partial \mathbf{x}_\Sigma} \mathbf{n}_\Gamma^{\left( d_f \right)} \right\|_2^{-1} \\
  & + A \left( {\partial \left| {\partial \mathbf{x}_\Gamma \over \partial \mathbf{x}_\Sigma} \right| \left\| {\partial \mathbf{x}_\Gamma \over \partial \mathbf{x}_\Sigma} \mathbf{n}_\Gamma^{\left( d_f \right)} \right\|_2^{-1} \over \partial d_f} \delta d_f + {\partial \left| {\partial \mathbf{x}_\Gamma \over \partial \mathbf{x}_\Sigma} \right| \left\| {\partial \mathbf{x}_\Gamma \over \partial \mathbf{x}_\Sigma} \mathbf{n}_\Gamma^{\left( d_f \right)} \right\|_2^{-1} \over \partial \nabla_\Sigma d_f} \cdot \nabla_\Sigma \delta d_f \right) \,\mathrm{d}\Sigma + \int_{\partial\Sigma} \left( {\partial B \over \partial \mathbf{u}} \cdot \delta \mathbf{u} + {\partial B \over \partial p} \delta p \right) \\
  & \left\| \left( \mathbf{n}_\Sigma \times \nabla_\Sigma d_f \right) \times \left( \mathbf{n}_\Sigma - \nabla_\Sigma d_f\right) \right\|_2 \left\| \left( \partial \mathbf{x}_\Gamma \over \partial \mathbf{x}_\Sigma \right)^{-1} \left[ \left( \mathbf{n}_\Sigma \times \nabla_\Sigma d_f \right) \times \left( \mathbf{n}_\Sigma - \nabla_\Sigma d_f\right)\right] \right\|_2^{-1} \,\mathrm{d}l_{\partial\Sigma} \\
  & + \int_{\partial\Sigma} B { \partial \left\| \left( \mathbf{n}_\Sigma \times \nabla_\Sigma d_f \right) \times \left( \mathbf{n}_\Sigma - \nabla_\Sigma d_f\right) \right\|_2 \left\| \left( \partial \mathbf{x}_\Gamma \over \partial \mathbf{x}_\Sigma \right)^{-1} \left[ \left( \mathbf{n}_\Sigma \times \nabla_\Sigma d_f \right) \times \left( \mathbf{n}_\Sigma - \nabla_\Sigma d_f\right)\right] \right\|_2^{-1} \over \partial d_f } \delta d_f \\
  & + B { \partial \left\| \left( \mathbf{n}_\Sigma \times \nabla_\Sigma d_f \right) \times \left( \mathbf{n}_\Sigma - \nabla_\Sigma d_f\right) \right\|_2 \left\| \left( \partial \mathbf{x}_\Gamma \over \partial \mathbf{x}_\Sigma \right)^{-1} \left[ \left( \mathbf{n}_\Sigma \times \nabla_\Sigma d_f \right) \times \left( \mathbf{n}_\Sigma - \nabla_\Sigma d_f\right)\right] \right\|_2^{-1} \over \partial \nabla_\Sigma d_f } \cdot \nabla_\Sigma \delta d_f \,\mathrm{d}l_{\partial\Sigma} \\
  & + \int_\Sigma \Bigg\{ \rho \left( \delta \mathbf{u} \cdot \nabla_\Gamma^{\left( d_f \right)} \right) \mathbf{u} \cdot \mathbf{u}_a + \rho \left( \mathbf{u} \cdot \nabla_\Gamma^{\left( d_f,\delta d_f \right)} \right) \mathbf{u} \cdot \mathbf{u}_a + \rho \left( \mathbf{u} \cdot \nabla_\Gamma^{\left( d_f \right)} \right) \delta \mathbf{u} \cdot \mathbf{u}_a \\
  & + {\eta\over2} \left( \nabla_\Gamma^{\left( d_f \right)} \delta \mathbf{u} + \nabla_\Gamma^{\left( d_f \right)} \delta \mathbf{u}^\mathrm{T} \right) : \left( \nabla_\Gamma^{\left( d_f \right)} \mathbf{u}_a + \nabla_\Gamma^{\left( d_f \right)} \mathbf{u}_a^\mathrm{T} \right) + {\eta\over2} \left( \nabla_\Gamma^{\left( d_f, \delta d_f \right)} \mathbf{u} + \nabla_\Gamma^{\left( d_f, \delta d_f \right)} \mathbf{u}^\mathrm{T} \right) \\
  & : \left( \nabla_\Gamma^{\left( d_f \right)} \mathbf{u}_a + \nabla_\Gamma^{\left( d_f \right)} \mathbf{u}_a^\mathrm{T} \right) + {\eta\over2} \left( \nabla_\Gamma^{\left( d_f \right)} \mathbf{u} + \nabla_\Gamma^{\left( d_f \right)} \mathbf{u}^\mathrm{T} \right) : \left( \nabla_\Gamma^{\left( d_f, \delta d_f \right)} \mathbf{u}_a + \nabla_\Gamma^{\left( d_f, \delta d_f \right)} \mathbf{u}_a^\mathrm{T} \right) \\
  & - \delta p \, \mathrm{div}_\Gamma^{\left( d_f \right)} \mathbf{u}_a - p \, \mathrm{div}_\Gamma^{\left( d_f, \delta d_f \right)} \mathbf{u}_a - p_a \, \mathrm{div}_\Gamma^{\left( d_f, \delta d_f \right)} \mathbf{u} - p_a \, \mathrm{div}_\Gamma^{\left( d_f \right)} \delta \mathbf{u} + {\partial \alpha \over \partial \gamma_p} {\partial \gamma_p \over \partial \gamma_f} \mathbf{u} \cdot \mathbf{u}_a \delta \gamma_f + \alpha \delta \mathbf{u} \cdot \mathbf{u}_a \\
  & + \delta \lambda \mathbf{u}_a \cdot {\nabla_\Sigma d_f + \mathbf{n}_\Sigma \over \left\| \nabla_\Sigma d_f + \mathbf{n}_\Sigma \right\|_2 } + \lambda \mathbf{u}_a \cdot \Bigg({\nabla_\Sigma \delta d_f \over \left\| \nabla_\Sigma d_f + \mathbf{n}_\Sigma \right\|_2 } - {\nabla_\Sigma d_f + \mathbf{n}_\Sigma \over \left( \nabla_\Sigma d_f + \mathbf{n}_\Sigma \right)^2 } { \left( \nabla_\Sigma d_f + \mathbf{n}_\Sigma \right) \cdot \nabla_\Sigma \delta  d_f \over \left\| \nabla_\Sigma d_f + \mathbf{n}_\Sigma \right\|_2 } \Bigg) \\
  & + \lambda_a \delta \mathbf{u} \cdot {\nabla_\Sigma d_f + \mathbf{n}_\Sigma \over \left\| \nabla_\Sigma d_f + \mathbf{n}_\Sigma \right\|_2 } + \lambda_a \mathbf{u} \cdot \left({\nabla_\Sigma \delta d_f \over \left\| \nabla_\Sigma d_f + \mathbf{n}_\Sigma \right\|_2 } - {\nabla_\Sigma d_f + \mathbf{n}_\Sigma \over \left( \nabla_\Sigma d_f + \mathbf{n}_\Sigma \right)^2 } { \left( \nabla_\Sigma d_f + \mathbf{n}_\Sigma \right) \cdot \nabla_\Sigma \delta  d_f \over \left\| \nabla_\Sigma d_f + \mathbf{n}_\Sigma \right\|_2 } \right) \Bigg\} \\
  & \left| {\partial \mathbf{x}_\Gamma \over \partial \mathbf{x}_\Sigma} \right| \left\| {\partial \mathbf{x}_\Gamma \over \partial \mathbf{x}_\Sigma} \mathbf{n}_\Gamma^{\left( d_f \right)} \right\|_2^{-1} + \Bigg[ \rho \left( \mathbf{u} \cdot \nabla_\Gamma^{\left( d_f \right)} \right) \mathbf{u} \cdot \mathbf{u}_a + {\eta\over2} \left( \nabla_\Gamma^{\left( d_f \right)} \mathbf{u} + \nabla_\Gamma^{\left( d_f \right)} \mathbf{u}^\mathrm{T} \right) : \left( \nabla_\Gamma^{\left( d_f \right)} \mathbf{u}_a + \nabla_\Gamma^{\left( d_f \right)} \mathbf{u}_a^\mathrm{T} \right) \\
  & - p \, \mathrm{div}_\Gamma^{\left( d_f \right)} \mathbf{u}_a - p_a \, \mathrm{div}_\Gamma^{\left( d_f \right)} \mathbf{u} + \alpha \mathbf{u} \cdot \mathbf{u}_a + \lambda \mathbf{u}_a \cdot {\nabla_\Sigma d_f + \mathbf{n}_\Sigma \over \left\| \nabla_\Sigma d_f + \mathbf{n}_\Sigma \right\|_2} + \lambda_a \mathbf{u} \cdot {\nabla_\Sigma d_f + \mathbf{n}_\Sigma \over \left\| \nabla_\Sigma d_f + \mathbf{n}_\Sigma \right\|_2} \Bigg] \\
  & \Bigg( {\partial \left| {\partial \mathbf{x}_\Gamma \over \partial \mathbf{x}_\Sigma} \right| \left\| {\partial \mathbf{x}_\Gamma \over \partial \mathbf{x}_\Sigma} \mathbf{n}_\Gamma^{\left( d_f \right)} \right\|_2^{-1} \over \partial d_f} \delta d_f + {\partial \left| {\partial \mathbf{x}_\Gamma \over \partial \mathbf{x}_\Sigma} \right| \left\| {\partial \mathbf{x}_\Gamma \over \partial \mathbf{x}_\Sigma} \mathbf{n}_\Gamma^{\left( d_f \right)} \right\|_2^{-1} \over \partial \nabla_\Sigma d_f} \cdot \nabla_\Sigma \delta d_f \Bigg) + \bigg[ r_f^2 \bigg( \nabla_\Gamma^{\left( d_f \right)} \delta \gamma_f \cdot \nabla_\Gamma^{\left( d_f \right)} \gamma_{fa} \\
  & + \nabla_\Gamma^{\left( d_f, \delta d_f \right)} \gamma_f \cdot \nabla_\Gamma^{\left( d_f \right)} \gamma_{fa} + \nabla_\Gamma^{\left( d_f \right)} \gamma_f \cdot \nabla_\Gamma^{\left( d_f, \delta d_f \right)} \gamma_{fa} \bigg) + \delta \gamma_f \gamma_{fa} - \delta \gamma \gamma_{fa} \bigg] \\
  & \left| {\partial \mathbf{x}_\Gamma \over \partial \mathbf{x}_\Sigma} \right| \left\| {\partial \mathbf{x}_\Gamma \over \partial \mathbf{x}_\Sigma} \mathbf{n}_\Gamma^{\left( d_f \right)} \right\|_2^{-1} + \left( r_f^2 \nabla_\Gamma^{\left( d_f \right)} \gamma_f \cdot \nabla_\Gamma^{\left( d_f \right)} \gamma_{fa} + \gamma_f \gamma_{fa} - \gamma \gamma_{fa} \right) \\
  & \left( {\partial \left| {\partial \mathbf{x}_\Gamma \over \partial \mathbf{x}_\Sigma} \right| \left\| {\partial \mathbf{x}_\Gamma \over \partial \mathbf{x}_\Sigma} \mathbf{n}_\Gamma^{\left( d_f \right)} \right\|_2^{-1} \over \partial d_f} \delta d_f + {\partial \left| {\partial \mathbf{x}_\Gamma \over \partial \mathbf{x}_\Sigma} \right| \left\| {\partial \mathbf{x}_\Gamma \over \partial \mathbf{x}_\Sigma} \mathbf{n}_\Gamma^{\left( d_f \right)} \right\|_2^{-1} \over \partial \nabla_\Sigma d_f} \cdot \nabla_\Sigma \delta d_f \right) \\
  & + r_m^2 \nabla_\Sigma \delta d_f \cdot \nabla_\Sigma d_{fa} + \delta d_f d_{fa} - A_d \delta d_m d_{fa} \,\mathrm{d}\Sigma
\end{split}
\end{equation}
with the satisfication of the constraints in Eq. \ref{equ:ConstraintForAugmentedLagrangian}
and
\begin{equation}\label{equ:ConstraintForVariationalAugmentedLagrangian}
  \left\{\begin{split}
  & \delta \mathbf{u} \in\left(\mathcal{H}\left(\Sigma\right)\right)^3~\mathrm{with}~ \delta \mathbf{u} = \mathbf{0}, ~ {\forall \mathbf{x} \in l_{v,\Sigma} } \\
  & \delta p \in \mathcal{H}\left(\Sigma\right)~\mathrm{with}~ \delta p = 0,~ \forall \mathbf{x} \in \mathcal{P}_\Sigma \\
  & \delta \lambda \in \mathcal{L}^2\left(\Sigma\right)~\mathrm{with}~ \delta \lambda = 0,~ \forall \mathbf{x} \in l_{v,\Sigma}
  \end{split}\right.. \\
\end{equation}

According to the Karush-Kuhn-Tucker conditions of the PDE constrained optimization problem \cite{HinzeSpringer2009}, the first order variational of the augmented Lagrangian to the variables $\mathbf{u}$, $p$ and $\lambda$ can be set to be zero as
\begin{equation}\label{equ:WeakAdjEquSNSEqu}
\begin{split}
  & \int_\Sigma \bigg[ {\partial A \over \partial \mathbf{u}} \cdot \delta \mathbf{u} + {\partial A \over \partial \nabla_\Gamma^{\left( d_f \right)} \mathbf{u} } : \nabla_\Gamma^{\left( d_f \right)} \delta \mathbf{u} + {\partial A \over \partial p } \delta p + \rho \left( \delta \mathbf{u} \cdot \nabla_\Gamma^{\left( d_f \right)} \right) \mathbf{u} \cdot \mathbf{u}_a + \rho \left( \mathbf{u} \cdot \nabla_\Gamma^{\left( d_f \right)} \right) \delta \mathbf{u} \cdot \mathbf{u}_a \\
  & + {\eta\over2} \left( \nabla_\Gamma^{\left( d_f \right)} \delta \mathbf{u} + \nabla_\Gamma^{\left( d_f \right)} \delta \mathbf{u}^\mathrm{T} \right) : \left( \nabla_\Gamma^{\left( d_f \right)} \mathbf{u}_a + \nabla_\Gamma^{\left( d_f \right)} \mathbf{u}_a^\mathrm{T} \right) - \delta p \, \mathrm{div}_\Gamma^{\left( d_f \right)} \mathbf{u}_a - p_a \, \mathrm{div}_\Gamma^{\left( d_f \right)} \delta \mathbf{u} \\
  & + \alpha \delta \mathbf{u} \cdot \mathbf{u}_a + \left( \delta \lambda \mathbf{u}_a + \lambda_a \delta \mathbf{u} \right) \cdot {\nabla_\Sigma d_f + \mathbf{n}_\Sigma \over \left\| \nabla_\Sigma d_f + \mathbf{n}_\Sigma \right\|_2 } \bigg] \left| {\partial \mathbf{x}_\Gamma \over \partial \mathbf{x}_\Sigma} \right| \left\| {\partial \mathbf{x}_\Gamma \over \partial \mathbf{x}_\Sigma} \mathbf{n}_\Gamma^{\left( d_f \right)} \right\|_2^{-1} \\
  & \,\mathrm{d}\Sigma + \int_{\partial\Sigma} \left( {\partial B \over \partial \mathbf{u}} \cdot \delta \mathbf{u} + {\partial B \over \partial p} \delta p \right) \left\| \left( \mathbf{n}_\Sigma \times \nabla_\Sigma d_f \right) \times \left( \mathbf{n}_\Sigma - \nabla_\Sigma d_f\right) \right\|_2 \\
  & \left\| \left( \partial \mathbf{x}_\Gamma \over \partial \mathbf{x}_\Sigma \right)^{-1} \left[ \left( \mathbf{n}_\Sigma \times \nabla_\Sigma d_f \right) \times \left( \mathbf{n}_\Sigma - \nabla_\Sigma d_f\right)\right] \right\|_2^{-1} \,\mathrm{d}l_{\partial\Sigma} = 0;
\end{split}
\end{equation}
the first order variational of the augmented Lagrangian to the variable $\gamma_f$ can be set to be zero as
\begin{equation}\label{equ:WeakAdjEquSPDEFilterGa}
\begin{split}
  & \int_\Sigma \Bigg( {\partial A \over \partial \gamma_p } {\partial \gamma_p \over \partial \gamma_f } \delta \gamma_f + {\partial \alpha \over \partial \gamma_p} {\partial \gamma_p \over \partial \gamma_f} \mathbf{u} \cdot \mathbf{u}_a \delta \gamma_f + r_f^2 \nabla_\Gamma^{\left( d_f \right)} \delta \gamma_f \cdot \nabla_\Gamma^{\left( d_f \right)} \gamma_{fa} + \delta \gamma_f \gamma_{fa} \Bigg) \\
  & \left| {\partial \mathbf{x}_\Gamma \over \partial \mathbf{x}_\Sigma} \right| \left\| {\partial \mathbf{x}_\Gamma \over \partial \mathbf{x}_\Sigma} \mathbf{n}_\Gamma^{\left( d_f \right)} \right\|_2^{-1} \,\mathrm{d}\Sigma = 0;
\end{split}
\end{equation}
and the first order variational of the augmented Lagrangian to the variable $d_f$ can be set to be zero as
\begin{equation}\label{equ:WeakAdjEquSPDEFilterDm}
\begin{split}
  & \int_\Sigma \Bigg[ {\partial A \over \partial \nabla_\Gamma^{\left( d_f \right)} \mathbf{u} } : \nabla_\Gamma^{\left( d_f, \delta d_f \right)} \mathbf{u} + \rho \left( \mathbf{u} \cdot \nabla_\Gamma^{\left( d_f,\delta d_f \right)} \right) \mathbf{u} \cdot \mathbf{u}_a + {\eta\over2} \left( \nabla_\Gamma^{\left( d_f, \delta d_f \right)} \mathbf{u} + \nabla_\Gamma^{\left( d_f, \delta d_f \right)} \mathbf{u}^\mathrm{T} \right) \\
  & : \left( \nabla_\Gamma^{\left( d_f \right)} \mathbf{u}_a + \nabla_\Gamma^{\left( d_f \right)} \mathbf{u}_a^\mathrm{T} \right) + {\eta\over2} \left( \nabla_\Gamma^{\left( d_f \right)} \mathbf{u} + \nabla_\Gamma^{\left( d_f \right)} \mathbf{u}^\mathrm{T} \right) : \left( \nabla_\Gamma^{\left( d_f, \delta d_f \right)} \mathbf{u}_a + \nabla_\Gamma^{\left( d_f, \delta d_f \right)} \mathbf{u}_a^\mathrm{T} \right) \\
  & - p \, \mathrm{div}_\Gamma^{\left( d_f, \delta d_f \right)} \mathbf{u}_a - p_a \, \mathrm{div}_\Gamma^{\left( d_f, \delta d_f \right)} \mathbf{u} + \left( \lambda \mathbf{u}_a + \lambda_a \mathbf{u} \right) \cdot \bigg({\nabla_\Sigma \delta d_f \over \left\| \nabla_\Sigma d_f + \mathbf{n}_\Sigma \right\|_2 } - {\nabla_\Sigma d_f + \mathbf{n}_\Sigma \over \left( \nabla_\Sigma d_f + \mathbf{n}_\Sigma \right)^2 } \\
  & { \left( \nabla_\Sigma d_f + \mathbf{n}_\Sigma \right) \cdot \nabla_\Sigma \delta  d_f \over \left\| \nabla_\Sigma d_f + \mathbf{n}_\Sigma \right\|_2 } \bigg) + r_f^2 \left( \nabla_\Gamma^{\left( d_f, \delta d_f \right)} \gamma_f \cdot \nabla_\Gamma^{\left( d_f \right)} \gamma_{fa} + \nabla_\Gamma^{\left( d_f \right)} \gamma_f \cdot \nabla_\Gamma^{\left( d_f, \delta d_f \right)} \gamma_{fa} \right) \\
  & + \left( \nabla_\Gamma^{\left(d_f, \delta d_f \right)} d_{\boldsymbol\tau_\Gamma} \cdot \nabla_\Gamma^{\left( d_f \right)} d_{\boldsymbol\tau_\Gamma a} + \nabla_\Gamma^{\left( d_f \right)} d_{\boldsymbol\tau_\Gamma} \cdot \nabla_\Gamma^{\left( d_f, \delta d_f \right)} d_{\boldsymbol\tau_\Gamma a} \right) - \bigg( \nabla_\Gamma^{\left( d_f, \delta d_f \right)} f \cdot \nabla_\Gamma^{\left( d_f \right)} f_a \\
  & + \nabla_\Gamma^{\left( d_f \right)} f \cdot \nabla_\Gamma^{\left( d_f, \delta d_f \right)} f_a \bigg) \Bigg] \left| {\partial \mathbf{x}_\Gamma \over \partial \mathbf{x}_\Sigma} \right| \left\| {\partial \mathbf{x}_\Gamma \over \partial \mathbf{x}_\Sigma} \mathbf{n}_\Gamma^{\left( d_f \right)} \right\|_2^{-1} + \Bigg[ A + \rho \left( \mathbf{u} \cdot \nabla_\Gamma^{\left( d_f \right)} \right) \mathbf{u} \cdot \mathbf{u}_a \\
  & + {\eta\over2} \left( \nabla_\Gamma^{\left( d_f \right)} \mathbf{u} + \nabla_\Gamma^{\left( d_f \right)} \mathbf{u}^\mathrm{T} \right) : \left( \nabla_\Gamma^{\left( d_f \right)} \mathbf{u}_a + \nabla_\Gamma^{\left( d_f \right)} \mathbf{u}_a^\mathrm{T} \right) - p \, \mathrm{div}_\Gamma^{\left( d_f \right)} \mathbf{u}_a - p_a \mathrm{div}_\Gamma^{\left( d_f \right)} \mathbf{u} \\
  & + \alpha \mathbf{u} \cdot \mathbf{u}_a + \left( \lambda \mathbf{u}_a + \lambda_a \mathbf{u} \right) \cdot {\nabla_\Sigma d_f + \mathbf{n}_\Sigma \over \left\| \nabla_\Sigma d_f + \mathbf{n}_\Sigma \right\|_2} + \left( r_f^2 \nabla_\Gamma^{\left( d_f \right)} \gamma_f \cdot \nabla_\Gamma^{\left( d_f \right)} \gamma_{fa} + \gamma_f \gamma_{fa} - \gamma \gamma_{fa} \right) \\
  & + \Big( \nabla_\Gamma^{\left( d_f \right)} d_{\boldsymbol\tau_\Gamma} \cdot \nabla_\Gamma^{\left( d_f \right)} d_{\boldsymbol\tau_\Gamma a} + d_{\boldsymbol\tau_\Gamma} d_{\boldsymbol\tau_\Gamma a} - f d_{\boldsymbol\tau_\Gamma a} - \nabla_\Gamma^{\left( d_f \right)} f \cdot \nabla_\Gamma^{\left( d_f \right)} f_a - f_a \Big) \Bigg] \\
  & \left( {\partial \left| {\partial \mathbf{x}_\Gamma \over \partial \mathbf{x}_\Sigma} \right| \left\| {\partial \mathbf{x}_\Gamma \over \partial \mathbf{x}_\Sigma} \mathbf{n}_\Gamma^{\left( d_f \right)} \right\|_2^{-1} \over \partial d_f} \delta d_f + {\partial \left| {\partial \mathbf{x}_\Gamma \over \partial \mathbf{x}_\Sigma} \right| \left\| {\partial \mathbf{x}_\Gamma \over \partial \mathbf{x}_\Sigma} \mathbf{n}_\Gamma^{\left( d_f \right)} \right\|_2^{-1} \over \partial \nabla_\Sigma d_f} \cdot \nabla_\Sigma \delta d_f \right) + r_m^2 \nabla_\Sigma \delta d_f \cdot \nabla_\Sigma d_{fa} + \delta d_f d_{fa} \,\mathrm{d}\Sigma \\
  & + \int_{\partial\Sigma} B { \partial \left\| \left( \mathbf{n}_\Sigma \times \nabla_\Sigma d_f \right) \times \left( \mathbf{n}_\Sigma - \nabla_\Sigma d_f\right) \right\|_2 \left\| \left( \partial \mathbf{x}_\Gamma \over \partial \mathbf{x}_\Sigma \right)^{-1} \left[ \left( \mathbf{n}_\Sigma \times \nabla_\Sigma d_f \right) \times \left( \mathbf{n}_\Sigma - \nabla_\Sigma d_f\right)\right] \right\|_2^{-1} \over \partial d_f } \\
  & \delta d_f + B { \partial \left\| \left( \mathbf{n}_\Sigma \times \nabla_\Sigma d_f \right) \times \left( \mathbf{n}_\Sigma - \nabla_\Sigma d_f\right) \right\|_2 \left\| \left( \partial \mathbf{x}_\Gamma \over \partial \mathbf{x}_\Sigma \right)^{-1} \left[ \left( \mathbf{n}_\Sigma \times \nabla_\Sigma d_f \right) \times \left( \mathbf{n}_\Sigma - \nabla_\Sigma d_f\right)\right] \right\|_2^{-1} \over \partial \nabla_\Sigma d_f } \\
  & \cdot \nabla_\Sigma \delta d_f \,\mathrm{d}l_{\partial\Sigma} = 0. \\
\end{split}
\end{equation}
The constraints in Eqs. \ref{equ:ConstraintForAugmentedLagrangian} and \ref{equ:ConstraintForVariationalAugmentedLagrangian} are imposed to Eq. \ref{equ:WeakAdjEquSNSEqu}. Further, the adjoint sensitivity of $J$ is derived from
\begin{equation}\label{equ:AdjSensitivityGaDmVariationalForm}
\begin{split}
\delta J = \int_\Sigma - \gamma_{fa} \delta \gamma \left| {\partial \mathbf{x}_\Gamma \over \partial \mathbf{x}_\Sigma} \right| \left\| {\partial \mathbf{x}_\Gamma \over \partial \mathbf{x}_\Sigma} \mathbf{n}_\Gamma^{\left( d_f \right)} \right\|_2^{-1} - A_d d_{fa} \delta d_m \,\mathrm{d}\Sigma.
\end{split}
\end{equation}

Without losing the arbitrariness of $\delta \mathbf{u}$, $\delta p$, $\delta \lambda$, $\delta \gamma_f$, $\delta d_f$, $\delta \gamma$ and $\delta d_m$, one can set $\delta \mathbf{u} = \tilde{\mathbf{u}}_a$ with $\forall \tilde{\mathbf{u}}_a \in \left(\mathcal{H}\left(\Sigma\right)\right)^3$, $\delta p = \tilde{p}_a$ with $\forall \tilde{p}_a \in \mathcal{H}\left(\Sigma\right)$, $\delta \lambda = \tilde{\lambda}_a$ with $\forall \tilde{\lambda}_a \in \mathcal{L}^2\left(\Sigma\right)$, $\delta \gamma_f = \tilde{\gamma}_{fa}$ with $\forall \tilde{\gamma}_{fa} \in \mathcal{H}\left(\Sigma\right)$, $\delta d_f = \tilde{d}_{fa}$ with $\forall \tilde{d}_{fa} \in \mathcal{H}\left(\Sigma\right)$, $\delta \gamma = \tilde{\gamma}$ with $\forall \tilde{\gamma} \in \mathcal{L}^2\left(\Sigma\right)$ and $\delta d_m = \tilde{d}_m$ with $\forall \tilde{d}_m \in \mathcal{L}^2\left(\Sigma\right)$, to derive the adjoint system composed of Eqs. \ref{equ:AdjSensitivityGaDm}, \ref{equ:AdjSurfaceNavierStokesEqusJObjective}, \ref{equ:AdjPDEFilterJObjectiveGa} and \ref{equ:AdjPDEFilterJObjectiveDm}.

\subsection{Adjoint analysis for area constraint} \label{sec:AdjointAnalysisAreaConstraint}

Based on the variational formulations of the surface-PDE filters in Eqs. \ref{equ:VariationalFormulationPDEFilterBaseManifold} and \ref{equ:VariationalFormulationPDEFilter},  the augmented Lagrangian of the pattern area $s\left|\Gamma\right|$ can be derived as
\begin{equation}\label{equ:AugmentedLagrangianAreaConstr1}
\begin{split}
  \widehat{s\left|\Gamma\right|} = & \int_\Sigma \left( \gamma_p + r_f^2 \nabla_\Gamma^{\left( d_f \right)} \gamma_f \cdot \nabla_\Gamma^{\left( d_f \right)} \gamma_{fa} + \gamma_f \gamma_{fa} - \gamma \gamma_{fa} \right) \left| {\partial \mathbf{x}_\Gamma \over \partial \mathbf{x}_\Sigma} \right| \left\| {\partial \mathbf{x}_\Gamma \over \partial \mathbf{x}_\Sigma} \mathbf{n}_\Gamma^{\left( d_f \right)} \right\|_2^{-1} \\
  & + r_m^2 \nabla_\Sigma d_f \cdot \nabla_\Sigma d_{fa} + d_f d_{fa} - A_d \left( d_m - {1\over2} \right) d_{fa} \,\mathrm{d}\Sigma.
\end{split}
\end{equation}
Based on the transformed operators in Eq. \ref{equ:TransformedTangentialOperator} and its first order variational in Eq. \ref{equ:FirstOrderVarisForTangentialOperators}, the first order variational of $\widehat{s\left|\Gamma\right|}$ can be derived as
\begin{equation}\label{equ:1stVarAugmentedLagrangianAreaConstr1}
\begin{split}
  \delta \widehat{s\left|\Gamma\right|} = & \int_\Sigma \bigg[ {\partial \gamma_p \over \partial \gamma_f} \delta \gamma_f + r_f^2 \Big( \nabla_\Gamma^{\left( d_f \right)} \delta \gamma_f \cdot \nabla_\Gamma^{\left( d_f \right)} \gamma_{fa} + \nabla_\Gamma^{\left( d_f, \delta d_f \right)} \gamma_f \cdot \nabla_\Gamma^{\left( d_f \right)} \gamma_{fa} + \nabla_\Gamma^{\left( d_f \right)} \gamma_f \cdot \nabla_\Gamma^{\left( d_f, \delta d_f \right)} \gamma_{fa} \Big) \\
  & + \delta \gamma_f \gamma_{fa} - \delta \gamma \gamma_{fa} \bigg] \left| {\partial \mathbf{x}_\Gamma \over \partial \mathbf{x}_\Sigma} \right| \left\| {\partial \mathbf{x}_\Gamma \over \partial \mathbf{x}_\Sigma} \mathbf{n}_\Gamma^{\left( d_f \right)} \right\|_2^{-1} + \left( \gamma_p + r_f^2 \nabla_\Gamma^{\left( d_f \right)} \gamma_f \cdot \nabla_\Gamma^{\left( d_f \right)} \gamma_{fa} + \gamma_f \gamma_{fa} - \gamma \gamma_{fa} \right) \\
  & \left( { \partial \left| {\partial \mathbf{x}_\Gamma \over \partial \mathbf{x}_\Sigma} \right| \left\| {\partial \mathbf{x}_\Gamma \over \partial \mathbf{x}_\Sigma} \mathbf{n}_\Gamma^{\left( d_f \right)} \right\|_2^{-1} \over \partial d_f} \delta d_f + { \partial \left| {\partial \mathbf{x}_\Gamma \over \partial \mathbf{x}_\Sigma} \right| \left\| {\partial \mathbf{x}_\Gamma \over \partial \mathbf{x}_\Sigma} \mathbf{n}_\Gamma^{\left( d_f \right)} \right\|_2^{-1} \over \partial \nabla_\Sigma d_f} \cdot \nabla_\Sigma \delta d_f \right) \\
  & + r_m^2 \nabla_\Sigma \delta d_f \cdot \nabla_\Sigma d_{fa} + \delta d_f d_{fa} - A_d \delta d_m d_{fa} \,\mathrm{d}\Sigma.
\end{split}
\end{equation}
According to the Karush-Kuhn-Tucker conditions of the PDE constrained optimization problem \cite{HinzeSpringer2009}, the first order variational of the augmented Lagrangian to the variable $\gamma_f$ can be set to be zero as
\begin{equation}\label{equ:AdjEquAreaConstr1GaInform}
\begin{split}
  & \int_\Sigma \left( {\partial \gamma_p \over \partial \gamma_f} \delta \gamma_f + r_f^2 \nabla_\Gamma^{\left( d_f \right)} \delta \gamma_f \cdot \nabla_\Gamma^{\left( d_f \right)} \gamma_{fa} + \delta \gamma_f \gamma_{fa} \right) \left| {\partial \mathbf{x}_\Gamma \over \partial \mathbf{x}_\Sigma} \right| \left\| {\partial \mathbf{x}_\Gamma \over \partial \mathbf{x}_\Sigma} \mathbf{n}_\Gamma^{\left( d_f \right)} \right\|_2^{-1} \,\mathrm{d}\Sigma = 0;
\end{split}
\end{equation}
and the first order variational of the augmented Lagrangian to the variable $d_f$ can be set to be zero as
\begin{equation}\label{equ:AdjEquAreaConstr1DmInform}
\begin{split}
  & \int_\Sigma r_f^2 \left( \nabla_\Gamma^{\left( d_f, \delta d_f \right)} \gamma_f \cdot \nabla_\Gamma^{\left( d_f \right)} \gamma_{fa} + \nabla_\Gamma^{\left( d_f \right)} \gamma_f \cdot \nabla_\Gamma^{\left( d_f, \delta d_f \right)} \gamma_{fa} \right) \left| {\partial \mathbf{x}_\Gamma \over \partial \mathbf{x}_\Sigma} \right| \left\| {\partial \mathbf{x}_\Gamma \over \partial \mathbf{x}_\Sigma} \mathbf{n}_\Gamma^{\left( d_f \right)} \right\|_2^{-1} \\
  & + \left( \gamma_p + r_f^2 \nabla_\Gamma^{\left( d_f \right)} \gamma_f \cdot \nabla_\Gamma^{\left( d_f \right)} \gamma_{fa} + \gamma_f \gamma_{fa} - \gamma \gamma_{fa} \right) \\
  & \left( { \partial \left| {\partial \mathbf{x}_\Gamma \over \partial \mathbf{x}_\Sigma} \right| \left\| {\partial \mathbf{x}_\Gamma \over \partial \mathbf{x}_\Sigma} \mathbf{n}_\Gamma^{\left( d_f \right)} \right\|_2^{-1} \over \partial d_f} \delta d_f + { \partial \left| {\partial \mathbf{x}_\Gamma \over \partial \mathbf{x}_\Sigma} \right| \left\| {\partial \mathbf{x}_\Gamma \over \partial \mathbf{x}_\Sigma} \mathbf{n}_\Gamma^{\left( d_f \right)} \right\|_2^{-1} \over \partial \nabla_\Sigma d_f} \cdot \nabla_\Sigma \delta d_f \right) \\
  & + r_m^2 \nabla_\Sigma \delta d_f \cdot \nabla_\Sigma d_{fa} + \delta d_f d_{fa} \,\mathrm{d}\Sigma = 0.
\end{split}
\end{equation}
Further, the adjoint sensitivity of $s\left|\Gamma\right|$ is derived from
\begin{equation}\label{equ:AdjSensAreaConstr1Inform}
\begin{split}
  \delta \widehat{s\left|\Gamma\right|} = & \int_\Sigma - \delta \gamma \gamma_{fa} \left| {\partial \mathbf{x}_\Gamma \over \partial \mathbf{x}_\Sigma} \right| \left\| {\partial \mathbf{x}_\Gamma \over \partial \mathbf{x}_\Sigma} \mathbf{n}_\Gamma^{\left( d_f \right)} \right\|_2^{-1} - A_d \delta d_m d_{fa} \,\mathrm{d}\Sigma.
\end{split}
\end{equation}

Based on the variational formulations of the surface-PDE filters in Eqs. \ref{equ:VariationalFormulationPDEFilterBaseManifold} and \ref{equ:VariationalFormulationPDEFilter},  the augmented Lagrangian of the area of the implicit 2-manifold can be derived as
\begin{equation}\label{equ:AugmentedLagrangianAreaConstr11}
\begin{split}
  \widehat{\left|\Gamma\right|} = & \int_\Sigma \left( 1 + r_f^2 \nabla_\Gamma^{\left(d_f\right)} \gamma_f \cdot \nabla_\Gamma^{\left(d_f\right)} \gamma_{fa} + \gamma_f \gamma_{fa} - \gamma \gamma_{fa} \right) \left| {\partial \mathbf{x}_\Gamma \over \partial \mathbf{x}_\Sigma} \right| \left\| {\partial \mathbf{x}_\Gamma \over \partial \mathbf{x}_\Sigma} \mathbf{n}_\Gamma^{\left( d_f \right)} \right\|_2^{-1} \\
  & + r_m^2 \nabla_\Sigma d_f \cdot \nabla_\Sigma d_{fa} + d_f d_{fa} - A_d \left( d_m - {1\over2} \right) d_{fa} \,\mathrm{d}\Sigma.
\end{split}
\end{equation}
Based on the transformed operators in Eq. \ref{equ:TransformedTangentialOperator} and its first order variational in Eq. \ref{equ:FirstOrderVarisForTangentialOperators}, the first order variational of $\widehat{\left|\Gamma\right|}$ can be derived as
\begin{equation}\label{equ:1stVarAugmentedLagrangianAreaConstr1}
\begin{split}
  \delta \widehat{\left|\Gamma\right|} =
  & \int_\Sigma \bigg[ r_f^2 \left( \nabla_\Gamma^{\left(d_f\right)} \delta \gamma_f \cdot \nabla_\Gamma^{\left(d_f\right)} \gamma_{fa} + \nabla_\Gamma^{\left(d_f, \delta d_f\right)} \gamma_f \cdot \nabla_\Gamma^{\left(d_f\right)} \gamma_{fa} + \nabla_\Gamma^{\left(d_f\right)} \gamma_f \cdot \nabla_\Gamma^{\left(d_f, \delta d_f\right)} \gamma_{fa} \right) \\
  & + \delta \gamma_f \gamma_{fa} - \delta \gamma \gamma_{fa} \bigg] \left| {\partial \mathbf{x}_\Gamma \over \partial \mathbf{x}_\Sigma} \right| \left\| {\partial \mathbf{x}_\Gamma \over \partial \mathbf{x}_\Sigma} \mathbf{n}_\Gamma^{\left( d_f \right)} \right\|_2^{-1} + \left( 1 + r_f^2 \nabla_\Gamma^{\left(d_f\right)} \gamma_f \cdot \nabla_\Gamma^{\left(d_f\right)} \gamma_{fa} + \gamma_f \gamma_{fa} - \gamma \gamma_{fa} \right) \\
  & \left( { \partial \left| {\partial \mathbf{x}_\Gamma \over \partial \mathbf{x}_\Sigma} \right| \left\| {\partial \mathbf{x}_\Gamma \over \partial \mathbf{x}_\Sigma} \mathbf{n}_\Gamma^{\left( d_f \right)} \right\|_2^{-1} \over \partial d_f} \delta d_f + { \partial \left| {\partial \mathbf{x}_\Gamma \over \partial \mathbf{x}_\Sigma} \right| \left\| {\partial \mathbf{x}_\Gamma \over \partial \mathbf{x}_\Sigma} \mathbf{n}_\Gamma^{\left( d_f \right)} \right\|_2^{-1} \over \partial \nabla_\Sigma d_f} \cdot \nabla_\Sigma \delta d_f \right) \\
  & + r_m^2 \nabla_\Sigma \delta d_f \cdot \nabla_\Sigma d_{fa} + \delta d_f d_{fa} - A_d \delta d_m d_{fa} \,\mathrm{d}\Sigma.
\end{split}
\end{equation}
According to the Karush-Kuhn-Tucker conditions of the PDE constrained optimization problem \cite{HinzeSpringer2009}, the first order variational of the augmented Lagrangian to the variable $\gamma_f$ can be set to be zero as
\begin{equation}\label{equ:AdjEquAreaConstr11GaInform}
\begin{split}
  & \int_\Sigma \left( r_f^2 \nabla_\Gamma^{\left(d_f\right)} \delta \gamma_f \cdot \nabla_\Gamma^{\left(d_f\right)} \gamma_{fa} + \delta \gamma_f \gamma_{fa} \right) \left| {\partial \mathbf{x}_\Gamma \over \partial \mathbf{x}_\Sigma} \right| \left\| {\partial \mathbf{x}_\Gamma \over \partial \mathbf{x}_\Sigma} \mathbf{n}_\Gamma^{\left( d_f \right)} \right\|_2^{-1} \,\mathrm{d}\Sigma = 0;
\end{split}
\end{equation}
and the first order variational of the augmented Lagrangian to the variable $d_f$ can be set to be zero as
\begin{equation}\label{equ:AdjEquAreaConstr11DmInform}
\begin{split}
  & \int_\Sigma r_f^2 \left( \nabla_\Gamma^{\left(d_f, \delta d_f\right)} \gamma_f \cdot \nabla_\Gamma^{\left(d_f\right)} \gamma_{fa} + \nabla_\Gamma^{\left(d_f\right)} \gamma_f \cdot \nabla_\Gamma^{\left(d_f, \delta d_f\right)} \gamma_{fa} \right) \left| {\partial \mathbf{x}_\Gamma \over \partial \mathbf{x}_\Sigma} \right| \left\| {\partial \mathbf{x}_\Gamma \over \partial \mathbf{x}_\Sigma} \mathbf{n}_\Gamma^{\left( d_f \right)} \right\|_2^{-1} \\
  & + \left( 1 + r_f^2 \nabla_\Gamma^{\left(d_f\right)} \gamma_f \cdot \nabla_\Gamma^{\left(d_f\right)} \gamma_{fa} + \gamma_f \gamma_{fa} - \gamma \gamma_{fa} \right) \\
  & \left( { \partial \left| {\partial \mathbf{x}_\Gamma \over \partial \mathbf{x}_\Sigma} \right| \left\| {\partial \mathbf{x}_\Gamma \over \partial \mathbf{x}_\Sigma} \mathbf{n}_\Gamma^{\left( d_f \right)} \right\|_2^{-1} \over \partial d_f} \delta d_f + { \partial \left| {\partial \mathbf{x}_\Gamma \over \partial \mathbf{x}_\Sigma} \right| \left\| {\partial \mathbf{x}_\Gamma \over \partial \mathbf{x}_\Sigma} \mathbf{n}_\Gamma^{\left( d_f \right)} \right\|_2^{-1} \over \partial \nabla_\Sigma d_f} \cdot \nabla_\Sigma \delta d_f \right) \\
  & + r_m^2 \nabla_\Sigma \delta d_f \cdot \nabla_\Sigma d_{fa} + \delta d_f d_{fa} \,\mathrm{d}\Sigma = 0.
\end{split}
\end{equation}
Further, the adjoint sensitivity of $\left|\Gamma\right|$ is derived from
\begin{equation}\label{equ:AdjSensVolConstr1Inform}
\begin{split}
  \delta \widehat{\left|\Gamma\right|} = & \int_\Sigma - \delta \gamma \gamma_{fa} \left| {\partial \mathbf{x}_\Gamma \over \partial \mathbf{x}_\Sigma} \right| \left\| {\partial \mathbf{x}_\Gamma \over \partial \mathbf{x}_\Sigma} \mathbf{n}_\Gamma^{\left( d_f \right)} \right\|_2^{-1} - A_d \delta d_m d_{fa} \,\mathrm{d}\Sigma.
\end{split}
\end{equation}

Without losing the arbitrariness of $\delta \gamma_f$, $\delta d_f$, $\delta \gamma$ and $\delta d_m$, one can set $\delta \gamma_f = \tilde{\gamma}_{fa}$ with $\forall \tilde{\gamma}_{fa} \in \mathcal{H}\left(\Sigma\right)$, $\delta d_f = \tilde{d}_{fa}$ with $\forall \tilde{d}_{fa} \in \mathcal{H}\left(\Sigma\right)$, $\delta \gamma = \tilde{\gamma}$ with $\forall \tilde{\gamma} \in \mathcal{L}^2\left(\Sigma\right)$ and $\delta d_m = \tilde{d}_m$ with $\forall \tilde{d}_m \in \mathcal{L}^2\left(\Sigma\right)$, to derive the adjoint systems composed of Eqs. \ref{equ:AdjSensAreaConstr1}, \ref{equ:AdjEquAreaConstr1Ga}, \ref{equ:AdjEquAreaConstr1Dm} and Eqs. \ref{equ:AdjSensAreaConstr11}, \ref{equ:AdjEquAreaConstr11Ga}, \ref{equ:AdjEquAreaConstr11Dm}, respectively. Then, the adjoint sensitivity of the area fraction $s$ can be derived from Eq. \ref{equ:AdjSensitivityGaDmAreaConstr}.

\subsection{Adjoint analysis for volume constraint} \label{sec:AdjointAnalysisVolumeConstraint}

Based on the variational formulations of the surface-PDE filter in Eq. \ref{equ:VariationalFormulationPDEFilterBaseManifold},  the augmented Lagrangian of the volume fraction $v$ can be derived as
\begin{equation}\label{equ:AugmentedLagrangianVolConstr1}
\begin{split}
  \hat{v} = & \int_\Sigma {1 \over \left| \Sigma \right|} d_f + r_m^2 \nabla_\Sigma d_f \cdot \nabla_\Sigma d_{fa} + d_f d_{fa} - A_d \left( d_m - {1\over2} \right) d_{fa} \,\mathrm{d}\Sigma.
\end{split}
\end{equation}
The first order variational of $\hat{v}$ can be derived as
\begin{equation}\label{equ:1stVarAugmentedLagrangianVolConstr1}
\begin{split}
  \delta \hat{v} = & \int_\Sigma {1 \over \left| \Sigma \right|} \delta d_f + r_m^2 \nabla_\Sigma \delta d_f \cdot \nabla_\Sigma d_{fa} + \delta d_f d_{fa} - A_d \delta d_m d_{fa} \,\mathrm{d}\Sigma.
\end{split}
\end{equation}
According to the Karush-Kuhn-Tucker conditions of the PDE constrained optimization problem \cite{HinzeSpringer2009}, the first order variational of the augmented Lagrangian to the variable $d_f$ can be set to be zero as
\begin{equation}\label{equ:AdjEquVolConstr1DmInform}
\begin{split}
  & \int_\Sigma {1 \over \left| \Sigma \right|} \delta d_f + r_m^2 \nabla_\Sigma \delta d_f \cdot \nabla_\Sigma d_{fa} + \delta d_f d_{fa} \,\mathrm{d}\Sigma = 0.
\end{split}
\end{equation}
Further, the adjoint sensitivity of $v$ is derived from
\begin{equation}\label{equ:AdjSensVolConstr1Inform}
\begin{split}
  \delta \hat{v} = & \int_\Sigma - A_d \delta d_m d_{fa} \,\mathrm{d}\Sigma.
\end{split}
\end{equation}

Without losing the arbitrariness of $\delta d_f$ and $\delta d_m$, one can set $\delta d_f = \tilde{d}_{fa}$ with $\forall \tilde{\gamma}_{fa} \in \mathcal{H}\left(\Sigma\right)$, and $\delta d_m = \tilde{d}_m$ with $\forall \tilde{d}_m \in \mathcal{L}^2\left(\Sigma\right)$, to derive the adjoint system composed of Eqs. \ref{equ:AdjSensitivityGaDmVolConstr} and \ref{equ:AdjEquVolConstrDm}.

\end{document}